\documentclass[12pt]{article}

\usepackage{amsmath}
\usepackage{amssymb}
\usepackage{latexsym}
\usepackage{epsf}
\usepackage{supertabular}
\usepackage{rotating}
\usepackage{multicol}
\usepackage{multirow}
\usepackage{epsfig}
\usepackage{fancyhdr}       

\newif\ifger
\gerfalse

\newtheorem{theorem}{Theorem}[section]

\newtheorem{lemma}[theorem]{Lemma}

\newtheorem{corollary}[theorem]{Corollary}

\newtheorem{definition}[theorem]{Definition}

\newtheorem{remark}[theorem]{Remark}

\newtheorem{proposition}[theorem]{Proposition}

\newtheorem{algorithm}{Algorithm}

\newcommand{\mm}{\phantom{mm}}

\def\calP{\mathcal{P}}

\def\calC{{\mathfrak C}}

\def\calD{\mathcal{D}}
\newcommand{\lcm}{\mathrm{lcm}}

\newcommand{\cC}{{\mathfrak C}}

\newcommand{\cL}{{\cal L}}
\newcommand{\cM}{{\cal M}}

\newcommand{\cP}{{\cal P}}

\newcommand{\cS}{{\cal S}}

\newcommand{\spec}{{\rm spec}}

\def\PG{{\rm PG}}

\def\Soc{{\rm Soc}}
\def\Aut{{\rm Aut}}

\newcommand{\Fix}{{\rm Fix}}
\newcommand{\Sym}{{\rm Sym}}
\newcommand{\Alt}{{\rm Alt}}
\newcommand{\PSL}{{\rm PSL}}
\newcommand{\PGL}{{\rm PGL}}
\newcommand{\PSSL}{{\rm P\Sigma L}}
\newcommand{\PGGL}{{\rm P\Gamma L}}

\newcommand{\GL}{{\rm GL}}

\newcommand{\ASL}{{\rm ASL}}
\newcommand{\AGL}{{\rm AGL}}

\newcommand{\AGGL}{{\rm A\Gamma L}}

\newcommand{\rk}{{\rm rk}}

\def\De{\Delta}

\def\al{\alpha}

\def\lamo{\lambda_0}
\def\lam{\lambda}

\title{Linear spaces with a line-transitive point-imprimitive
automorphism group and Fang-Li parameter gcd({\it k,r}) at most eight
\footnote{This research was supported by the Australian Research Council
grants DP0557587 and DP0209706. The sixth author was partly supported by
the NSF of Guangdong Province.}}

\author{Anton Betten, Anne Delandtsheer, Maska Law\\
Alice C. Niemeyer, Cheryl E. Praeger, Shenglin Zhou}

\begin{document}
\maketitle

\begin{abstract} In 1991, Weidong Fang and Huiling Li proved that there are only finitely many
non-trivial linear spaces that admit a line-transitive, point imprimitive group action, for a given
value of $\gcd(k,r)$, where $k$ is the line size and $r$ is the number of lines on a point.  The aim
of this paper is to make that result effective. We obtain a classification of all linear spaces with
this property having $\gcd(k,r)\leq 8$. To achieve this we collect together existing theory, and
prove additional theoretical restrictions of both a combinatorial and group theoretic nature. These
are organised into a series of algorithms that, for $\gcd(k,r)$ up to a given maximum value, return
a list of candidate parameter values and candidate groups. We examine in detail each of the
possibilities returned by these algorithms for $\gcd(k,r)\leq 8$, and complete the classification in
this case. 



\end{abstract}

\parskip 1ex

\begin{center}
\begin{tabular}{|ll|}\hline
1&Introduction\\
2&Basic facts about line-transitive, point-imprimitive linear spaces\\
3&Outline of the search strategy\\
4&Parameter restrictions and 2-step imprimitivity\\
5&Sifting for primitive top and bottom groups\\
6&Quasiprimitive and $G$-normal sifts\\
7&Results of algorithms applied for $k^{(r)}\leq8$\\
8&Analysing remaining {\sc Lines}\\
A&Appendix: {\sc ParameterList(8)}\\
\hline
\end{tabular}
\end{center}

\section{Introduction}

A {\em finite linear space} $\cS=(\cP,\cL)$ consists of a finite set
$\cP$ of points and a non-empty set $\cL$ of distinguished subsets of
$\cP$ called lines such that any two points lie in exactly one line and
each line contains at least two points. A linear space is said to be
\emph{trivial} if it has only one line, or if all its lines have only
two points; otherwise it is called non-trivial.  The {\em automorphism
group} $\Aut(\cS)$ of $\cS$ consists of all permutations of $\cP$ that
leave $\cL$ invariant. We are interested in {\em line-transitive} linear
spaces $\cS,$ that is, those for which $\Aut(\cS)$ acts transitively on
$\cL$. In this case, the size of lines is constant, say $k$, and to
avoid trivial cases we assume that $2<k<v$. Also, by a result of Richard
Block~\cite{Block67},\cite{Block68}, if a subgroup $G \le \Aut(\cS)$ is
transitive on the lines of $\cS,$ then it is also transitive on points.
It is possible for a line-transitive group $G$ to leave invariant a
non-trivial partition of the point set (see
Subsection~\ref{sub:imprim}), and in this case we say that $G$ is
\emph{point-imprimitive} on $\cS$. 

Let $v,$ $b$ and $r$ be the number of points, the number of lines, and the number of lines through a point, respectively.
In 1991, Weidong Fang and Huiling Li~\cite{FangLi91,FangLi93} proved that for a given value of $k^{(r)}:=\gcd(k,r)$,
there are only finitely many line-transitive, point-imprimitive linear spaces.  Moreover, in~\cite{LiLiu01}, Huiling Li
and Weijun Liu proved that, if $G$ is line-primitive and $k^{(r)}\leq 10$, then $G$ is point-primitive. The aim of this
paper is to demonstrate that the `Fang-Li bound' can be made effective without the additional assumption of
line-primitivity.  We present a collection of tests, based on the theory of point-imprimitive, line-transitive linear
spaces $\cS$, to restrict both the parameters of such spaces, and the structure of a line-transitive, point-imprimitive
subgroup $G$. These tests were organised into a series of algorithms, which we describe in the paper. The algorithms were
then applied to produce a list of all candidate parameters and groups for pairs $(\cS, G)$, in the case where $G$ is
line-transitive and point-imprimitive on $\cS$, and the Fang-Li parameter $k^{(r)}$ is at most eight. We then dealt with
each of these possibilities and classified all line-transitive, point-imprimitive linear spaces with $k^{(r)}\leq8$.

\begin{theorem}\label{main}
Suppose that $\cS$ is a linear space  with
$v$ points, line size $k$, and $r$ lines on each point, that admits a 
line-transitive, point-imprimitive subgroup of automorphisms, and is 
such that the Fang-Li parameter  $k^{(r)}\leq8$. Then $\cS$ is the 
Desarguesian projective plane $\PG(2,4)$ or $\PG(2,7)$, or the 
Mills design or Colbourn-McCalla design, both with 
$(v,k)=(91,6)$, or one of the $467$ designs constructed in {\rm\cite{NNOPP}}, 
all with $(v,k)=(729,8)$.
\end{theorem}

This result may also be viewed as a strengthening of a classification of Camina and 
Mischke~\cite{CaminaMischke96}: we replace their condition $k\leq 8$ with the much weaker 
condition gcd$(k,r)\leq 8$ and prove that no additional examples arise.

\section{Basic facts about line-transitive, point-imprimitive linear spaces}
\label{sect:basic}

The following equalities and inequalities are considered basic, for a linear 
space with parameters $v,r,k$ and $b\geq2$, as defined above.

\begin{align}
vr &= kb,  \label{eqn:vrkb} \\
v-1 &= r(k-1), \label{eqn:vrk} \\
b &= \binom{v}{2} / \binom{k}{2} \label{eqn:b} \\
\intertext{and}
v &\le b
\text{\;\; (Fisher's inequality), \;\;}
\label{eqn:fisher}
\end{align}
with equality if and only if $\cS$ is a {\em projective plane,} that is, 
a space where $v =n^2+n+1=b$ and  $k=r=n+1$ for some integer $n.$

\subsection{Point partitions and the Delandtsheer--Doyen parameters}
\label{sub:imprim}

A partition $\cC$ of a finite set $X$ is a set of non-empty pairwise disjoint subsets whose union
equals the set $X.$ We write $\cC = \{ C_1,\ldots, C_d \}$ and call the $C_i$ the {\em classes }
of $\cC.$ A group $G$ acting on a set $X$ is said to {\em leave the partition $\cC$ invariant} if,
for all $g \in G$ and all $C \in \cC,$ the image $C^g$ also is a class of $\cC.$ A transitive
group $G$ of permutations of $X$ is said to act {\em imprimitively} on $X$ if there exists a
$G$-invariant partition of $X$ which is not trivial, that is to say, the partition $\cC$ neither
consists of only one class nor does it contain only one-element classes. Otherwise, the group is
said to be {\em primitive} on $X.$

\bigskip

Linear spaces admitting a line-transitive, point-imprimitive automorphism group deserve special
attention due to the following result, which shows that the number of points is bounded above by a
function of the line size $k.$

\smallskip

\begin{theorem}[Delandtsheer and Doyen~\cite{DelandtsheerDoyen89}]\label{thm:DD} Let 
$\cS = ( \cP, \cL )$ be a linear space admitting a line-transitive point-imprimitive automorphism
group $G$. Let $\cC = \{C_1, \ldots , C_d\}$ be a $G$-invariant partition of $\cP$ with $d$
classes of size $c$. Let $x$ be the number of inner pairs of a line $\lambda,$ that is, pairs of
points $\{\alpha,\beta\} \subseteq \lambda$ such that $\alpha$ and $\beta$ lie in the same class
of $\cC$. Then there exists another positive integer $y$ such that 
\begin{align}\label{eqn:DD} c =\frac{\binom{k}{2}-x}{y} \quad \text{and} \quad 
d=\frac{\binom{k}{2}-y}{x}.
\end{align}
\end{theorem}

We call the pair $(x,y)$ the {\em Delandtsheer-Doyen parameters} 
corresponding to~$\cC.$
Note that the above equalities are equivalent to
\begin{align}\label{eqn:DDreversed}
x = \frac{\binom{k}{2}(c-1)}{cd-1}
\quad \text{and} \quad
y = \frac{\binom{k}{2}(d-1)}{cd-1}.
\end{align}
Hence the parameters $(c,d,k)$ and $(x,y,k)$ mutually determine each other.
Moreover, if $(c,d,k)$ corresponds to $(x,y,k)$, then the triple
$(d,c,k)$ corresponds to $(y,x,k)$. We emphasise that this second 
correspondence is purely arithmetic: if there is a line-transitive 
point-imprimitive linear space with Delandtsheer--Doyen parameters $(x,y)$ 
corresponding to some point partition, there may or may not exist an example
for which $(y,x)$ are the Delandtsheer--Doyen parameters.

\subsection{The Fang-Li parameters of linear spaces}\label{sub:fangli}

In 1991, Fang and Li introduced some more parameters of linear spaces.
Let $\cS = (\cP,\cL)$ be a linear space admitting a line-transitive
automorphism group $G.$
Note that (\ref{eqn:vrk}) implies that $v$ and $r$ are coprime.
Put
\begin{align*}
k^{(v)} := \gcd(k,v), \quad
k^{(r)} := \gcd(k,r), \quad
b^{(v)} := \gcd(b,v), \quad
b^{(r)} := \gcd(b,r).
\end{align*}
Then the parameters $v, k, b, r$ factorise naturally, as can be seen from 
equations (\ref{eqn:vrkb}) -- (\ref{eqn:b}).

\begin{lemma}\label{prop:fl1} 
\begin{enumerate}
\item 
$k = k^{(v)} \cdot k^{(r)}$ and $b=b^{(v)} \cdot b^{(r)}$ and
\item 
$v = b^{(v)} \cdot k^{(v)}$ and $r = b^{(r)} \cdot k^{(r)}.$
\end{enumerate}
\end{lemma}

Fisher's inequality in terms of the Fang-Li parameters becomes:
\begin{align}\label{fisherfangli}
k^{(v)} \le b^{(r)},
\end{align}
with equality if and only if the linear space is a projective plane, and in 
this case we have in fact that $k^{(v)}=b^{(r)}=1$.

\bigskip

If $G$ preserves a non-trivial partition $\cC$
of the point set $\cP$ with
$d$ classes of size $c,$ much more can be said:

\smallskip

\begin{proposition}[Fang and Li~\cite{FangLi91}]\label{prop:fl2} 
There exist positive integers $\gamma$ and $\delta$ such that
\begin{enumerate}
\item 
\label{flcd}
$c=\gamma b^{(r)} + 1$ and $d = \delta b^{(r)} + 1;$
\item 
$x = \gamma k^{(v)} / 2$ and $y = \delta k^{(v)} / 2$;
In particular, $(\gamma, c-1, x)=\gamma(1,b^{(r)},\frac{k^{(v)}}{2})$ and 
$(\delta, d-1, y)=\delta(1,b^{(r)},\frac{k^{(v)}}{2})$;
\item 
$\gamma + \delta + \gamma\delta b^{(r)} = k^{(r)} (k-1)$ and $\gamma \delta < {k^{(r)}}^2$;
\item 
$k^{(v)}$ divides $(\gamma + k^{(r)})(\delta + k^{(r)})$.
\end{enumerate}
\end{proposition}

\bigskip

The parameters $k^{(v)}, k^{(r)}, b^{(v)}, b^{(r)}, \gamma, \delta$ are called
the {\em Fang-Li parameters} of the linear space corresponding to the partition
$\cC$. Note that parts (iii) and (iv) of Proposition~\ref{prop:fl2} imply that
$k$ is bounded above by a function of $k^{(r)}$, and then Theorem~\ref{thm:DD}
implies that $v$ also is bounded above by a function of $k^{(r)}$.

\subsection{Partition refinements}\label{sub:partitions}

Recall the ordering of (unordered) partitions of sets: A partition $\cC$ of a set $X$ is {\em
contained} in a partition $\cC'$ if every class of $\cC$ is contained in a class of $\cC'.$ In
this case, we say that $\cC$ {\em refines} $\cC'$ and that $\cC'$ is {\em coarser} than $\cC.$ The
set of all partitions of $X$ equipped with this ordering forms a lattice, called the {\em
partition lattice} of $X$. In this paper, partitions left invariant by an imprimitive group $G$
acting on $X$ are of special importance. Note that the join and meet of partitions invariant under
a group $G$ are again $G$-invariant. Hence the set of $G$-invariant partitions forms a sublattice
of the partition lattice of $X$.

We use the following notation: for a group $G$ acting on a set $X$ and a subset $C\subseteq X$, 
we define the setwise and pointwise stabilisers of $C$ in $G$ by $G_C$ and $G_{(C)}$ respectively. 
Note that $G_C$ acts on $C$ with kernel $G_{(C)}$, and we use $G^C$ to denote the permutation 
group on $C$ induced by $G_C$. Clearly, $G^C\cong G_C/G_{(C)}$.

A transitive group $G$ of permutations of a set $X$ leaving invariant a non-trivial partition
$\cC$ of $X$ gives rise to two further permutation groups. For a class $C \in \cC,$ we have the
group $G^C$ induced on $C$ (this group is sometimes called the {\em bottom group} and is
independent of the choice of $C$ up to permutational isomorphism). The permutation group which $G$
induces on the set of classes is denoted by $G^\cC$ (this group is sometimes called the {\em top
group}).  By a standard result about permutation groups, $G$ can be embedded in the wreath product
$G^C\wr G^\cC$, and this wreath product is the largest subgroup of permutations on $X$ leaving
$\cC$ invariant, and inducing the same top group and bottom group as $G$. We will therefore assume
that $G\leq G^C\wr G^\cC$.

In any family of partitions of $X$ that contains the two trivial partitions, we call a partition
{\em minimal} if the only strict refinement in the family is the discrete partition with all
classes of size 1; and a partition $\cC$ in the family is {\em maximal} if the only partition in
the family that is strictly coarser than $\cC$ is the all-in-one partition (with a single class).

For a minimal $G$-invariant partition, the group $G^C$ induced on a class $C$ is primitive. For a
maximal $G$-invariant partition $\cC$, the group $G^\cC$ induced on the set of classes is
primitive. If $\cC$ is both minimal and maximal as a $G$-invariant partition of $X$, then we say
that the action of $G$ on $X$ is 2-\emph{step imprimitive} relative to $\cC$.  Minimality and
maximality of a point-partition $\cC$ invariant under the action of a line-transitive automorphism
group of a linear space can sometimes be infered from the parameters, as demonstrated by the next
two results from \cite{PraegerTuan02}.

\begin{theorem}[Praeger and Tuan~\cite{PraegerTuan02}, Theorem 1.2]\label{thm:refinement} Suppose 
that $G$ is a line-transit\-ive point-im\-pri\-mi\-tive automorphism group
of a linear space. Let $\cC$ be a non-trivial $G$-invariant partition
with $d$ classes of size $c$ and let $x$ and $y$ be the 
Delandtsheer--Doyen parameters with respect to $\cC.$
Suppose there is another $G$-invariant partition $\cC'$ refining $\cC$
and let this partition have $d'$ classes of size $c'$ with 
Delandtsheer--Doyen parameters $x'$ and $y'.$ Then
\begin{enumerate}
\item
$c'$ divides $x - x'$ and $d$ divides $y' - y$;
\item
$x \ge 3x' + 1$ and $y' \ge 3y + 1$;
\item
$2r (xy - x'y') = (x'-x + y' - y) k;$
\item
$\displaystyle \frac{2r}{k} = \frac{c-c'}{x-x'}.$
\end{enumerate}
\end{theorem}
\begin{proof}
(i)-(iii) refer to  ~\cite[Theorem 1.2]{PraegerTuan02}. (iv) Double
counting the number of $\cC$-inner, $\cC'$-outer pairs gives
$b(x-x')=d'\frac{c'(c-c')}{2}=v\frac{c-c'}{2}$. It follows that
$\frac{2r}{k} =\frac{2b}{v}= \frac{c-c'}{x-x'}$.
\end{proof}

\bigskip

\begin{corollary}[Praeger and Tuan~\cite{PraegerTuan02}]\label{cor:maxmin} If $G$ is a line-transitive point-imprimitive automorphism group
of a linear space and $\cC$ is a $G$-invariant partition of the point set
with Delandtsheer-Doyen parameters $x$ and $y,$ then
\begin{enumerate}
\item
if $\binom{k}{2} > (x-2) xy + x + y$ then $G^C$ is primitive;
\item
if $\binom{k}{2} > (y-2) xy + x + y$ then $G^\cC$ is primitive.
\end{enumerate}
In particular, if $x \le 4$ then $G^C$ is primitive, and if $y \le 4$ then
$G^\cC$ is primitive.
\end{corollary}

\subsection{Normal partitions and quasiprimitive groups}\label{sub:normal}

The kernel of $G$ on $\cC$ is the subgroup $G_{(\cC)}$ of elements
$g \in G$ with $C^g = C$ for each $C\in\cC$. Thus
$G^\cC \simeq G/G_{(\cC)}$.
We say that $\cC$ is {\em $G$-normal} if 
$G_{(\cC)}$ is transitive on each of the classes of $\cC.$
Note that the set of orbits of a normal subgroup $N$ of $G$
always forms an invariant partition with
$N \le G_{(\cC)}$ transitive on each of the classes.
However, in general, not every
$G$-invariant partition arises as the set of orbits of a normal subgroup.
If no non-trivial 
$G$-normal partition exists, then every normal subgroup distinct
from the identity subgroup acts transitively.
In this case we say that $G$ is {\em quasiprimitive}
on the set $X.$
We then have $G^\cC  \simeq G$ since
the kernel $G_{(\cC)}$ on the set of classes is trivial in this case.

Summarizing the above remarks:\  \emph{for an imprimitive permutation group $G$
on $X$, either $G$ is quasiprimitive on $X$, or there is a non-trivial  
$G$-normal partition of $X$.}

We now turn to the context in which $X$ is the point set $\cP$ 
of a linear space and
$G$ is a line-transitive, point-imprimitive group of automorphisms. 
In a study in \cite{CaminaPraeger01} of the case where $G$ was assumed to be 
point-quasiprimitive (that is, quasiprimitive on $\cP$), it was
proved that $G$ must be {\em almost simple}, that is, 
there is a non-abelian simple group $T$ such that $T \le G \le \Aut(T).$
However, no examples of this kind are known. In fact, in 
all the known examples, there is a non-trivial $G$-normal point-partition 
(see \cite[Section 5]{CaminaPraeger01}).

\begin{theorem}[Camina and Praeger~\cite{CaminaPraeger01}]\label{thm:CP01} Let $G$ be a 
line-transitive point-imprim\-itive\\ group of automorph\-isms of a linear space. Then either
\begin{enumerate}
\item
there exists a non-trivial $G$-normal point-partition; \ \ or
\item
$G$ is point-quasiprimitive and almost simple.
\end{enumerate}
\end{theorem}

In an earlier study of the $G$-normal case, Camina and Praeger showed that an
intransitive normal subgroup $N$ must act faithfully on each of its orbits $C$,
that is to say, $N_{(C)}=1$ so that $N\simeq N^C$.

\begin{theorem}[Camina and Praeger~\cite{CaminaPraeger93}]\label{thm:CP93} If $G$ is a 
line-transitive group of automorph\-isms of a linear space $\cS = (\cP, \cL)$, and $N$ is a normal 
subgroup that is intransitive on $\cP$, then, for each $N$-orbit $C$ in $\cP$, 
\begin{enumerate}
\item $N$  acts faithfully on $C$; \ \ and in particular,
\item if $N$ is abelian then $|N|=|C|$ is odd.
\end{enumerate}
\end{theorem}

\section{Outline of the search strategy}\label{sect:outline}

In this section we describe briefly our approach to the search for line-transitive,
point-imprimitive linear spaces for which the Fang-Li parameter $k^{(r)}$ is at most some maximum
value $k^{(r)}_{\max}$. In particular, we will apply this approach in the case where
$k^{(r)}_{\max}=8$. Throughout the rest of the paper we assume the following.

\medskip\noindent
{\sc Hypothesis.}\label{hyp}\quad 
Let $\cS=(\cP,\cL)$ be a non-trivial linear space, admitting a line-transitive
point-imprimit\-ive automorphism group $G\leq\Aut (\cS).$ Let $v, k, b, r$ be as in 
(\ref{eqn:vrkb})--(\ref{eqn:fisher}). Let $\cC$ be a non-trivial
$G$-invariant partition of $\cP$ with $d$ classes of size $c$, and 
let $x$ and $y$ be the Delandtsheer--Doyen parameters, and $\gamma, \delta$ 
be the Fang-Li parameters, corresponding to $\cC$.

\smallskip
First we make explicit our aims and the nature of the output, and then we describe how the search 
proceeds in five broad steps.

\subsection{Aim and nature of the output}

\bigskip
\noindent
{\sc Aim:}\quad  Find each linear space $\cS=(\cP,\cL)$ such that 
the Fang-Li parameter $k^{(r)}\leq k^{(r)}_{\max}$, and there 
exists a group $G\leq \Aut(\cS)$ that is line-transitive 
and point-imprimitive. 

\bigskip
By Theorem~\ref{thm:CP01}, the search must take into account the following
two cases.

\smallskip\noindent
{\sc Case 1.}\quad $G$ is point-quasiprimitive and almost simple;
and

\noindent
{\sc Case 2.}\quad  $G$  leaves invariant a 
non-trivial $G$-normal point-partition.

\medskip
If $G$ is point-quasiprimitive then $G$ is faithful on all non-trivial 
$G$-invariant partitions, so we search only for the maximal ones. 
Thus the output from searching these cases will be as follows.

\smallskip\noindent
{\sc Output from Case 1.}\quad  all $(\cS, G, \cC)$, where $G$ is point-quasiprimitive and almost 
simple, and $\cC$ is a maximal $G$-invariant point-partition.

\noindent
{\sc Output from Case 2.}\quad  all $(\cS, G, \cC)$, where $\cC$ is a non-trivial
$G$-normal partition. 

Note that the output from these two cases will produce all pairs $(\cS, G)$ where $G$ is 
line-transitive and point-imprimitive on a linear space $\cS$. The search will not necessarily 
identify every non-trivial $G$-invariant partition for a given $(\cS, G)$. In particular, if $G$ is 
not point-quasiprimitive but acts faithfully on a non-trivial point partition $\cC$ then this 
partition will not be identified. See also Remark~\ref{rem:1}.

In the first part of the search we treat these two cases together, since the tests we apply are
valid whether or not $G$ is point-quasiprimitive. 

\subsection{Descriptions of {\sc Steps} in Search}

Suppose that $\cS=(\cP,\cL)$ is a non-trivial
linear space, $G$ is a line-transitive group of automorphisms, and $\cC$ is a non-trivial
$G$-invariant partition of $\cP$ with $d$ classes of size $c$. In Section~\ref{sect:params} we
define some additional parameters that give extra information and restrictions, and extract a series
of tests that form Algorithm~\ref{alg1}. We also present Algorithm~\ref{alg2} that tests for certain
sufficient conditions under which the partition $\cC$ is guaranteed to be minimal or maximal.

\smallskip\noindent
{\sc Step 1. (Fang-Li Parameter Sift)}\quad \\
Apply Algorithm~\ref{alg1} and output 
{\sc ParameterList}$(k^{(r)}_{\max})$, a list of parameter values that 
pass all these tests with $k^{(r)}\leq k^{(r)}_{\max}$. 
Each {\sc Line} of  {\sc ParameterList}$(k^{(r)}_{\max})$ gives values of 
the following parameters: 

\begin{enumerate}
\item  $d$ and $c,$ and hence $v = d c;$
\item  the Fang-Li parameters
$k^{(v)}$, $k^{(r)}$, $b^{(v)},$ $b^{(r)}$, as in Subsection~\ref{sub:fangli},
and hence $k,$ $r,$ $b,$;
\item the Delandtsheer-Doyen parameters $x$ and $y$ as in 
(\ref{eqn:DDreversed});
\item the Fang-Li parameters $\gamma$ and $\delta$, as in 
Proposition~\ref{prop:fl2};
\item the intersection type $(0^{d_0},1^{d_1},\ldots, k^{d_k})$ and hence also 
$\spec\,\cS$, as defined in Subsection~\ref{sub:inter}.\;
\item an upper bound $t_{\max}$ for the transitivity of the `top group' $G^\cC$, as defined in 
Definition~\ref{def:tmax}, (that is to say, if $G^\cC$ is $t$-transitive then $t\leq t_{\max}$).
\end{enumerate}

\smallskip\noindent {\sc Step 2. (Two-Step Imprimitivity Test)}\quad 
To each {\sc Line} of {\sc 
ParameterList}$(k^{(r)}_{\max})$ apply
Algorithm~\ref{alg2} and partition the list from {\sc Step} 1, {\sc ParameterList}$(k^{(r)}_{\max})$, into two parts: {\sc
ParameterListA}$(k^{(r)}_{\max})$ consists of those {\sc Lines} for which $G$ is guaranteed to be 2-step imprimitive
relative to $\cC$, while {\sc ParameterListB}$(k^{(r)}_{\max})$ contains all the remaining {\sc Lines}.

\smallskip
It turned out that for all the {\sc Lines} in {\sc ParameterList}$(8)$, the partition $\cC$ was
both minimal and maximal, that is to say, {\sc ParameterListB}$(8)$ was empty. Thus in each {\sc
Line} of {\sc ParameterListA}$(8)$ we have that the top group $H:=G^\cC$ is primitive and also the
bottom group $L:=G^C$ (where $C\in\cC$) is primitive, and $G\leq L\wr H$.

Next we enhanced the information contained in each {\sc Line} of {\sc
ParameterListA}$(k^{(r)}_{\max})$ by applying tests that restricted the possibilities for the
groups $H$ and $L$. In some cases these tests eliminated the {\sc Line} as no possibilities for
one of $H$, $L$ remained. These tests are described in Section~\ref{sect:groups} and form
Algorithms~\ref{alg:topgroup} and \ref{alg:bottomgroup}.  These algorithms may also be applied to
candidate parameter values in {\sc ParameterListB}$(k^{(r)}_{\max})$ to determine candidate top
groups in the case of a maximal $G$-invariant partition, or for candidate bottom groups in the
case of a minimal $G$-invariant partition. The remaining {\sc Steps} 4 and 5 that we describe, 
however,
are focused on dealing with {\sc ParameterListA}$(k^{(r)}_{\max})$. 

Some of the tests depend on the availability of a complete list of all primitive groups of a given degree, and we had
available to us, through the computer system {\sf GAP} \cite{GAP4} such lists for degrees less than 2,500. Thus the next
step could be applied completely only in those cases where both the number $d$ of classes and the class size $c$ were
less than 2,500.

\smallskip\noindent
{\sc Step 3. (Top and Bottom Group Sifts and Grid Test)}\quad Apply Algorithm~\ref{alg:topgroup} to 
give 
information in a list called {\sc TopGroups} about the
primitive `top group' $H=G^{\cC}$. For each {\sc Line} in {\sc ParameterListA}$(k^{(r)}_{\max})$ for which a list {\sc
Prim}$(d)$ of all primitive permutation groups of degree $d$ is available, {\sc TopGroups} will contain an explicit list
of candidates for $H$. If {\sc Prim}$(d)$ is not available then {\sc TopGroups} will contain this information, and
possibly some restrictions on $H$.

Apply Algorithm~\ref{alg:bottomgroup} to give information in a list called {\sc BottomGroups} about the primitive `bottom
group' $L=G^{C}$. For each {\sc Line} in {\sc ParameterListA}$(k^{(r)}_{\max})$ for which a list {\sc Prim}$(c)$ of all
primitive permutation groups of degree $c$ is available, {\sc BottomGroups} will contain an explicit list of candidates
for $L$.  If {\sc Prim}$(c)$ is not available then {\sc BottomGroups} will contain this information, and possibly some
restrictions on $L$.

Then, in Algorithm~\ref{alg:step3and4} we do the following. If, for a {\sc Line} of {\sc
ParameterlistA}$(k^{(r)}_{\max})$, one of {\sc TopGroups} or {\sc BottomGroups} is empty then this
{\sc Line} is removed from the list {\sc ParameterListA}$(k^{(r)}_{\max})$. At this point we will
have, in {\sc ParameterListA}$(k^{(r)}_{\max})$, {\sc Lines} representing all candidates for
parameters for any $(\cS, G, \cC)$ satisfying the {\sc Hypothesis}. 

We make one further pass through this list before dividing our search into the two {\sc Cases}. For
each {\sc Line} representing a possible $(\cS,G, \cC)$, we define {\sc PossibleGrid} to be `yes' if
either $c=d$, or there is another {\sc Line} in {\sc ParameterListA}$(k^{(r)}_{\max})$ corresponding
to the same value of $k$ and to a point-partition with $c$ classes of size $d$. Otherwise the value
of {\sc PossibleGrid} is defined as `no'. We add the value of {\sc PossibleGrid} to the {\sc Line},
(see Remark~\ref{rem:1}). This completes {\sc Step 3}.

\smallskip
In the final {\sc Steps} we divide the search into the two {\sc Cases}. For each {\sc Line} of the 
list
{\sc ParameterListA}$(k^{(r)}_{\max})$ for which an explicit list {\sc TopGroups} is available
(which means for us, $d<2,500$), Algorithm~\ref{alg:qptopgroup} produces candidates for almost
simple point-quasiprimitive groups $G$, and for each {\sc Line} of {\sc
ParameterListA}$(k^{(r)}_{\max})$ for which an explicit list {\sc BottomGroups} is available (for
us this means $c<2,500$), Algorithm~\ref{alg:gnormalbottomgroup} produces candidates for the
primitive bottom group in the case where $\cC$ is a $G$-normal partition.

\smallskip\noindent
{\sc Step 4. (Quasiprimitive and $G$-Normal Sifts)}\quad Apply Algorithm~\ref{alg:qptopgroup} to 
each {\sc Line} of {\sc ParameterListA}$(k^{(r)}_{\max})$ to
obtain a list {\sc QuasiprimTopGroups}: if {\sc TopGroups} contains an explicit list of groups, then {\sc
QuasiprimTopGroups} is a list of candidate primitive top groups in the case where $G$ is quasiprimitive, and otherwise
is a list containing restrictions for this case. Apply Algorithm~\ref{alg:gnormalbottomgroup} to each {\sc Line} to
obtain a list {\sc GNormalBottomGroups}: if {\sc BottomGroups} contains an explicit list of groups, then {\sc
GNormalBottomGroups} is a list of candidate bottom groups in the case where $\cC$ is $G$-normal, and otherwise is a list
containing restrictions for this case. 

\medskip
Our final step we state only for $k^{(r)}_{\max}=8$.

\smallskip\noindent
{\sc Step 5. (Analysing Remaining Lines)}\quad We make a {\sc Line}-by-{\sc Line} consideration of 
the output {\sc 
ParameterListA}$(8)$ of {\sc Step} 4 to determine all the possible linear spaces, and thereby 
complete the proof of Theorem~\ref{main}.

\begin{remark}\label{rem:1}{\rm Finally we make some comments on the need for the parameter {\sc
PossibleGrid}, and on its name.  Consider a situation in which a group $G$ leaves invariant two
non-trivial partitions of $\cP$, namely a partition $\cC$ with $d$ classes of size $c$, and a
partition $\cC'$ with $c$ classes of size $d$. Suppose also that the associated groups $G^\cC,
G^{\cC'}$, and the groups induced on classes are all primitive. Then each class of $\cC'$ consists
of one point from each of the classes of $\cC$, and each class of $\cC$ consists of one point from
each of the classes of $\cC'$. Thus the point set $\cP$ can be identified with the Cartesian
product $\cC\times\cC'$, and $G$ is a subgroup of the full stabiliser $\Sym(\cC)\times \Sym(\cC')$
of this `grid structure'. It is possible that the group $G$ acts faithfully on one of these
partitions, say $\cC'$, and is not faithful on the other partition $\cC$. (For a simple example in
the group setting, take $|\cC|=2$ and $|\cC'|=5$. Then $S_2 \times S_5$ contains a transitive
subgroup isomorphic to $S_5$ with these properties.) In this case, $G$ is not quasiprimitive, and
$\cC$ is $G$-normal but $\cC'$ is not $G$-normal. Thus after the next steps the triple $(\cS, G,
\cC)$ will be identified in {\sc Case} 2, but the triple $(\cS, G,\cC')$ will not appear in either
{\sc case}. This is in line with our aim to find all possible $\cS$. However, when identifying the
$G$-normal triple $(\cS, G, \cC')$, the information that there was a possible grid structure could
be useful at {\sc Step 5}.  In other {\sc Lines} it could also be important at {\sc Step 5} to
know that it is impossible for such a grid structure to be preserved by the group.
}\end{remark}

\section{Parameter restrictions and 2-step imprimitivity}\label{sect:params}

Let $(\cS, G, \cC)$ be as in the {\sc Hypothesis}. In this section we present several results that 
restrict the parameters and provide sufficient conditions for 2-step imprimitivity. These results 
are used to design the {\sc Fang-Li Parameter Sift} and the {\sc Two-Step Imprimitivity Test} in 
Subsections~\ref{sec:flsift} and \ref{sub:2step}.

\bigskip
First we record some equalities and inequalities relating these parameters. The last statement of 
part (ii) was proved in \cite{DelandtsheerNiemeyeretal01}.

\begin{theorem}[Praeger and Tuan~\cite{PraegerTuan02}, Theorem 1.1]\label{thm:CNP}
Assume that the {\sc Hypothesis} holds. 
Then
\begin{enumerate}
\item 
$c = \frac{2xr}{k} + 1,$
$d = \frac{2yr}{k} + 1$,
and
$\frac{b}{v} = \frac{r}{k} = \frac{\binom{k}{2} - x - y}{2xy}.$
\item 
$\binom{k}{2} \ge 2xy + x + y,$
$c \ge 2x + 1$
and
$d \ge 2y + 1,$
and equality holds in one of these relations if and only if equality holds in all three,
and this occurs if and only if $\cS$ is a projective plane.
In particular, if $c = 3$ or $d=3$ then $\cS$ is a projective plane.
\item 
At least one of $k \ge 2x$ and $k \ge 2y$ holds. Moreover
\begin{enumerate}
\item
if $k \ge 2x$ then $k - 2x \ge 2,$ $y \le \binom{k - 2x}{2},$ and $d \ge 2k -2x - 1 > k;$
\item
if $k \ge 2y$ then $k - 2y \ge 2,$ $x \le \binom{k - 2y}{2},$ and $c \ge 2k -2y - 1 > k.$
\end{enumerate}
In particular $k < \max\{c,d\}.$
\end{enumerate}
\end{theorem}

\smallskip
Under additional conditions,
equality holds in the bound for $y$ in (iii) (a), namely,
$y = \binom{k-2x}{2}.$ See Theorem~\ref{thm:normal}.

\subsection{Intersection numbers and spectrum}\label{sub:inter}

Let $\lambda \in \cL.$
Define the intersection numbers
\begin{equation}\label{eqn:di}
d_i = \Big|  \big\{ \; C \in \cC: \;
|C \cap \lambda|=i \big\}\Big|,
\end{equation}
which are independent of the choice of the line $\lambda$.
The {\em intersection type} is the vector $(0^{d_0},1^{d_1},\ldots, k^{d_k})$
and the {\em spectrum} is the set of non-zero intersection sizes
\begin{equation}\label{eqn:spec}
\spec\,\cS := \{ i > 0 \mid d_i \neq 0 \}.
\end{equation}
We sometimes write $\spec_\cC \, \cS$ if we need to specify the partition $\cC.$

\smallskip

Of particular interest are the smallest and the largest elements of this set:
\begin{eqnarray*}
i_{\min} &:=& \min(\spec \, \cS),
\qquad \text{and} \qquad
i_{\max} := \max(\spec \, \cS).
\end{eqnarray*}

\smallskip
For a point $\al$, $C(\al)$ denotes the class of $\cC$ containing $\al$.
A point $\al$ and a line $\lambda$ are said to be \emph{$i$-incident} 
if $|\lambda \cap C(\al) | = i$.
A class $C$ and a line $\lambda$ are \emph{$i$-incident} if 
$|C \cap \lambda | = i.$
Transitivity of $G$ on classes and on points implies that the following 
numbers are well defined, that is, they do not depend on the choice of the 
class $C \in \cC$ and the choice of the point $\al \in \cP:$
\begin{align*}
b_i & = \text{ the number of lines which are $i$-incident with a class $C$}, \\
r_i & = \text{ the number of lines which are $i$-incident with a point $\al$}.
\end{align*}

\smallskip

Recall that for two points $\alpha$ and $\beta,$
$\lambda(\alpha,\beta)$ denotes the unique line joining $\alpha$
with $\beta$.

\smallskip

Next we give a few properties of $\spec\,\cS$. As was noted in
\cite{Delandtsheer89}, part (i) of Theorem~\ref{thm:HM} below can be derived
from the proof of \cite[Proposition 3]{HigmanMcLaughlin61}, but for
completeness we give a short proof. For any subset $F \subseteq \cP$
containing at least two points, define

\begin{equation}\label{eqn:induced}
\cL|_F := \{ \lambda \cap F \;:\; \lambda \in \cL, \; |\lambda \cap F| \ge 2\}.
\end{equation}

Then the induced incidence structure $\cS|_F := (F, \cL|_F)$ is a
(possibly trivial) linear space, called the {\em linear space induced on
$F.$} Under certain conditions on the parameters, $\cS|_F$ has constant
line size.

\begin{theorem}\label{thm:HM} 
Assume that the {\sc Hypothesis} holds. Then
\begin{enumerate}
\item $|\spec \, \cS| \ge 2;$
\item $c\not\in\spec\,\cS$;
\item if $\spec \, \cS = \{ 1,h\}$ with $h \ge 2,$ then the induced 
linear space $\cS|_C$ has constant line size $h$. (If $h$ is $2$, this 
structure is essentially a complete graph with $c$ vertices.)
\end{enumerate}
\end{theorem}

\begin{proof}
Suppose that $\spec\, \cS = \{h \} $. Since two points of $C$ lie on a
line, $h\geq2$. For a point $\al\in C$, the set 
$\lambda_1,\dots,\lambda_r$ of lines containing $\al$ gives rise to a
partition $\{ (C\cap\lambda_i)\setminus\{\al\}\,:\,1\leq i\leq r\}$ of
$C\setminus\{\al\}$, whence $c =r (h-1) + 1$. Hence $v=r(k-1)+1 =
d(r(h-1)+1)$, which implies that $r\big((k-1)-d(h-1)\big) = d-1$. Since
$d \geq 2$ and $h\geq2$, this implies that $r \leq d-1 < d \leq k-1 $,
so that $ r < k$, contradicting (\ref{eqn:vrkb}) and (\ref{eqn:fisher}).
Thus part (i) is proved. 

Part (ii) is proved in \cite[Corollary 2.2]{DelandtsheerNiemeyeretal01},
and part (iii) is proved in \cite[Proposition
2.3(iii)]{DelandtsheerNiemeyeretal01}.
\end{proof}

\smallskip

The next lemma is an immediate consequence of a result of Camina and Siemons~\cite[Lemma
2]{CaminaSiemons89a}. It gives sufficient conditions for $\cS|_F$ (as defined 
above by (\ref{eqn:induced})) to be line-transitive.

\begin{lemma}[Camina and Siemons~\cite{CaminaSiemons89a}]\label{lem:CaminaSiemons} Assume that 
$\cS=(\cP,\cL)$ is a non-trivial linear space admitting a line-transitive automorphism group $G$. 
Let $\lambda \in\cL$ and $H\leq G_{\lambda}$ such that, for $F:=\Fix_{\cP}(H)$,

\begin{enumerate}
\item $2\leq |F\cap \lambda| <|F|$, and
\item if $K\leq G_\lambda$ and $|\Fix_{\cP}(K)\cap \lambda|\geq 2$, and 
$H$ and $K$ are conjugate in $G$, then $H$ and $K$ are conjugate in 
$G_\lambda$.
\end{enumerate}

Then $\cS|_F$ has constant line size and $N_G(H)$ acts line-transitively 
on $\cS|_F$.
\end{lemma}

\begin{corollary} \label{cor:CaminaSiemons}
Assume that the {\sc Hypothesis} holds. Let $\lambda\in\cL$ and $p$ be a 
prime dividing $|G_{\lambda}|$. Let $P$ be a Sylow $p$-subgroup of 
$G_\lambda$ and $F:=\Fix_{\cP}(P)$. Suppose that $2\leq 
|F\cap\lambda|<|F|$. Then 
\begin{enumerate}
\item $N_G(P)$ is line-transitive on $\cS|_{F}$;
\item $\cC|_F:=\{C\cap F:C\in\cC,C\cap F\ne\emptyset\}$ is an 
$N_G(P)$-invariant partition of $F$;
\item $|F|=f.|C\cap F|$ where $f=|\cC|_F|,\,C\cap F\in\cC|_F$, and $|C\cap F|\geq 3$.
\end{enumerate}
\end{corollary}

\begin{proof}
Part (i) follows from Lemma~\ref{lem:CaminaSiemons}, and so by \cite{Block68}, $N_G(P)$ is
transitive on $F$. Consider $C\cap F$, for some class $C\in\cC$ such that $C\cap F\ne\emptyset$.
Then $(C\cap F)^g=C^g\cap F^g=C^g\cap F$ for any $g\in N_G(P)$. Since $C$ is a block of
imprimitivity for $G$, $C^g=C$ or $C^g\cap C=\emptyset$, and so $(C\cap F)^g=C\cap F$ or $(C\cap
F)^g\cap(C\cap F)=\emptyset$. Thus $\cC|_F$ forms an $N_G(P)$-invariant partition of $F$, and in
particular, $|F|=f.|C\cap F|$ where $f$ is the number of classes in $\cC$ that contain fixed points
of $P$, and hence are fixed by $P$.  Moreover, $|C\cap F|\ne 2$ by Theorem~\ref{thm:CNP}(ii) so 
$|C\cap F|\geq 3$.
\end{proof}

Finally in this subsection we collect some useful arithmetical
relationships between these parameters. Some parts were proved
in~\cite{DelandtsheerNiemeyeretal01} and \cite{PraegerTuan02}.
\smallskip

\begin{sloppypar}
\begin{proposition}\label{prop:intersection}
The following all hold.
\begin{enumerate}
\item 
$\displaystyle d = \sum_{i=0}^k d_i,$
$\displaystyle b = \sum_{i=0}^k b_i,$
$\displaystyle cr = \sum_{i=1}^k ib_i,$
$\displaystyle r = \sum_{i=1}^k r_i,$
$\displaystyle k = \sum_{i=1}^k i d_i,$
$\displaystyle x = \sum_{i=2}^k \binom{i}{2}d_i,$
\item 
$\displaystyle c-1 = \sum_{i=1}^k (i-1)r_i,$ $\displaystyle v-c =
\sum_{i=1}^k (k-i)r_i,$ $\displaystyle \binom{c}{2} = \sum_{i=1}^k
\binom{i}{2} b_i,$ $\displaystyle \binom{v-c}{2} = \sum_{i=0}^k
\binom{k-i}{2} r_i,$ $\displaystyle c (v-c) = \sum_{i=1}^{k-1} i
(k-i)r_i,$
\item 
$\displaystyle b d_i = d b_i$, $\displaystyle c r_i=i b_i $,  and
$\displaystyle v r_i = b i d_i$ for $0 \le i \le k,$
\item 
$\displaystyle \frac{c d_i}{k^{(v)}} = \frac{b_i}{b^{(r)}}$
and
$\displaystyle r_i' := \frac{r_i}{b^{(r)}} = \frac{i d_i}{k^{(v)}}$
are integers for $0 \le i \le k,$
\item 
$\displaystyle \beta := \frac{d}{\gcd(d,b^{(v)})} \; \big| \; d_i$
for $0 \le i \le k,$ hence for all these $i,$
$\displaystyle d_i' := \frac{d_i}{\beta},$
$\displaystyle \frac{k}{\beta} = \sum_{i=0}^k i \cdot d_i',$ 
and
$\displaystyle \frac{2 x}{\beta} = \sum_{i=2}^k i (i-1) d_i'$ are all integers. 
\item 
$\displaystyle \alpha_i := \frac{k^{(v)}}{\gcd(i,k^{(v)})} \;\Big|\; d_i$ 
for $0\le i\le k,$
and hence $\displaystyle \frac{\alpha_i }{\gcd(\alpha_i, \beta)} \;\Big|\; 
d_i'.$
\item 
$i_{\max} < \min \{ \sqrt{c} + \frac{1}{2},\; \sqrt{2x} + 1\}.$
\item 
Given a point $\alpha,$ the number of points $\beta \neq \alpha$ such that $\lambda(\alpha,\beta)$
is $1$-incident with both $\alpha$ and $\beta$ is $r_1(d_1-1) = \frac{r}{k}d_1(d_1 - 1)$
(we call such configurations 1-incident point-line-point triples).
\item 
$d_1 \ge k - 2x$ with equality if and only if $\spec \, \cS = \{ 1,2\},$
that is, if and only if 
the intersection type is $(0^{d-k+x},1^{k-2x}, 2^x, 3^0, \ldots, k^0).$
In particular,  if $k \ge 2x$ then $d_1 > 0.$
\end{enumerate}
\end{proposition}
\end{sloppypar}

\begin{proof}
(i)\quad  The equations follow from the definitions of the $d_i$,
$b_i$ and $r_i$ (also see~\cite{DelandtsheerNiemeyeretal01}).

(ii)\quad The first equality is proved in ~\cite[Proposition
2.4(iv)]{DelandtsheerNiemeyeretal01}. The second one follows
from this and the fact that $v-1=r(k-1)$. We obtain the third
equality by counting, for a fixed class $C$, the number of pairs
$(\{\alpha,\beta\},\lambda)$ for which $\alpha,\beta\in C\cap
\lambda$. For the fourth equation, we count, for a fixed class $C$,
the number of pairs $(\{\alpha,\beta\},\lambda)$ where
$\alpha,\beta\in \lambda$ and $\alpha,\beta\in {\cal
P}\setminus\{C\}$. Similarly, counting, for a fixed class $C$, the
number of pairs $(\{\alpha,\beta\},\lambda)$ where
$\alpha,\beta\in \lambda$, and $\alpha \in C$ , $\beta\in {\cal
P}\setminus\{C\}$ in two ways yields the last equality.

(iii) and (iv) are proved in~\cite[Proposition
2.3]{DelandtsheerNiemeyeretal01}.

(v)\quad By (iii), $\displaystyle\frac{db_i}{b}=d_i$. Noting that
$d$ (as a divisor of $v$) is prime to $b^{(r)}$, we must have that
$\beta:=\displaystyle\frac{d}{gcd(d,b^{(v)})}$ divides $d_i$ for 
$0\leq i\leq k$. Note that this number does not depend on $i$. The
equation $k=\sum_{i=1}^{k}id_i$ from (i) reduces to
$\frac{k}{\beta}=\sum_{i=1}^{k}id_i'$. Similarly, the equation
$2x=\sum_{i=2}^{k}i(i-1)d_i$ of (i) becomes
$\frac{2x}{\beta}=\sum_{i=2}^{k}i(i-1)d_i'$.

(vi)\quad By (iv), $r_i'=\frac{id_i}{k^{(v)}}$ is an integer, and
therefore $\frac{k^{(v)}}{gcd(k^{(v)},i)}$ a divisor of $d_i$.

(vii) is proved in~\cite[Lemma 3.1]{DelandtsheerNiemeyeretal01}.

(viii) is proved in~\cite[Lemma 2.3]{PraegerTuan02}.

(ix) By part (i), $k-2x=\sum_{i=1}^kid_i-\sum_{i=2}^ki(i-1)d_i =d_1-\sum_{i=2}^ki(i-2)d_i\leq d_1$.
Equality holds if and only if $d_i=0$ for all $i\geq3$, that is, if and only if
$\spec\,\cS\subseteq\{1,2\}$. Since $|\spec\,\cS|\geq2$, by Theorem~\ref{thm:HM}~(i), this
condition is equivalent to $\spec\,\cS=\{1,2\}$, which in turn is equivalent to the intersection
type being $(0^{d-k+x},1^{k-2x}, 2^x, 3^0, \ldots, k^0).$ Now suppose that $k\geq 2x$. If $k>2x$
then $d_1\geq k-2x>0$, so suppose that $k=2x$. If $d_1=0$ then $d_1=k-2x=0$, and we have just shown
that in this case the intersection type is $(0^{d-k+x},1^{k-2x}, 2^x, 3^0, \ldots, k^0)$, which
implies that $|\spec\,\cS|=1$, contradicting Theorem~\ref{thm:HM}. 
\end{proof}

\bigskip
Hence the intersection type $(0^{d_0},1^{d_1},\ldots, k^{d_k})$ determines the
numbers $r_i$ and $b_i$ via (iii):
\begin{align}\label{eqn:ribifromdi}
r_i = \frac{b i d_i}{v} = \frac{r i d_i}{k}, \quad \text{and} \quad b_i = \frac{cr_i}{i} = \frac{cr}{k} d_i.
\end{align}
In particular, the numbers $b_i$ are proportional to the $d_i.$

\subsection{Parameter restrictions using the top group}

Here we derive an upper bound for the transitivity of the top group $G^\cC$ that depends only on the
intersection type $(0^{d_0}, 1^{d_1}, \ldots, k^{d_k})$. Recall the definitions of the $d_i$ and the
spectrum $\spec\, \cS$ in (\ref{eqn:di}) and (\ref{eqn:spec}).

\begin{definition}\label{def:tmax}{\rm For a given intersection type
$(0^{d_0},1^{d_1},\ldots, k^{d_k})$, and non-empty subset $S\subseteq \spec \, 
\cS$, set $d(S):= \sum_{i\in S} d_i.$ Define $t_{\max}$ to be the largest 
positive integer $t$ such that, for all $S\subseteq \spec \, \cS$, and all
positive integers $h\leq\min\{t, d(S)\}$,
\[
\prod_{j=0}^{h-1} (d-j) \quad \mbox{divides} \quad 
b \prod_{j=0}^{h-1} (d(S)-j).
\]
}
\end{definition}

For $h=1$ and any non-empty subset $S$ of $\spec \, \cS$, the displayed condition in
Definition~\ref{def:tmax} is simply `$d$ divides $b\,d(S)$', and the truth of this follows
immediately from the definition of $d(S)$ and from the equalities $d b_i=b d_i$ which hold for all
$i\in S$ by Proposition~\ref{prop:intersection}~(iii).  Hence the condition in
Definition~\ref{def:tmax} holds for $t=1$, and so $t_{\max}$ is well-defined. We prove next that
$t_{\max}$ is an upper bound for the transitivity of $G^\cC$.

\begin{lemma}\label{lem:maxtrans}
If $G^\cC$ is $t$-transitive then $t\leq t_{\max}$.
\end{lemma}

\begin{proof}
As remarked above, $t_{\max}\geq1$, so if $t=1$ then the result holds. 
Assume that $t\geq2$.
Let $S\subseteq \spec \, \cS$, set $s:=d(S)$ as in Definition~\ref{def:tmax}, and let 
$h$ be such that  $h\leq\min\{t, s\}$.
We double count the set
\[
\cM = \Big\{ ((C_1,\ldots,C_h), \lambda) \in \cC^h \times \cL \;:\;
C_i \neq C_j \;\; \text{for} \;\; i \neq j,
\;\;
|C_i \cap \lambda| \in S \;\; \text{for} \;\; i=1,\ldots,h \Big\}.
\]
Since $G$ acts $t$-transitively on the set $\cC$ of classes,
it also acts transitively on the set of $h$-tuples of pairwise distinct
classes. Thus the number of possibilities for the line $\lambda$ in the second 
component, for a given $h$-tuple $(C_1,\dots,C_h)$,
does not depend on the choice of the $h$-tuple.
Call this number $n,$ so that the cardinality of the set $\cM$ is
$d(d-1)\cdots (d-h+1)\cdot n$.

On the other hand, given any line $\lambda,$ the number of choices
of $h$-tuples of pairwise distinct classes $(C_1,...,C_h),$
each intersecting $\lambda$ in a number of points belonging to $S$,
is $s(s-1)\cdots (s-h+1),$ so  $ |\cM|= b 
\cdot \prod_{j=0}^{h-1}(s-j).$ Hence $\prod_{j=0}^{h-1}(d-j)$
divides $b \cdot \prod_{j=0}^{h-1}(s-j).$ Thus the displayed condition of 
Definition~\ref{def:tmax} holds for $t$, and so $t\leq t_{\max}$.
\end{proof}

\bigskip
Although the next result will not be used until we consider more detailed
group theoretic properties, we derive here the following extension of 
Lemma~\ref{lem:maxtrans}.  

\begin{lemma}\label{lem:altd}
If $G^\cC \ge \Alt_d,$ then
$\displaystyle{\lcm_{i=1}^k \binom{d}{d_i}}$ divides $b.$
\end{lemma}

\begin{proof}
By assumption, $G^\cC$ is $t$-transitive, where $t=d-2$ if $G^\cC$ is $\Alt_d$,
and $t=d$ if $G^\cC=\Sym_d$. 
Let $i$ be such that $0<i\leq k$. By Theorem~\ref{thm:HM}~(i), $d_i<d$.
If $d_i=d-1$, then by Proposition~\ref{prop:intersection}~(iii), $d$ 
divides $bd_i=b(d-1)$, and hence $d=\binom{d}{d_i}$ divides $b$. Suppose now 
that $d_i\leq d-2$. Choose $S=\{i\}\subseteq \spec\,\cS$ so that 
$d(S)=d_i\leq d-2\leq t$, and choose $h =d_i= \min\{t,d(S)\}$. Then,
by Lemma~\ref{lem:maxtrans} and Definition~\ref{def:tmax},  ${
\prod_{j=0}^{d_i-1} (d-j)}$ divides ${
b\prod_{j=0}^{d_i-1} (d_i-j)}$. Now
${ \prod_{j=0}^{d_i-1} (d-j) = \frac{d!}{(d-d_i)!}}$
and ${ b\prod_{j=0}^{d_i-1} (d_i-j)} =bd_i!.$ It follows that
$\binom{d}{d_i}$ divides $b$. As this is the case for all $i$ with
$0<i\leq k$, the result follows.
\end{proof}

Finally in this subsection we record a test for determining $t_{\max}$, that 
essentially just checks the condition of Definition~\ref{def:tmax}.

\begin{algorithm}\label{algtmax}{\rm ({\sc TMax})\\
{\sc Input:}\quad The values $d,c,k$ and an intersection type $(0^{d_0},1^{d_1},\ldots, k^{d_k})$.

\smallskip\noindent
{\sc Output:}\quad $t_{\max}$ and $\spec\,\cS$.

\smallskip\noindent
compute\ $\spec\,\cS :=\{i>0\;:\;d_i\ne0\}$;\\
compute \ $\mathcal{T} :=\{ \sum_{i\in S}d_i \;:\; S\subseteq\spec\,\cS,\,S\ne\emptyset\}$;\\
set $T:=0$;\\
for $t=1,\ldots,d$\\
\mm if, for some $s\in\mathcal{T}$ with $s\geq t$, \\
\mm \mm $d(d-1)\ldots(d-t+1)$ does not divide $bs(s-1)\ldots(s-t+1)$\\
\mm \mm skip to instruction $(\ast)$; \\
\mm else set $T:=T+1$;\\ 
($\ast$) set $t_{\max}:=T$ and return $t_{\max}$ and $\spec\,\cS$.
}
\end{algorithm}

\begin{lemma}\label{lem:algtmax} Let $d,c,k$ be given, and let $(0^{d_0},1^{d_1},\ldots, k^{d_k})$ be a corresponding
intersection type with $\sum_i id_i=k$. Then, along with $\spec\,\cS$, either Algorithm~\ref{algtmax} returns $0$ and
there are no linear spaces satisfying the {\sc Hypothesis} with these parameters, or it returns the correct value of
$t_{\max}$. \end{lemma}

\begin{proof} Let $\cS, \mathcal{T}$ be as in Algorithm~\ref{algtmax}. Suppose that $t'$ is the integer returned by
Algorithm~\ref{algtmax}. If $t'=0$ then the parameter $T$ is not increased during the run of the `if loop' with $t=1$,
and hence, for some $s\in\mathcal{T}$, $d$ does not divide $bs$, contradicting Proposition~\ref{prop:intersection}~(iii)
(see the discussion preceding Lemma~\ref{lem:maxtrans}). Thus in this case there is no linear space with these parameters
satisfying the {\sc Hypothesis}. 

Suppose now that $t'>0$. This means that each run of the `if loop' in which the parameter $t\leq t'$ finishes with $T$
being increased by $1$, so that at the end of all these runs of the `if loop' we have $T=t'$. Also, if $t'<d$, then in
the run of the `if loop' with $t=t'+1$ the value of $T$ is not increased and some instance of the divisibility condition
fails.  Let $s\in\mathcal{T}$, let $h$ satisfy $1\leq h\leq\min\{t',s\}$, and let $c(s,h)$ denote the condition
`$d(d-1)\ldots(d-h+1)$ divides $bs(s-1)\ldots(s-h+1)$'. Then in the run of the `if loop' with $t=h$ we would have
verified that $c(s,h)$ holds. Thus $t'\leq t_{\max}$ by Definition ~\ref{def:tmax}. If $t'=d$ then (since $t_{\max}\leq
d$), we must have $t'=t_{\max}$, so suppose that $t'<d$. As mentioned above, in the run of the `if loop' with $t=t'+1$,
$c(s,t'+1)$ fails for some $s\in\mathcal{T}$ with $s\geq t'+1$. Hence $t'=t_{\max}$ in this case also. \end{proof}

\subsection{{\sc Fang-Li Parameter Sift}}\label{sec:flsift}

In this subsection we describe Algorithm~\ref{alg1} that uses the results presented so far to sift for feasible parameter
sets for line-transitive, point-imprimitive linear spaces satisfying the {\sc Hypothesis} given at the beginning of this
section. Applying this algorithm will complete {\sc Step} 1 of the search procedure outlined in
Section~\ref{sect:outline}.

We apply tests to determine feasible values for the parameters $k, d, c$ and also the Fang--Li and the
Delandtsheer--Doyen parameters. In addition, we compute feasible intersection types, and the value of $t_{\max}$. We
restrict to the cases where
\begin{align}\label{restrictions}
k^{(r)} \le k^{(r)}_{\max} := 8.
\end{align}

The output of Algorithm~\ref{alg1} is a list {\sc ParameterList}$(k^{(r)}_{\max})$ of parameter sequences, called {\sc
Lines}, that satisfy the conditions given so far on the parameters of line-transitive, point-imprimitive linear spaces
where the Fang--Li parameter $k^{(r)}$ is at most some given upper bound $k^{(r)}_{\max}$.

\begin{algorithm}\label{alg1}{\rm ({\sc ParameterList})\\
{\sc Input:}\quad An upper bound $k^{(r)}_{\max}$ for the Fang--Li parameter 
$k^{(r)}$.

\smallskip\noindent
{\sc Output:}\quad The list {\sc ParameterList}$(k^{(r)}_{\max})$.

\smallskip\noindent
set {\sc ParameterList}$(k^{(r)}_{\max}):=$ an empty list; \\
for $k^{(r)} = 1 , \ldots, k^{(r)}_{\max}$  \\
\mm for $\gamma = 1, \ldots , {k^{(r)}}^2-1$ (Proposition~\ref{prop:fl2}(iii))\\
\mm \mm for $\delta = 1, \ldots$ with $\gamma\delta < {k^{(r)}}^2$ (Proposition~\ref{prop:fl2}(iii)) \\
\mm \mm \mm for all $k^{(v)}$ dividing $(\gamma + k^{(r)})(\delta + k^{(r)})$
(Proposition~\ref{prop:fl2}(iv))\\
\mm \mm \mm \mm $k := k^{(v)} \cdot k^{(r)}$ (Lemma~\ref{prop:fl1}(i)); \\
\mm \mm \mm \mm if $2 \nmid \gamma k^{(v)}$ skip to next $k^{(v)}$;\\
\mm \mm \mm \mm $x := \frac{\gamma k^{(v)}}{2}$ (Proposition~\ref{prop:fl2}(ii));\\
\mm \mm \mm \mm if $2 \nmid \delta k^{(v)}$ skip to next $k^{(v)}$;\\
\mm \mm \mm \mm $y := \frac{\delta k^{(v)}}{2}$ (Proposition~\ref{prop:fl2}(ii));\\
\mm \mm \mm \mm if $\gamma\delta \nmid (k^{(r)} (k - 1) - \gamma - \delta)$ skip to next $k^{(v)}$;\\
\mm \mm \mm \mm $b^{(r)} := \frac{k^{(r)} (k - 1) - \gamma - \delta}{\gamma\delta}$ (Proposition~\ref{prop:fl2}(iii));\\
\mm \mm \mm \mm if $k^{(v)} > b^{(r)}$ skip to next $k^{(v)}$ (Fisher's inequality (\ref{fisherfangli})); \\
\mm \mm \mm \mm $v := (\gamma + \delta + \gamma\delta b^{(r)}) \cdot b^{(r)} + 1$;\\
\mm \mm \mm \mm if $k^{(v)} \nmid v$ skip to next $k^{(v)}$;\\
\mm \mm \mm \mm $b^{(v)} := \frac{v}{k^{(v)}}$ (Lemma~\ref{prop:fl1}(i));\\
\mm \mm \mm \mm $b := b^{(v)} \cdot b^{(r)}$ (Lemma~\ref{prop:fl1}(i));\\
\mm \mm \mm \mm $r := k^{(r)} \cdot b^{(r)}$ (Lemma~\ref{prop:fl1}(ii));\\
\mm \mm \mm \mm $c := \gamma b^{(r)} + 1$ (Proposition~\ref{prop:fl2}(i));\\
\mm \mm \mm \mm $d := \delta b^{(r)} + 1$ (Proposition~\ref{prop:fl2}(i));\\
\mm \mm \mm \mm if $k \ge 2x$ \\
\mm \mm \mm \mm \mm if $y > \binom{k - 2x}{2}$ skip to next $k^{(v)}$ (Theorem~\ref{thm:CNP}(c)(i)); \\
\mm \mm \mm \mm if $k \ge 2y$ \\
\mm \mm \mm \mm \mm if $x > \binom{k - 2y}{2}$ skip to next $k^{(v)}$ (Theorem~\ref{thm:CNP}(c)(ii)); \\
\mm \mm \mm \mm compute all intersection types $(0^{d_0},1^{d_1},\ldots, k^{d_k})$ \\
\mm \mm \mm \mm \mm with respect to Proposition~\ref{prop:intersection}(v)-(vii) and Theorem~\ref{thm:HM};\\
\mm \mm \mm \mm if there are no intersection types, skip to next $k^{(v)}$; \\
\mm \mm \mm \mm for each intersection type $(0^{d_0},1^{d_1},\ldots, k^{d_k})$\\
\mm \mm \mm \mm \mm compute $t_{\max}$ and $\spec\,\cS$ (Algorithm~\ref{algtmax}); \\
\mm \mm \mm \mm \mm if $t_{\max}=0$ skip to next intersection type (Lemma~\ref{lem:algtmax});\\
\mm \mm \mm \mm \mm insert in {\sc ParameterList}$(k^{(r)}_{\max})$ the {\sc Line} \\
\mm \mm \mm \mm \mm \mm 
$(d,c,x,y,\gamma,\delta,k^{(v)},k^{(r)},b^{(v)},b^{(r)},(0^{d_0},1^{d_1},\ldots,k^{d_k}),t_{\max},\spec\,\cS)$; \\
\mm \mm \mm \mm \mm according to the value of $v$;\\
return {\sc ParameterList}$(k^{(r)}_{\max})$.
}
\end{algorithm}

Applying Algorithm~\ref{alg1} with $k^{(r)}_{\max}=8$ produced a list {\sc ParameterList}(8)
containing 1207 {\sc Lines}. We refer to {\sc Line} $i$ of {\sc ParameterList}(8) as {\sc
ParamList}${}_8(i)$. As noted in {\sc Step} 1 in Section~\ref{sect:outline}, {\sc
ParamList}${}_8(i)$ will give values of the following parameters: \[
(d,c,x,y,\gamma,\delta,k^{(v)}, k^{(r)}, b^{(v)}, b^{(r)}, (0^{d_0},1^{d_1},\ldots, k^{d_k}),
t_{\max}, \spec\, \cS). \] We give a shortened version of the list in Appendix~\ref{appx}, namely
we omit from each {\sc Line} the intersection type $(0^{d_0},1^{d_1},\ldots, k^{d_k})$, $t_{\max}$
and $\spec\, \cS$. For those {\sc Lines} that survive the further tests of {\sc Steps} $2-4$, the 
information about the intersection type, $t_{\max}$ and $\spec\, \cS$ is needed for the final 
analysis. Full information about these {\sc Lines} is presented in Section~\ref{sect:results}.

\subsection{Testing for $2$-step imprimitivity}\label{sub:2step}

As we noted in Subsection~\ref{sub:partitions}, information in Theorem~\ref{thm:refinement} and
Corollary~\ref{cor:maxmin} provides sufficient conditions for the partition $\cC$ to be maximal or minimal, and hence
also provides sufficient conditions for $G$ to be 2-step imprimitive relative to $\cC$ (that is, $\cC$ is both maximal
and minimal). We formalise these sufficient conditions in Algorithm~\ref{alg2}. Applying this algorithm will complete
{\sc Step} 2 of the search procedure described in Section~\ref{sect:outline}, the {\sc Two-Step 
Imprimitivity Test}. 

We need only the values for $(d,c,x,y,k^{(v)},k^{(r)})$ for each {\sc Line} of {\sc ParameterList}$(k^{(r)}_{\max})$, so 
we will denote 
this
subsequence of parameters of {\sc ParamList}${}_{k^{(r)}_{\max}}(i)$ by {\sc ParamList}${}_{k^{(r)}_{\max}}'(i)$.

\begin{algorithm}\label{alg2}{\rm ({\sc TwoStepImprimitive})\\
{\sc Input:}\quad {\sc ParameterList}$(k^{(r)}_{\max})$.

\smallskip\noindent
{\sc Output:}\quad {\sc ParameterListA}$(k^{(r)}_{\max})$ containing those {\sc Lines} 
where $G$ is guaranteed to be 2-step imprimitive relative to $\cC$, and
{\sc ParameterListB}$(k^{(r)}_{\max})$ containing the remaining {\sc Lines} of
{\sc ParameterList}$(k^{(r)}_{\max})$ (sometimes with the information that $\cC$
is minimal or maximal).

\smallskip\noindent
set {\sc ParameterListA}$(k^{(r)}_{\max}):=$ an empty list;\\ 
set {\sc ParameterListB}$(k^{(r)}_{\max}):=$ {\sc ParameterList}$(k^{(r)}_{\max})$;\\
let $n$ be the number of {\sc Lines} of {\sc ParameterListB}$(k^{(r)}_{\max})$;\\
for $i=1,\dots,n$ apply the following tests to {\sc ParamList}${}_{k^{(r)}_{\max}}'(i)=(d,c,x,y,k^{(v)},k^{(r)})$\\
\mm if at least one of the following conditions holds\\
\mm \mm (i)\ $x\leq 4$,\\
\mm \mm (ii)\ $\binom{k}{2} > (x-2)xy+x+y$,\\
\mm \mm (iii)\ there is no $j\leq n$ such that, if {\sc 
ParamList}${}_{k^{(r)}_{\max}}'(j)=(d',c',x',y',k^{(v)'},k^{(r)'})$,\\
\mm \mm \mm then $k'=k$, $c'$ is a proper divisor of $c$, and $d'=cd/c'$,\\
\mm then record that $\cC$ is minimal for {\sc ParamList}${}_{k^{(r)}_{\max}}(i)$;\\
\mm if at least one of the following conditions holds\\
\mm \mm (i)\ $y\leq 4$,\\
\mm \mm (ii)\ $\binom{k}{2} > (y-2)xy+x+y$,\\
\mm \mm (iii)\ there is no $j\leq n$ such that, if {\sc 
ParamList}${}_{k^{(r)}_{\max}}'(j)=(d',c',x',y',k^{(v)'},k^{(r)'})$,\\ 
\mm \mm \mm then $k'=k$, $d'$ is a proper divisor of $d$, and $c'=cd/d'$,\\
\mm then record that $\cC$ is maximal for {\sc ParamList}${}_{k^{(r)}_{\max}}(i)$;\\
\mm if $\cC$ is both minimal and maximal for {\sc ParamList}${}_{k^{(r)}_{\max}}(i)$ \\
\mm \mm then remove {\sc ParamList}${}_{k^{(r)}_{\max}}(i)$ from {\sc ParameterListB}$(k^{(r)}_{\max})$\\
\mm \mm and add it to {\sc ParameterListA}$(k^{(r)}_{\max})$;\\
return {\sc ParameterListA}$(k^{(r)}_{\max})$ and {\sc ParameterListB}$(k^{(r)}_{\max})$.
}
\end{algorithm}

\begin{lemma}\label{lem:test2step}
The output of Algorithm~{\rm\ref{alg2}} is correct.
\end{lemma}

\begin{proof}
If $\cC$ is not minimal then there exists a proper $G$-invariant
refinement $\cC'$ of $\cC$ with $d'$ classes of size $c'$ where $c'$
is a proper divisor of $c$. Similarly if $\cC$ is not maximal then
$\cC$ is a proper refinement of a $G$-invariant partition $\cC'$
with $d'$ classes of size $c'$ where $d'$ is a proper divisor of
$d$. Thus parts (iii) are sufficient conditions for minimality and maximality 
respectively. The sufficiency of parts (i) and (ii) follows from 
Corollary~\ref{cor:maxmin}.
\end{proof}

\bigskip
Algorithm~\ref{alg2} was applied to {\sc ParameterList}$(8)$ (obtained from 
Algorithm~\ref{alg1}). It turned out that the list {\sc ParameterListB}$(8)$ 
returned was empty, and hence that for every {\sc Line} of {\sc ParameterList}$(8)$,
the group $G$ was guaranteed to be 2-step imprimitive relative to $\cC$. 

\begin{theorem}\label{thm:2step}
If the {\sc Hypothesis} holds and $k^{(r)}\leq8$, then the group $G$ is $2$-step
imprimitive relative to $\cC$.
\end{theorem}

\section{Sifting for primitive top and bottom groups}\label{sect:groups}

In this section we give some additional properties of the top group $H=G^\cC$ and the bottom group
$L=G^C$. In Subsection~\ref{sub:topbottom} we present three algorithms for determining possibilities
for these groups in the cases where $\cC$ is maximal or minimal respectively, and determining the
value of the parameter {\sc PossibleGrid}. Applying these algorithms will complete {\sc Step} 3 of
the search procedure described in Section~\ref{sect:outline}, the {\sc Top and Bottom Group Sifts}.

\subsection{Subdegrees}

Assume $G$ is a transitive group acting on a set $X$ and let $\al\in X$.  The orbits of the
stabiliser $G_\al$ in $X$ are called the {\em suborbits of $G$ relative to $\al$} and all
suborbits apart from $\{\al\}$ are called non-trivial (even if they have length 1).  The lengths
of the suborbits are the {\em subdegrees} of $G$, and a subdegree is called non-trivial if it is
the length of a non-trivial suborbit.  The subdegrees are independent of the choice of $\al$
because of the transitivity of $G$. The {\em rank} of the group, abbreviated as $\rk \, G,$ is the
number of suborbits. We have upper bounds on the ranks of $G^\cC$ and $G^C$ in terms of the
parameters of $\cS$.

\smallskip

\begin{lemma}[Delandtsheer, Niemeyer and Praeger~\cite{DelandtsheerNiemeyeretal01}, Proposition 
2.6] \label{lem:subdegree} Assume\\ that the {\sc Hypothesis} holds.
\begin{enumerate}
\item
The number $b^{(r)}$ divides each non-trivial subdegree of $G^\cC$, and in 
particular, $\rk \, G^\cC \le 1 + \delta.$
\item
The number $b^{(r)}$ divides each non-trivial subdegree of $G^C$, and in
particular, $\rk \, G^C \le 1 + \gamma.$
\item Moreover, for $\al\in\cP$, 
$b^{(r)}$ divides each non-trivial subdegree of $G^\cP$
and each orbit length of $G_\al$ in $\{\lambda\in\cL\;:\;\al\in\lambda\}$. 
\end{enumerate}
\end{lemma}

\bigskip
Next we give a link between the transitivity of the bottom group and the 
spectrum. A permutation group $G$ on $X$ is {\it $t$-homogeneous} if it acts transitively on the set of $t$-element 
subsets of $X$. For all $t$, a $t$-transitive group is $t$-homogeneous, while for $t\geq 2$, a transitive $t$-homogeneous 
group is $(t-1)$-transitive (see~\cite[Section 9.4]{DixonMortimer}).

\begin{lemma}\label{lem:twotrans}
Assume that the {\sc Hypothesis} holds.
\begin{enumerate}
\item
If $G^C$ is $3$-homogeneous then $\spec \, \cS = \{1,2\}.$
\item
If $G^C$ is $2$-homogeneous then $\spec \, \cS = \{1,h\}$ with $h \ge 2.$
\end{enumerate}
\end{lemma}

\begin{proof} (i) Consider the linear space $\cS|_C$ defined at (\ref{eqn:induced}), which admits
$G^C$ as an automorphism group. In any linear space, three points are either collinear or they form
a triangle, that is, the three pairs lie on three different lines. Since $G^C$ is transitive on
3-subsets of points of $C,$ only one of these types of three point sets is possible. If all sets of
three points were collinear, then all points of $C$ would be contained in one line of $\cL$, and
hence $c\in\spec\, \cS$, contradicting Theorem~\ref{thm:HM}(ii). Hence no line of $\cL|_C$ contains
more than two points, which means that $\spec \, \cS$ is a subset of $\{1,2\}.$ On the other hand,
by Theorem~\ref{thm:HM}(i), the cardinality of $\spec \, \cS$ is at least $2,$ so $\spec \, \cS =
\{1,2\}.$

(ii) Here, $G^C$ is an automorphism group of $\cS|_C$ which is transitive on unordered pairs of points of $C.$ Since any
two points of $C$ determine a line of $\cL|_C$, and since $G^C$ is transitive on pairs from $C$, all lines of $\cL|_C$
have the same size, say $h$, with $h \ge 2.$ Thus $\spec \, \cS = \{1,h\}$ with $h \ge2.$ \end{proof}

\subsection{Sifting for primitive top and bottom groups}\label{sub:topbottom}

In this subsection we describe separately procedures for restricting the possible top groups $H=G^\cC$ and bottom groups
$L=G^C$ in the cases where $\cC$ is maximal and minimal respectively. Recall that the $G$-invariance of $\cC$ implies
that $G\leq L\wr H$, which is why we refer to $H$ and $L$ as the top and bottom groups of the wreath product. The reason
for separating these procedures is that the procedures are applicable for searches over larger ranges of values of
$k^{(r)}$, and in particular in circumstances where it is not known that both the top group and the bottom group are
primitive. We then present Algorithm~\ref{alg:step3and4} that completes {\sc Step 3} of the overview of the search
described in Section~\ref{sect:outline}.

Our procedures use tables of primitive groups of a given degree, if available. Currently such lists are available in the
computer algebra package GAP \cite{GAP4} for all degrees less than $2,500$. A few tests derived from the theoretical
results presented so far in this paper do not require these lists, but if the relevant list is available, then additional
tests are applied to each primitive group in the list. 

Recall that the parameters and the intersection type are fixed at this stage. First we give the algorithm to sift for
$G^\cC$ (the top group) in the case where $\cC$ is maximal, that is, where $G^\cC$ is a primitive group of degree $d$. By
the \emph{transitivity} of a primitive permutation group $H$ we mean the maximum integer $t$ such that $H$ is
$t$-transitive.

\begin{algorithm}\label{alg:topgroup}{\rm ({\sc TopGroupSift})\\
{\sc Input:}\quad A {\sc Line} from {\sc ParameterList}$(k^{(r)}_{\max})$ for which $\cC$ is
maximal, and the list {\sc Prim}$(d)$ of all primitive permutation groups of degree $d$, if
available. 

\smallskip\noindent
{\sc Output:}\quad {\sc TopGroups}; this is a list of candidate primitive groups $G^\cC$ for that 
{\sc Line} if {\sc Prim}$(d)$ is available, and otherwise contains an entry `{\sc 
Prim}$(d)$ unavailable' (and possibly some other information about $G^\cC$).

\smallskip\noindent
set {\sc TopGroups}$:=$ an empty list;\\
if {\sc Prim}$(d)$ is not available\\
\mm set {\sc TopGroups}[1]\,:= `{\sc Prim}$(d)$ unavailable';\\
\mm if ${\lcm_{i>0} \binom{d}{d_i}} \nmid b$ add to {\sc TopGroups}\ \  `$G^{\cC}\ne \Alt_d,\Sym_d$' 
(Lemma~\ref{lem:altd});\\
else for each $H\in$ {\sc Prim}$(d)$ \\
\mm let $t$ be the transitivity of $H$; \\
\mm if $t > t_{\max}$ skip to next $H$ (Lemma~\ref{lem:maxtrans}); \\
\mm if $b^{(r)}$ does not divide all non-trivial subdegrees of $H$ skip to next $H$ (Lemma~\ref{lem:subdegree}(i)); \\
\mm if $H \ge \Alt_d$ and ${\lcm_{i>0} \binom{d}{d_i}} \nmid b$ skip to next $H$ (Lemma~\ref{lem:altd});\\
\mm add $H$ to {\sc TopGroups};\\
return {\sc TopGroups}.
}
\end{algorithm}

\bigskip

Now we give the algorithm to sift for $G^C$ (the bottom group) in the case where $\cC$ is minimal, that is, where $G^C$ is a
primitive group of degree $c$.

\begin{algorithm}\label{alg:bottomgroup}{\rm ({\sc BottomGroupSift})\\
{\sc Input:}\quad A {\sc Line} from {\sc ParameterList}$(k^{(r)}_{\max})$ for which $\cC$ is minimal, and the list {\sc 
Prim}$(c)$ of all primitive permutation groups of degree $c$, if available. 

\smallskip\noindent
{\sc Output:}\quad {\sc BottomGroups}; this is a list of candidate primitive groups $G^C$ for that
{\sc Line} if {\sc Prim}$(c)$ is available, and otherwise contains an entry `{\sc Prim}$(c)$
unavailable' (and possibly some other information about $G^C$).

\smallskip\noindent
set {\sc BottomGroups}$:=$ an empty list;\\
if {\sc Prim}$(c)$ is not available\\
\mm set {\sc BottomGroups}[1]\,:= `{\sc Prim}$(c)$ unavailable';\\
\mm if $\spec\,\cS\ne\{1,h\}$ with $h\geq2$\\  
\mm \mm add to {\sc BottomGroups}\ \  `$G^C$ is not $2$-homogeneous' (Lemma~\ref{lem:twotrans});\\
\mm if $\spec\,\cS=\{1,h\}$ for some $h\geq3$\\ 
\mm \mm add to {\sc BottomGroups}\ \  `$G^C$ is not $3$-homogeneous' (Lemma~\ref{lem:twotrans});\\
else for each $L\in$ {\sc Prim}$(c)$ \\
\mm if $b^{(r)}$ does not divide all non-trivial subdegrees of $L$ skip to next $L$ (Lemma~\ref{lem:subdegree}(ii)); \\
\mm let $t$ be the transitivity of $L$; \\
\mm if $t \ge 3$ and $\spec \, \cS \neq \{ 1,2\}$ skip to next $L$ (Lemma~\ref{lem:twotrans});\\
\mm if $t \ge 2$ then \\
\mm \mm if $|\spec \, \cS| \neq 2$ skip to next $L$ (Lemma~\ref{lem:twotrans}); \\
\mm \mm if $1 \not\in \spec \, \cS$ skip to next $L$ (Lemma~\ref{lem:twotrans}); \\
\mm add $L$ to {\sc BottomGroups};\\
return {\sc BottomGroups}.
}
\end{algorithm}

To complete {\sc Step} 3 of the search, we apply Algorithms~\ref{alg:topgroup} and
\ref{alg:bottomgroup} to each of the {\sc Lines} of {\sc ParameterListA}$(k^{(r)}_{\max})$, and
remove a {\sc Line} if at least one of the lists {\sc TopGroups} or {\sc BottomGroups} is empty for
that {\sc Line}.  Then we determine the value of {\sc PossibleGrid} for each of the {\sc Lines}
remaining. These tests are recorded in the algorithm below.

\begin{algorithm}\label{alg:step3and4}{\rm ({\sc Top\&BottomGroups}) \\
{\sc Input:}\quad  {\sc ParameterListA}$(k^{(r)}_{\max})$ and lists {\sc Prim}$(n)$ (where 
available) of primitive\\ permutation groups of degree $n$, for all $n$ equal to an entry $d$ or $c$ in some {\sc Line} of \\
{\sc ParameterListA}$(k^{(r)}_{\max})$.

\smallskip\noindent
{\sc Output:}\quad An enhanced  {\sc ParameterListA}$(k^{(r)}_{\max})$ containing lists
{\sc TopGroups} and {\sc BottomGroups}, and the value of {\sc PossibleGrid} 
for each surviving {\sc Line}. 

\smallskip\noindent
for each {\sc Line} of {\sc ParameterListA}$(k^{(r)}_{\max})$\\
\mm compute {\sc TopGroups} (Algorithm~\ref{alg:topgroup});\\
\mm \mm if {\sc TopGroups} is empty then\\
\mm \mm \mm remove this {\sc Line} from {\sc ParameterListA}$(k^{(r)}_{\max})$ and skip to next {\sc Line};\\
\mm \mm else append {\sc TopGroups} to this {\sc Line};\\
\mm compute {\sc BottomGroups} (Algorithm~\ref{alg:bottomgroup});\\
\mm \mm if {\sc BottomGroups} is empty then\\
\mm \mm \mm remove this {\sc Line} from {\sc ParameterListA}$(k^{(r)}_{\max})$ and skip to next {\sc Line};\\
\mm \mm else append {\sc BottomGroups} to this {\sc Line};\\
for each (surviving) {\sc Line} of {\sc ParameterListA}$(k^{(r)}_{\max})$\\
\mm set {\sc PossibleGrid} to be `no';\\
\mm if either $c=d$, or there is another (surviving) {\sc Line} of {\sc ParameterListA}$(k^{(r)}_{\max})$\\
\mm \mm corresponding to the same value of $k$ and a partition with $c$ classes of size $d$, then\\
\mm \mm reset {\sc PossibleGrid} to be `yes';\\
\mm append {\sc PossibleGrid} to this {\sc Line};\\
return {\sc ParameterListA}$(k^{(r)}_{\max})$.
}
\end{algorithm}
  
Note that at this point, in each surviving {\sc Line} of {\sc ParameterListA}$(k^{(r)}_{\max})$, the
entries {\sc TopGroups} and {\sc BottomGroups} are both non-empty.

\section{Quasiprimitive and $G$-normal sifts}\label{sect:stepfour}

This section comprises three subsections which together enable us to complete {\sc Step} 4 of the
search strategy described in the overview in Section~\ref{sect:outline}. In Subsection~\ref{subqp}
we present an algorithm to search for candidate groups $G$ in the cases where $G$ is quasiprimitive.
In Subsection~\ref{sec:normal} we prove several theorems concerning groups in the $G$-normal case,
and in Subsection~\ref{sub:normalsift} we bring together the restrictions these theorems provide
into an algorithm to refine the lists of feasible bottom groups in the case of a $G$-normal
partition. 

\subsection{Sifting for primitive top groups in the quasiprimitive case}
\label{subqp}

In this subsection we give an algorithm that gives more restrictions on the possibilities for top
groups in the case where $G$ is quasiprimitive on $\cP$ and $\cC$ is maximal. The algorithm
produces, for each {\sc Line} of {\sc ParameterListA}$(k^{(r)}_{\max})$ for which {\sc TopGroups}
contains an explicit list of groups, a (possibly empty) list {\sc QuasiprimTopGroups} of feasible
top groups in the case where $G$ is quasiprimitive; in other cases it gives restrictions on such a
list. This is the first part of {\sc Step} 4. Recall that, for a quaisprimitive group $G$, the
kernel $G_{(\cC)}$ of the action of $G$ on $\cC$ is trivial, and therefore $G \simeq G^\cC.$

\begin{algorithm}\label{alg:qptopgroup}{\rm ({\sc QuasiprimitiveSift})\\
{\sc Input:}\quad {\sc ParameterListA}$(k^{(r)}_{\max})$ from Algorithm~\ref{alg:step3and4}.

\smallskip\noindent
{\sc Output:}\quad An enhanced {\sc ParameterListA}$(k^{(r)}_{\max})$ such that each {\sc Line} has
an extra entry {\sc QuasiprimTopGroups} appended: if {\sc TopGroups} has first entry `{\sc
Prim}$(d)$ unavailable', then the list {\sc QuasiprimTopGroups} contains entries `{\sc Prim}$(d)$
unavailable' and `$G^\cC$ almost simple'; otherwise {\sc QuasiprimTopGroups} is a (possibly empty) 
list of candidate quasiprimitive groups $G$ for this {\sc Line}.

\smallskip\noindent
for each {\sc Line} from {\sc ParameterListA}$(k^{(r)}_{\max})$ \\
\mm set {\sc QuasiprimTopGroups}\ $:=$ {\sc TopGroups};\\
\mm if the first entry of {\sc QuasiprimTopGroups} is `{\sc Prim}$(d)$ unavailable' then \\
\mm \mm add to {\sc QuasiprimTopGroups} `$G^\cC$ almost simple' \\
\mm \mm skip to instruction $(\ast)$ ({\sc Case} 1 in Section~\ref{sect:outline}); \\
\mm else for each $H\in $ {\sc TopGroups} \\
\mm \mm if $H$ is not almost simple then remove $H$ from {\sc QuasiprimTopGroups} \\
\mm \mm \mm skip to next $H$ ({\sc Case} 1 in Section~\ref{sect:outline});\\
\mm \mm if $v \nmid |H|$ then remove $H$ from {\sc QuasiprimTopGroups} \\
\mm \mm \mm skip to next $H$ (since $H = G^\cC \simeq G$ is point transitive);\\
\mm $(\ast)$ append {\sc QuasiprimTopGroups} to this {\sc Line};\\
return {\sc ParameterListA}$(k^{(r)}_{\max})$.
}
\end{algorithm}

The result of applying Algorithm~\ref{alg:qptopgroup} in the case $k^{(r)}_{\max}=8$ was a list of
only 19 {\sc Lines} in which {\sc QuasiprimTopGroups} was not an empty list. These {\sc Lines} are
given in Section~\ref{sect:results}, with those parameter {\sc Lines} with $k\leq 8$ distinguished
(for reasons discussed there). See Table~\ref{tab:outputqp} and {\sc Line} 2 of
Table~\ref{tab:outputkleq8}.

\subsection{Further restrictions for the $G$-normal case}
\label{sec:normal}

Here we give some additional restrictions on the group $G$ in the case where $\cC$ is $G$-normal.
The algorithm given at the end of this subsection will complete {\sc Step} 4 described in the
overview in Section~\ref{sect:outline}. First we give several results, the first of which is from
\cite{PraegerTuan02}.  Recall that $d_1$ is the number of classes that meet a line in exactly one
point. A permutation group is {\it semiregular} if only the identity element fixes a point; it is
{\it regular} if it is both transitive and semiregular.

\begin{theorem}[Praeger and Tuan~\cite{PraegerTuan02}, Theorems 1.5 and 1.6] \label{thm:semiregular} Assume 
that the {\sc Hypothesis} holds and that $\cC$ is $G$-normal. 
\begin{enumerate}
\item[(a)] If $k > 2x + \frac{3}{2} +
\sqrt{4x - \frac{7}{4}},$ then $G_{(\cC)}$ is semiregular on points and lines, $|G_{(\cC)}| = c$ is
odd, and $d_1 > 0.$
\item[(b)] If $x\leq 8$, then $G_{(\cC)}$ has an abelian subgroup $S$ of index at most $2$ such that $S$ is normal in $G$, semiregular on points, and $|S|=c$ is odd.
\item[(c)] If either of the conditions of (a) or (b) holds, and if $\cC$ is minimal, then  $c$ is an odd prime power and $G^C$ is affine.
\end{enumerate} \end{theorem}

\begin{proof} Parts (a) and (b) follow from Theorems 1.5 and 1.6 of 
\cite{PraegerTuan02} respectively. Suppose now that $\cC$ is minimal and in part (a) set $S:=G_{(\cC)}$. Then $S^C\cong S$ by \cite[Theorem 1]{CaminaPraeger93}, and hence $S^C$ is an odd order semiregular normal subgroup of the
primitive group $G^C$. Thus  $S^C$ is elementary abelian,
$c$ is a prime power and $G^C$ is affine.
\end{proof}

Recall the definition in (\ref{eqn:induced}) of the linear space $\cS|_F= (F,\cL|_F)$ induced on a subset $F\subseteq
\cP$ of size at least two. Under certain conditions, $G_{(\cC),\alpha}$ fixes exactly one point of each class of the
partition $\cC$ and, for $F:= \Fix_\cP\,(G_{(\cC),\alpha})$, the induced linear space $\cS|_F$ has constant line size.
The next result extends a result of Praeger and Tuan~\cite{PraegerTuan02}.

The proof of part (a) uses the fact that, for a line-transitive, point-imprimitive group $G$, every
involution in $G$ fixes at least one point.  For if an involution in a line-transitive group $G$ has
no fixed points then it was shown in \cite[Lemma 4]{CaminaSiemons89a} that $k$ divides $v$, and then
by \cite{CaminaGagen84}, $G$ is flag-transitive and hence primitive on $\cP$. The proof also uses
the fact that a group of odd order is soluble. Recall that a point $\al$ and a line $\lambda$ are
called $i$-incident if $\lambda$ meets the class of $\cC$ containing $\al$ in exactly $i$ points;
also $C\in\cC$ is said to be $i$-incident with $\lambda$ if $|\lambda\cap C|=i$.

\begin{theorem}\label{thm:normal}
Assume that the {\sc Hypothesis} holds and that $\cC$ is  
$G$-normal and minimal. Let $\al\in C\in\cC$, and set 
$F:=\Fix_\cP(G_{(\cC),\al})$. Then 
\begin{enumerate}
\item[(a)] either
\begin{enumerate}
\item[(i)] $G_{(\cC)}=Z_p^a$ is semiregular of order $c=p^a$
for some odd prime $p$; or
\item[(ii)] $F\cap C=\{\al\}$, and the set of classes of $\calC$ that
contain some point of $F$ is a block of imprimitivity for
$G^\calC$.
In particular, if in addition $\calC$ is maximal, then either
\begin{enumerate}
\item[(ii-1)] $F=\{\al\}$, and $d_1\leq 1$; or
\item[(ii-2)] $F$ consists of exactly one point from each class
of $\calC$.
\end{enumerate}
\end{enumerate}
\item[(b)] In particular, if $k \ge 2x$, 
then $c = p^a$ for some odd prime $p$, and either 
\begin{enumerate}
\item[(i)] $G_{(\cC)}=Z_p^a$ is 
semiregular, or
\item[(ii)] $G_{(\cC)} = Z_p^a \cdot Z_2,$ $\spec\,\cS = \{1,2\},$
(a)(ii-2) holds for $F$,
$\cS|_F = (F,\cL|_F)$ is a linear space with lines of size $d_1$
admitting a line-transitive action by $N_G(G_{(\cC),\alpha})$. 
Moreover $d_1 = k - 2x \ge 2,$ $y = \binom{d_1}{2},$
$y$ divides $\binom{d}{2}$, and $d_1 -1 $ divides $d - 1.$
Also, for any pair $C,D$ of distinct classes of $\cC,$
$G_{C,D}$ fixes setwise disjoint subsets of
$\cC \setminus \{ C,D\}$ of sizes $d_1 - 2$ and $x.$
In particular if $d_1 = 2$ then $G^\cC$ is $2$-homogeneous.
\end{enumerate}
\end{enumerate}
\end{theorem}

\begin{proof} Part (b) was proved in \cite[Theorem 1.4]{PraegerTuan02}. Thus we only need to prove part (a). Let $N:=
G_{(\cC)}$, so $N$ is a non-trivial normal subgroup of $G$ and $\cC$ is the set of $N$-orbits in $\cP$. By the minimality
of $\cC$, the permutation group $G_C^C$ is primitive on $C$, and $N^C$ is a transitive normal subgroup.  Moreover, by
Theorem~\ref{thm:CP93}, $N$ is faithful on $C$.  Let $\al\in C$. Consider first the case where $N^C$ is regular.  Then
$N_\al=1$, and so $N$ is semiregular on $\cP$.  As discussed above, all involutions in $G$ fix at least one point, and it
follows that $|N|$ is odd. Hence $N$ is soluble.  Since $N^C\cong N$ is a soluble regular normal subgroup of the
primitive group $G_C^C$, it follows that $G_C^C$ is of affine type and $N^C$ is elementary abelian of order $|C|=p^a$ for
some odd prime $p$.  Thus (a)~(i) holds.

Suppose now that $N^C$ is not regular.  Since $N^C$ is a normal subgroup of the transitive group$G_C^C$, $\Fix_C(N_\al)$
is a block of imprimitivity for $G_C^C$, and since $G_C^C$ is primitive it follows that $\Fix_C(N_\al)=\{\al\}$. 
Therefore $F=\Fix_\cP(N_\al)$ consists of at most one point from each class of $\cC$.  If $\beta\in F$ then
$N_\al\subseteq N_\beta$, and since $|N_\al|=|N|/c=|N_\beta|$ we have $N_\al=N_\beta$. Moreover, if $g\in G$ is such that
$\al^g=\beta$, then $(N_\al)^g=N_{\al^g}=N_\beta=N_\al$, that is $g\in N_G(N_\al)$.  It follows that $M:=N_G(N_\al)$ is
transitive on $\Fix_\calP(N_\al)$.  Since $G_\al\leq M<G$, we deduce that $\Fix_\cP(N_\al)$ is a block of imprimitivity
for $G$ in $\cP$ and its setwise stabiliser is $M$.

Let $\mathcal{D}$ be the subset of $\cC$ consisting of those classes containing a point of $F$, and let
$D:=\bigcup_{C'\in\mathcal{D}}C'$.  Since $N$ fixes each class of $\cC$ setwise, $D$ is $N$-invariant, and since $F$ is
$M$-invariant, $D$ is also $M$-invariant.  Thus $NM$ leaves $D$ invariant, and since $N$ is transitive on each class of
$\cC$ and $M$ is transitive on $F$, it follows that $NM$ is transitive on $D$.  Again, since $G_\al<NM\leq G$ it follows
that $D$ is also a block of imprimitivity for $G$ in $\cP$.  It then follows easily that $D$ is a block of imprimitivity
for the action of $G$ on $\cC$. If $\cC$ is maximal then $G$ is primitive on $\cC$, so either $D=C$ or $D=\cP$, and hence
either $F=\{\al\}$ or (a) (ii-2) holds respectively. Finally if $F=\{\al\}$, we claim that $d_1\leq1$.  If $d_1\geq2$ and
$\al,\beta$ are 1-incident with a line $\lambda$, then $N_\lambda$ fixes both $\al$ and $\beta$.  However, by
\cite[Corollary 4.2]{PraegerTuan02}, $N_\lambda=N_\al=N_\beta$, and so $F$ contains both $\al$ and $\beta$ contrary to
our assumption.  \end{proof}

\bigskip

The following result explores case (a)~(ii), that is, when the kernel $G_{(\cC)}$ is not semiregular
on points. It makes use of various results from~\cite{BettenDelandtsheeretal02}.

\begin{theorem}\label{notL27} 
Assume that the {\sc Hypothesis} holds and that $\cC$ is  
$G$-normal and minimal  with $G_{(\cC)}$ not 
semiregular. Let  $m\in \spec\,\cS$ with $m\geq2$.
\begin{enumerate}
\item[(a)] Then the
$G_{(\cC)}$-orbits on lines have equal length $c/z$ where $z$
divides $\gcd(c,m)$, and each prime divisor of $|G_{(\cC),\al}|$ is
at most $m(m-1)/z$. 
\item[(b)] If $m=2$, then $G_{(\cC),\al}$ is a Sylow $2$-subgroup of
$G_{(\cC)}$, and $\Fix_\calP(G_{(\cC),\al})$ consists of exactly one
point from each class of $\calC$. Moreover either
\begin{enumerate}
\item[(i)] $\Soc(G_{(\cC)})=Z_p^a$ (so $G^C$ is affine) and $c=p^a$ for 
some odd prime $p$ and $a\geq1$, or
\item[(ii)] $\Soc(G_{(\cC)})=\PSL(2,q)^a$, where $a\geq 2$ and
$(c,q)=(21^a, 7)$ or $(45^a,9)$.
\end{enumerate}
\item[(c)] If $m=3$ or $4$, then
$G_{(\cC),\al}$ is a $\{2,3\}$-group, and $G_{(\cC),\al}$ has an
orbit in $C$ of length $\ell$, where $\ell>1$ and $\ell$ divides $m(m-1)$.
\end{enumerate}
\end{theorem}

\begin{proof} Set $N:=G_{(\cC)}$.
By Theorem~\ref{thm:normal}, $\Fix_C(N_\al)=\{\al\}$.  Suppose that $\beta\in C\setminus\{\al\}$ and 
that the (unique) line $\lambda$ containing $\al$ and $\beta$ is $m$-incident with $C$.  By 
\cite[Proposition 4.1]{DelandtsheerNiemeyeretal01}, the $N$-orbits on lines have constant length 
$|N:N_\lambda|= c/z$ for some divisor $z$ of $c$. Since $|\lambda\cap C|=m$, we have, for each 
$\gamma\in \lambda\cap C$, that $m_\gamma:=|N_\lambda: N_{\lambda\,\gamma}|\leq m$.  Thus 
$c=|N:N_\gamma|$ divides $m_\gamma c/z$, and so $z$ divides $m_\gamma$. Since this is true for all 
such $\gamma$, it follows that $z$ divides $m$. Now $N_{\al\,\beta}$ fixes $\lambda$, and hence 
$N_{\al\,\beta} \leq N_{\lambda\,\al}\leq N_\al$. Thus the length $\ell$ of the $N_\al$-orbit 
containing $\beta$ is $\ell=|N_\al:N_{\al\,\beta}|=\ell_1\ell_2$, where 
$\ell_1=|N_\al:N_{\lambda\,\al}|$ and $\ell_2= |N_{\lambda\,\al}:N_{\al\,\beta}|$.  Now $\ell_2$ is 
the length of the $N_{\lambda\,\al}$-orbit containing $\beta$.  This is at most the length $m_\beta$ 
of the $N_\lambda$-orbit containing $\beta$ which, in turn, is contained in $(\lambda\cap 
C)\setminus\{\al\}$, whence $\ell_2\leq m_\beta\leq m-1$. Also
\[
\ell_1=|N_\al:N_{\lambda\,\al}|=\frac{|N:N_\lambda|\,|N_\lambda:N_{\lambda\,\al}|}{|N:N_\al|}
=\frac{1}{z}\, |N_\lambda:  N_{\lambda\,\al}|\leq\frac{m}{z}
\]
and hence $\ell=\ell_1\ell_2\leq m(m-1)/z$. By Theorem~\ref{thm:CP93}, $N$ is faithful on $C$, and 
hence $N_\al^C\cong N_\al\ne 1$. Thus $N^C$ is a non-regular normal subgroup of the primitive group 
$G_C^C$.  By \cite[Corollary~2.2\,(c)]{BettenDelandtsheeretal02}, each prime $p$ dividing $|N_\al|$ 
also divides the order of the permutation group induced by $N_\al$ on any of its orbits in
$C\setminus\{\al\}$.  Therefore $p$ is at most the minimum of the lengths of the $N_\al$-orbits in 
$C\setminus\{\al\}$, and in particular $p\leq m(m-1)/z$.  This completes the proof of part (a). We 
now use this strategy a little more carefully to prove the other parts.

Suppose that $m=2, 3$ or 4.  Since $N_{\al\,\beta} \leq N_{\lambda\,\al}\leq N_\al$, the permutation group induced by
$N_\al$ on its orbit containing $\beta$ is a subgroup of $S_{\ell_2}\wr S_{\ell_1}$, and since $\ell_1\leq m/z\leq 4$ and
$\ell_2\leq m-1\leq 3$, it follows that the only primes dividing the order of this induced permutation group are 2 or 3. 
It follows from \cite[Corollary~2.2\,(c)]{BettenDelandtsheeretal02} that $N_\al$ is a $\{2,3\}$-group. Since
$\Fix_C(N_\al)=\{\al\}$, all $N_\al$-orbits in $C\setminus\{\al\}$ have length at least 2, and the $N_\al$-orbit
containing $\beta$ has length $\ell=\ell_1\ell_2\leq m(m-1)/z$. If $m=2$, then we must have $\ell=2$. In this case, by
\cite[Theorem~2.3\,(c)]{BettenDelandtsheeretal02}, $N_\al$ is a Sylow 2-subgroup of $N$ and (since $N\cong N^C$) one of
(b)~(i) or (b)~(ii) holds (but possibly with $a=1$ in the case of (b)~(ii)). By Sylow's Theorems, the stabiliser in $N$
of an arbitrary point of $\cP$ is conjugate to $N_\al$, and hence $N_\al$ fixes a point from each class of $\cC$. Thus
(b) is proved, except for proving that $a\geq 2$ in case (b)~(ii) (see below for the proof).

Suppose next that $m=3$.  Then $\ell_1\leq 3$, $\ell_2\leq 2$, and $\ell=\ell_1\ell_2\leq 6/z$.  Thus either $\ell$
divides 6 or $\ell_1=\ell_2=2$ and $z=1$.  However in the latter case,
$|N_\lambda:N_{\lambda\,\al}|=|N_\al:N_{\lambda\,\al}|=\ell_1=2$, and since $|\lambda\cap C|=m=3$ it follows that
$N_\lambda$ fixes a point of $(\lambda\cap C)\setminus\{\al\}$, and hence that
$\ell_2=|N_{\lambda\,\al}:N_{\al\,\beta}|\leq m-2$, contradicting the fact that $\ell_2=2$. Now assume that $m=4$.  Here
$\ell_1\leq 4/z, \ell_2\leq \min\{3, m_\beta\}$ and $\ell=\ell_1\ell_2\leq 12/z$.  Hence $\ell$ divides 12, or
$\ell_1=4,\,\ell_2=2,\,z=1$, or $\ell_1=\ell_2=3,\,z=1$. In the third case, we obtain a contradiction using the same
argument as for $m=3$.  Thus we may assume that $\ell_1=4,\,\ell_2=2,\,z=1$. In this case $|N_\lambda:N_{\lambda\,\al}|=
|N_\al:N_{\lambda\,\al}|=\ell_1=4$ so $N_\lambda$ is transitive on $\lambda\cap C$.  Since
$|N_{\lambda\,\al}:N_{\al\,\beta}|= \ell_2=2$ it follows that $N_{\lambda\,\al}$ fixes two points of $\lambda\cap C$ and
interchanges the remaining two points.  Replacing $\beta$ with the point in $(\lambda\cap C)\setminus\{\al\}$ fixed by
$N_{\lambda\,\al}$, we obtain a new $N_\al$-orbit of length equal to the new value of $\ell$, namely 4. This completes
the proof of part (c).

It remains to prove that $a\geq2$ in part(b)~(ii).  Suppose then that (b)~(ii) holds with $a=1$.  We showed above that
the $N_\al$-orbit $\Delta:=\{\beta^g\,|\,g\in N_\al\}$ has length 2. Let $\lamo$ be the (unique) line containing
$\Delta$. Then $N_\al$ fixes $\lamo$ setwise, and hence $\lamo$ is a union of some $N_\al$-orbits.  Let
$\mathcal{D}:=\{\Delta^g\,|\,g\in G\}$ and $\mathcal{D}_C:=\{\Delta^g\,|\,g\in G_C\}$. By \cite[Lemma
3.2\,(c)]{BettenDelandtsheeretal02}, $|\mathcal{D}_C|=2c$ and it follows that $|\mathcal{D}|=2cd=2v$. Now the setwise
stabiliser of $\De$ satisfies $G_\De\leq G_{\lamo}<G$, and since $G$ is line-transitive it follows that each line
contains exactly $s:=|G_{\lamo}:G_\De|$ elements of $\mathcal{D}$, and that the stabiliser of a line acts transitively on
the $s$ elements of $\mathcal{D}$ it contains. This implies that $bs=|\mathcal{D}|=2v\leq 2b$ (by Fisher's inequality
(\ref{eqn:fisher})), so $s\leq 2$. Substituting $bs=2v$ into (\ref{eqn:b}) yields $v-1=2k(k-1)/s$.

Since $N_\al$ fixes a unique point in each class of $\cC$, it follows from \cite[Lemma 3.1]{BettenDelandtsheeretal02}
that $N_\al$ has $2d$ orbits of length $2$ (all members of $\calD$), $2d$ orbits of length $4$ (each the union of two
disjoint members of $\calD$), and either $d$ (if $q=7$) or $4d$ (if $q=9$) orbits of length $8$. Moreover, if $q=9$ then
exactly $2d$ of the $N_\al$-orbits of length $8$ are unions of pairwise disjoint members of $\calD$. Also, the union $D$
of a fixed point, say $\gamma$, of $N_\al$ and an $N_\al$-orbit of length 2 lying in the class $C(\gamma)$ is a block of
imprimitivity for $N$ in $C(\gamma)$, and $D$ contains three members of $\calD$. Thus, since $s\leq2$, $\lamo$ contains
no $N_\al$-orbits of length 4, and in the case $q=9$, $\lamo$ contains none of the $N_\al$-orbits of length 8 which are
unions of elements of $\mathcal{D}$.  Also $\lamo$ does not contain an $N$-block of imprimitivity of length 3 in any
class of $\cC$, and in particular $\al\not\in\lamo$. If $s=2$ then the second element $\De'$ of $\mathcal{D}$ contained
in $\lamo$ is also fixed setwise by $N_\al$. Since $N_\al$ fixes exactly one point in each class of $\cC$ and since
$\De'$ is contained in some class, it follows that $\De'$ is an $N_\al$-orbit of length 2. 

Since $M:=N_G(N_\al)$ is transitive on $F:=\Fix_\cP(N_\al)$, it follows that $M$ is transitive on $\cC$, and hence $M$
permutes the $N_\al$-orbits of length $8$ in orbits of length $d$ (or possibly $2d$ if $q=9$).  Now $G_\al\leq M$ and $M$
is transitive on the $d$ fixed points of $N_\al$, and so $|G:M|=c$.  Also $M$ is transitive on the set of $2d$ elements
of $\mathcal{D}$ which are $N_\al$-orbits.  Thus the $M$-orbit $\lamo^M$ containing $\lamo$ consists of $2d/s$ lines. Let
$\lamo$ contain $a_i$ of the $N_\al$-orbits of length $i$, for $i=1$ and $i=8$. Then $k=a_1+2s+8a_8$. Since at most $2d$
orbits of $N_\al$ of length 8 (or at most $d$ if $q=7$) lie in lines in $\lamo^M$, and since each such orbit is contained
in at most one line, it follows that 
and we deduce that either (i) $a_8=0$, or (ii) $a_8=1$, $s=2$, or (iii) $q=9$, and $a_8=s\leq 2$.  In particular
$a_8|\lamo^M|\leq 2d$ and hence $a_8\leq s$, so $k\leq a_1+10s$.

Suppose first that $\lamo^M$ contains all the lines fixed by $N_\al$.  Any line $\lam'$ containing a pair of points of 
$F=\Fix_\cP(N_\al)$ is fixed by $N_\al$, and by assumption $\lam'\in\lamo^M$. It follows that
$\lamo$ must contain at least one pair of points from $F$, so $a_1\geq2$.  There are $d(d-1)$ such ordered pairs of
points, and each of the $2d/s$ lines of $\lamo^M$ contains $a_1(a_1-1)$ of them.
Therefore $a_1(a_1-1)=(d-1)s/2$. We showed above that $v-1=2k(k-1)/s$, and so
\begin{eqnarray*}
21d-1\leq v-1&=&\frac{2}{s} k(k-1)\leq \frac{2}{s}(a_1+10s)(a_1+10s-1)\\
&=&\frac{2}{s}a_1(a_1-1) + 20(2a_1-1)+200s\\
&=&d-1 + 40a_1 -20 + 200s.
\end{eqnarray*}
Hence $d\leq 2a_1+10s-1$, so $a_1^2-a_1=(d-1)s/2\leq s(a_1+5s-1)$.  This
implies that $a_1\leq 3$ if $s=1$ and $a_1\leq 6$ if $s=2$.  For $s=1$, we compute the possibilities for $a_1=2$ or 3, $d=2a_1(a_1-1)+1, k=a_1+2+8a_8$
(with $a_8=0$ or 1), and $c=(2k(k-1)+1)/d$.  In no case do we find $c=21$ or
45, so we have a contradiction.  Hence $s=2$ and so $d=a_1(a_1-1)+1,
k=a_1+4+8a_8$ and $c=(k(k-1)+1)/d$.  We compute all possibilites for these
parameters with $a_1=2,\dots,6$ and $a_8=0, 1, 2$. The only case for which
$c$ turns out to be 21 or 45 is $(a_1,a_8,c)=(6,2,21)$. However we showed 
above that $a_8=2$ is only possible if $q=9$ and $c=45$. Thus we have a 
contradiction.

Therefore there is a line $\lam'$ fixed by $N_\al$ and not lying in $\lamo^M$, and hence containing no $N_\al$-orbit of
length 2.  The $s$ (at most 2) elements of $\mathcal{D}$ contained in $\lam'$ are fixed setwise by $N_\al$, and therefore
$s=2$ and $N_\al$ interchanges the two elements of $\mathcal{D}$ in $\lam'$.  Thus $\lam'$ contains a unique
$N_\al$-orbit of length 4.  By \cite[Lemma 3.1\,(a)]{BettenDelandtsheeretal02}, there are $2d$ such orbits and they are
permuted transitively by $M$, and hence $(\lam')^M$ forms a second $M$-orbit of lines fixed by $N_\al$, and
$|(\lam')^M|=2d$, while $|\lamo^M|=2d/s=d$.  Since any line fixed by $N_\al$ and not lying in $\lamo^M$ contains an
$N_\al$-orbit of length 4, these two $M$-orbits contain all of the lines fixed by $N_\al$. Suppose that $\lam'$ contains
$e_i$ of the $N_\al$-orbits of length $i$, for $i=1, 8$.  Then $k=e_1+4+8e_8$.  Each line containing two points of $F$ is
fixed by $N_\al$ and so lies in $\lamo^M$ or in $(\lam')^M$.  Thus $d(d-1)=da_1(a_1-1)+2de_1(e_1-1)$, that is,
\[ d-1=a_1(a_1-1)+2e_1(e_1-1). \]  
Also $k=a_1+4+8a_8=e_1+4+8e_8$. If $\lam'\supseteq F$, then $\lam'$ would be fixed setwise by $M$, which is not the case. 
Hence $\lam'\cap F\ne F$ so $e_1\leq d-1$. Similarly $a_1\leq d-1$. If $e_8=0$ then, using $v-1=2k(k-1)/s=k(k-1)$, we 
obtain 
\begin{eqnarray*}
21d-1\leq v-1&=& k(k-1)= (e_1+4)(e_1+3)\\
&=&e_1(e_1-1) + 8e_1+12\\
&\leq&\frac{d-1}{2} + 8(d-1)+12
\end{eqnarray*}
which is a contradiction. Hence $e_8\geq 1$, which implies that there are $2de_8$ orbits of $N_\al$ of length 8 contained
in lines in $(\lam')^M$. If $a_8=0$, then a similar argument leads to a contradiction. Hence $a_8\geq 1$, which implies
that there are $da_8$ orbits of $N_\al$ of length 8 contained in lines in $\lamo^M$. By \cite[Lemma
3.1\,(a)]{BettenDelandtsheeretal02}, there are $d$, $4d$ orbits of $N_\al$ of length 8, for $q=7,9$ respectively. Hence
$q=9$, $e_8=1$, and $1\leq a_8\leq 2$. However, this means that there are at least $3d$ of the $N_\al$-orbits of length 8
contained in lines in $(\lam')^M$ or $\lamo^M$, and by \cite[Lemma 3.1\,(c)]{BettenDelandtsheeretal02}, some of these are
unions of elements of $\mathcal{D}$. This contradicts the fact that each line in $(\lam')^M$ or $\lamo^M$ contains only
two elements of $\mathcal{D}$. \end{proof}

Our last result in this subsection is applied in particular in $G$-normal situations when the bottom group is 
2-transitive.

\begin{theorem}\label{not28}
Assume that the {\sc Hypothesis} holds and that $\cC$ is  
$G$-normal and minimal such that $L:=G^C$ is almost simple
and $S=\Soc(L)$ has at most two conjugacy classes of subgroups
isomorphic to $S_\alpha$ (where $\alpha\in C$). 
(In particular this holds if $L$ is $2$-transitive
and almost simple.) Then  there is a second $G$-invariant
partition $\calC'$ of $\calP$ with $c$ classes of size $d$ such that
$S\cong\Soc(G^{\calC'})$. Moreover, each class of $\cC'$ meets each class
of $\cC$ in $1$ point. 
\end{theorem}

\begin{proof} Since the socle of an almost simple 2-transitive permutation group of degree $c$ has at most two pairwise
inequivalent transitive representations of degree $c$, the condition on $\Soc(L)$ holds if $L$ is almost simple and
2-transitive. (This follows from the classification of finite 2-transitive groups, see~\cite{Cam}.)

Since $\calC$ is $G$-normal, $G_{(\cC)}\ne1$. Let $S:=\Soc(G_{(\cC)})$, so $S$ is normal in $G$. Since $\cC$ is minimal,
the classes of $\cC$ are the orbits of $S$ in $\cP$.  Hence by Theorem~\ref{thm:CP93}, $S$ is faithful on $C$. Since $L$
is almost simple and $S^C$ is normal in $L$, it follows that $S^C=\Soc(L)$ and so $S\cong S^C$ is simple. Let $\al\in\cP$
with $\al\in C$, and let $F$ denote the set of fixed points in $\cP$ of the stabilizer $S_\al$.  Since $S^C$ is normal in
the primitive group $L$, $F\cap C$ is a block of imprimitivity for $L$ in $C$, and so $F\cap C$ is either $C$ or
$\{\al\}$. However, since $L$ is primitive and almost simple, its socle $S^C$ is not regular, and hence $F\cap
C=\{\al\}$. Moreover, since $S$ is normal in $G$, $F$ is a block of imprimitivity for $G$ in $\cP$, and so determines a
non-trivial $G$-invariant partition, $\cC'$ say.

We claim that $F$ consists of exactly one point from each class of $\cC$. Define a relation on the
classes of $\cC$ by $C_1 \sim C_2$ if, for $\al_1\in C_1$, $S_{\al_1}$ fixes a point in $C_2$.
Clearly this is an equivalence relation and there are at most two equivalence classes by the
assumption on $\Soc(L)$.  Moreover, it follows, from the definition of $\sim$, that $\sim$ is
$G$-invariant.  Suppose that there are two $\sim$-equivalence classes, and consider the two-class
partition of $\cP$ where each class is the union of the $\cC$-classes in a $\sim$-equivalence
class. This is a $G$-invariant point partition with two classes, contradicting
Theorem~\ref{thm:CNP}~(ii).  Thus there is only one $\sim$-equivalence class, and hence $F$
contains at least one point from each class of $\cC$. The argument in the previous paragraph shows
conversely that $|F\cap C| \leq1$ for each class $C$, and hence $F$ consists of exactly one point
from each class of $\cC$, proving the claim. Finally, $L$ acts on $\cC'$ in the same way that it
acts on $C$. Then, as $S\triangleleft\, G$, we must have that $G^{\cC'}$ is almost simple with
socle $S^{\cC'}\cong S$. \end{proof}

\bigskip

\subsection{Sifting for primitive bottom groups in the $G$-normal
case}\label{sub:normalsift}

In this subsection we describe an algorithm for restricting the possibilities for the bottom group
$L=G^C$ in the case where the partition $\cC$ is $G$-normal and minimal. This will complete {\sc
Step} 4 of the overview in Section~\ref{sect:outline}.

\begin{algorithm}\label{alg:gnormalbottomgroup}{\rm ({\sc GNormalSift})\\
{\sc Input:}\quad {\sc ParameterListA}$(k^{(r)}_{\max})$ from Algorithm~\ref{alg:step3and4}.

\smallskip\noindent
{\sc Output:}\quad An enhanced {\sc ParameterListA}$(k^{(r)}_{\max})$ such that each {\sc Line} has
an extra entry {\sc GNormalBottomGroups} appended: the algorithm may prove that $\cC$ is not
$G$-normal and in this case {\sc GNormalBottomGroups} is an empty list; if this is not the case then
if {\sc BottomGroups} has first entry `{\sc Prim}$(c)$ unavailable', then {\sc GNormalBottomGroups}
contains `{\sc Prim}$(c)$ unavailable' and possibly some other information about the bottom groups
in the case when $\cC$ is $G$-normal; otherwise {\sc GNormalBottomGroups} is a (possibly empty) list
of candidate bottom groups for this {\sc Line} in the case when $\cC$ is $G$-normal.

\smallskip\noindent
for each {\sc Line} from {\sc ParameterListA}$(k^{(r)}_{\max})$ \\
\mm set {\sc GNormalBottomGroups}$:=$ {\sc BottomGroups};\\
\mm if $k \ge 2x$ and $c$ is not an odd prime power then \\
\mm \mm reset {\sc GNormalBottomGroups}$:=$ an empty list \\
\mm \mm skip to instruction $(\ast)$ (Theorem~\ref{thm:normal}~(b));\\
\mm if the first entry of {\sc GNormalBottomGroups} is `{\sc Prim}$(c)$ unavailable' then\\
\mm \mm if $k > 2x + \frac{3}{2}+\sqrt{4x - \frac{7}{4}}$,\quad or if $x\leq 8$\\
\mm \mm \mm add to {\sc GNormalBottomGroups}\ `$L$ is affine' 
(Theorem~\ref{thm:semiregular});\\ 
\mm \mm if $2\in\spec\,\cS$ \\
\mm \mm \mm if $c=c_0^a$ with $a\geq2$ and $c_0=45$ or $21$ \\
\mm \mm \mm \mm add to {\sc GNormalBottomGroups}\ `$\Soc(L)={\rm PSL}(2,q)^a$';\\
\mm \mm \mm or if $c$ is an odd prime power \\
\mm \mm \mm \mm add to {\sc GNormalBottomGroups}\ `$L$ is affine';  \\
\mm \mm \mm else reset {\sc GNormalBottomGroups}$:=$ an empty list \\
\mm \mm \mm \mm skip to instruction $(\ast)$ (Theorems~\ref{thm:normal}~(a) and~\ref{notL27}~(b));\\
\mm \mm if $3$ or $4\in\spec\,\cS$ \\
\mm \mm \mm add to {\sc GNormalBottomGroups} `$\Soc(L)_\alpha$ is a $\{2,3\}$-group' (Theorem~\ref{notL27}~(c));\\
\mm else \\
\mm \mm for each $L\in ${\sc BottomGroups} \\
\mm \mm \mm if $L$ is not affine \\
\mm \mm \mm \mm if $k > 2x + \frac{3}{2}+\sqrt{4x - \frac{7}{4}}$,
\quad or if $x\leq 8$\\
\mm \mm \mm \mm \mm remove $L$ from {\sc GNormalBottomGroups}\\
\mm \mm \mm \mm \mm  skip to next $L$ (Theorem~\ref{thm:semiregular});\\ 
\mm \mm \mm \mm if $2\in\spec\,\cS$ \\
\mm \mm \mm \mm \mm if $c$ is not $c_0^a$ with $a\geq 2$ and $c_0=45$ or $21$ \\
\mm \mm \mm \mm \mm \mm remove $L$ from {\sc GNormalBottomGroups} \\
\mm \mm \mm \mm \mm \mm skip to next $L$ (Theorem~\ref{notL27}~(b)); \\
\mm \mm \mm \mm if $3$ or $4\in\spec\,\cS$ \\
\mm \mm \mm \mm \mm if the point stabilizer $\Soc(L)_\alpha$ is not a $\{2,3\}$-group \\
\mm \mm \mm \mm \mm \mm remove $L$ from {\sc GNormalBottomGroups} \\
\mm \mm \mm \mm \mm \mm skip to next $L$ (Theorem~\ref{notL27}~(c), since $\Soc(L)_\alpha\leq G^C_{(\cC),\alpha}$);\\
\mm \mm \mm if $L$ is almost simple \\
\mm \mm \mm \mm if {\sc PossibleGrid} is `no' \\
\mm \mm \mm \mm \mm if $\Soc(L)$ has $\leq 2$ conjugacy classes of subgroups isomorphic to $\Soc(L)_\alpha$ \\
\mm \mm \mm \mm \mm \mm remove $L$ from {\sc GNormalBottomGroups} \\
\mm \mm \mm \mm \mm \mm skip to next $L$ (Theorem~$\ref{not28}$);\\
\mm $(\ast)$ append {\sc GNormalBottomGroups} to this {\sc Line};\\
return {\sc ParameterListA}$(k^{(r)}_{\max})$.
}
\end{algorithm}

\smallskip

The result of applying Algorithm~\ref{alg:gnormalbottomgroup} in the case $k^{(r)}_{\max}=8$ was a list
of 36 {\sc Lines} in which {\sc GNormalBottomGroups} was not an empty list. These {\sc Lines} are
given in Section~\ref{sect:results}, where those {\sc Lines} with $k\leq 8$ are distinguished (for
reasons discussed there). See Table~\ref{tab:outputgnormal} and Table~\ref{tab:outputkleq8}.

\section{Results of algorithms applied for $k^{(r)}\leq 8$}\label{sect:results}

In this section we present summaries of the output of
Algorithms~\ref{alg:step3and4},~\ref{alg:qptopgroup} and \ref{alg:gnormalbottomgroup} in the case
$k^{(r)}_{\max}=8$. In \cite{CaminaMischke96}, Camina and Mischke classified all line-transitive
and point-imprimitive linear spaces with $k\leq 8$, extending the classification of
\cite{NNOPP},\cite{KPP93}. We decided to include these cases in the parameter sift so that it could
provide a check of our sifting process. In light of their result, we present the surviving {\sc
Lines} from {\sc ParameterListA}$(8)$ in three separate tables: (i) those {\sc Lines} with $k\leq
8$; (ii) those {\sc Lines} with $k>8$ and which are potentially quasiprimitive; and (iii) those
{\sc Lines} with $k>8$ and which are potentially $G$-normal.

The entry in the column labeled {\sc Line} in each of these tables is the line number from {\sc
ParameterList}$(8)$, which is given in increasing order of $v$. The column labeled ``int type''
records the intersection type $(0^{d_0},1^{d_1},\ldots, k^{d_k})$ by omitting the entry $0^{d_0}$
and any entries $i^{d_i}$ for which $d_i=0$. This means that $\spec\,\cS$ can be read off
immediately by ignoring the non-zero superscripts $d_i$ of the remaining entries in the
intersection type. Following each table is a summary of relevant information obtained from the
algorithms concerning the candidate top and bottom groups in each case.


\subsection{Surviving {\sc Lines} with $k\leq 8$}\label{subsec:outputkleq8}

Table~\ref{tab:outputkleq8} displays the output of
Algorithms~\ref{alg:step3and4},~\ref{alg:qptopgroup} and \ref{alg:gnormalbottomgroup} with the
restriction $k\leq 8$. All {\sc Lines} form the output of Algorithm~\ref{alg:gnormalbottomgroup}, and
hence are potentially $G$-normal. The information obtained concerning the candidate groups is
displayed in Table~\ref{tab:outputgnormalgrpskleq8}. Algorithm~\ref{alg:qptopgroup} also produced {\sc
Line} 2 as the unique (potentially quasiprimitive) output, with candidate top group $G\cong
G^\cC\cong\PSL(3,2)<S_7$ and bottom group $G^C$ either $A_3$ or $S_3$.

\begin{table*}[!h]
\begin{center}
\begin{footnotesize}

\begin{tabular}{|cr@{\,$\cdot$\,}lr@{\,,\,}lr@{\,,\,}lr@{\,$\cdot$\,}lr@{\,$\cdot$\,}lcc|}
\hline
{\sc Line} & $d$ & $c$ & $(x$ & $y)$ & $(\gamma$ & $\delta)$ & $k^{(v)}$ & $k^{(r)}$ & $b^{(v)}$ & 
$b^{(r)}$ & int type & 
$t_{\max}$ \cr
\hline

1& 3&7 & (3&1) & (6&2) & 1&5 & 21&1 & $(1^2,3)$ &3\\
2& 7&3 & (1&3) & (2&6) & 1&5 & 21&1 & $(1^3,2)$ &2\\
3& 5&5 & (1&1) & (2&2) & 1&4 & 25&2 & $(1^2,2)$ &5\\
4& 3&19 & (9&1) & (18&2) & 1&8 & 57&1 & $(1,3,4)$ &3\\
5& 19&3 & (1&9) & (2&18) & 1&8 & 57&1 & $(1^6,2)$ &2\\
6& 9&9 & (1&1) & (2&2) & 1&5 & 81&4 & $(1^3,2)$ &2\\
7& 5&17 & (4&1) & (8&2) & 1&7 & 85&2 & $(1^2,2,3)$ &5\\
8& 17&5 & (1&4) & (2&8) & 1&7 & 85&2 & $(1^5,2)$ &1\\
9& 7&13 & (2&1) & (4&2) & 1&6 & 91&3 & $(1^2,2^2)$ &2\\
10& 13&7 & (1&2) & (2&4) & 1&6 & 91&3 & $(1^4,2)$ &2\\
12& 13&13 & (2&2) & (4&4) & 1&8 & 169&3 & $(1^4,2^2)$ &1\\
14& 9&25 & (3&1) & (6&2) & 1&8 & 225&4 & $(1^2,2^3)$ &2\\
15& 9&25 & (3&1) & (6&2) & 1&8 & 225&4 & $(1^5,3)$ &2\\
16& 25&9 & (1&3) & (2&6) & 1&8 & 225&4 & $(1^6,2)$ &2\\
35& 27&27 & (1&1) & (2&2) & 1&8 & 729&13 & $(1^6,2)$ &2\\

\hline
\end{tabular}

\caption{Surviving {\sc Lines} with $k\leq 8$\label{tab:outputkleq8}}
\end{footnotesize}
\end{center}
\end{table*}

\begin{table*}[!h]
\begin{center}
\begin{footnotesize}

\begin{tabular}{|c|l|l|}
\hline
{\sc Line} & $G^\cC$ & $G^C$ \cr
\hline
1& $A_{3}, S_{3}$ & $\PSL(3,2), C_7, D_{14}, 7:3, \AGL(1,7)$ \\
2& $\PSL(3,2), C_7, D_{14}, 7:3, \AGL(1,7)$ & $A_{3}, S_{3}$ \\
3& $A_5, S_5, D_{10}, \AGL(1,5)$ & $D_{10}, \AGL(1,5)$ \\
4& $A_{3}, S_{3}$ & $C_{19}, D_{38}, 19:3, 19:6, 19:9$ \\
5& $C_{19}, D_{38}, 19:3, 19:6, 19:9, \AGL(1,19)$ & $A_{3}, S_{3}$ \\
6& affine & affine \\
7& $A_5, S_5, D_{10}, \AGL(1,5)$ & $D_{34}, 17:4, 17:8$ \\
8& $D_{34}, 17:4, 17:8$ & $D_{10}, \AGL(1,5)$ \\
9& $\PSL(3,2), 7:3, \AGL(1,7)$ & $13:3, 13:6, \AGL(1,13)$ \\
10& $\PSL(3,3), 13:3, 13:6, \AGL(1,13)$ & $7:3, \AGL(1,7)$ \\
12& $13:3, 13:6$ & $13:3, 13:6, \AGL(1,13)$ \\
14& affine & affine \\
15& affine & affine; $(A_5\times A_5):2 \leq G^C \leq (S_5\times S_5):2$ \\
16& affine; $(A_5\times A_5):2 \leq G^\cC \leq (S_5\times S_5):2$ & affine \\
35& $3^3:13, \AGL(1,27), 3^3.13.3, \AGGL(1,27),$ & $3^3:13, \AGL(1,27), 3^3.13.3, \AGGL(1,27),$ \\
 & $\ASL(3,3), \AGL(3,3)$ & $\ASL(3,3), \AGL(3,3)$ \\ 
\hline
\end{tabular}

\caption{Candidate groups for surviving potentially $G$-normal {\sc Lines} with 
$k\leq8$\label{tab:outputgnormalgrpskleq8}}
\end{footnotesize}
\end{center}
\end{table*}


\subsection{Surviving potentially quasiprimitive {\sc Lines} with $k>8$}\label{subsec:outputqp}

Tables~\ref{tab:outputqp} and \ref{tab:outputqpgrps} display the output of Algorithm~\ref{alg:qptopgroup} with the
restriction that $k>8$. In particular, Table~\ref{tab:outputqp} displays the parameters of those {\sc Lines} in
{\sc ParameterListA}$(8)$ with $k>8$ for which {\sc QuasiprimTopGroups} was not an empty list in the output of
Algorithm~\ref{alg:qptopgroup} (including those for which $d\geq 2500$), with Table~\ref{tab:outputqpgrps} displaying
information about the candidate top and bottom groups for each surviving {\sc Line}. Recall that $G\cong G^\cC$ is almost 
simple (see Section~\ref{sect:outline}).

\begin{table*}[!h]
\begin{center}
\begin{footnotesize}

\begin{tabular}{|cccccccc|}
\hline
{\sc Line} & $d\cdot c$ & $(x,y)$ & $(\gamma,\delta)$ & $k^{(v)}\cdot k^{(r)}$ 
& $b^{(v)}\cdot b^{(r)}$ & int type & $t_{\max}$ \cr
\hline
23& $31\cdot16$ & $(2,4)$ & $(1,2)$ & $4\cdot3$ & $124\cdot15$ & $(1^8,2^2)$ &2\\
45& $133\cdot12$ & $(3,36)$ & $(1,12)$ & $6\cdot5$ & $266\cdot11$ & $(1^{24},2^3)$ &2\\
146& $528\cdot32$ & $(11,187)$ & $(1,17)$ & $22\cdot5$ & $768\cdot31$ & $(1^{88},2^{11})$ &1\\
670& $3112\cdot184$ & $(32,544)$ & $(1,17)$ & $64\cdot7$ & $8947\cdot183$ & $(1^{384},2^{32})$ &1\\
675& $3501\cdot176$ & $(36,720)$ & $(1,20)$ & $72\cdot7$ & $8558\cdot175$ & $(1^{432},2^{36})$ &1\\
707& $3510\cdot320$ & $(36,396)$ & $(1,11)$ & $72\cdot7$ & $15600\cdot319$ & $(1^{432},2^{36})$ &1\\
710& $7410\cdot240$ & $(76,2356)$ & $(1,31)$ & $152\cdot7$ & $11700\cdot239$ & $(1^{912},2^{76})$ &1\\
736& $3520\cdot783$ & $(144,648)$ & $(2,9)$ & $144\cdot7$ & $19140\cdot391$ & $(1^{720},2^{144})$ &2\\
743& $5720\cdot603$ & $(234,2223)$ & $(2,19)$ & $234\cdot7$ & $14740\cdot301$ & $(1^{1170},2^{234})$ &1\\
747& $5083\cdot848$ & $(52,312)$ & $(1,6)$ & $104\cdot7$ & $41446\cdot847$ & $(1^{624},2^{52})$ &2\\
1179& $6256\cdot696$ & $(64,576)$ & $(1,9)$ & $128\cdot7$ & $34017\cdot695$ & $(1^{768},2^{64})$ &2\\
1184& $3910\cdot1304$ & $(40,120)$ & $(1,3)$ & $80\cdot7$ & $63733\cdot1303$ & $(1^{480},2^{40})$ &2\\
1188& $9775\cdot544$ & $(100,1800)$ & $(1,18)$ & $200\cdot7$ & $26588\cdot543$ & $(1^{1200},2^{100})$ &1\\
1194& $10166\cdot536$ & $(104,1976)$ & $(1,19)$ & $208\cdot7$ & $26197\cdot535$ & $(1^{1248},2^{104})$ &1\\
1197& $3519\cdot1760$ & $(36,72)$ & $(1,2)$ & $72\cdot7$ & $86020\cdot1759$ & $(1^{432},2^{36})$ &2\\
1205& $17595\cdot464$ & $(180,6840)$ & $(1,38)$ & $360\cdot7$ & $22678\cdot463$ & $(1^{2160},2^{180})$ &1\\
1206& $3128\cdot3128$ & $(32,32)$ & $(1,1)$ & $64\cdot7$ & $152881\cdot3127$ & $(1^{384},2^{32})$ &2\\
1207& $3128\cdot3128$ & $(32,32)$ & $(1,1)$ & $64\cdot7$ & $152881\cdot3127$ & $(1^{408},2^8,3^8)$ &2\\
\hline
\end{tabular}

\caption{Surviving potentially quasiprimitive {\sc Lines} with $k>8$\label{tab:outputqp}}
\end{footnotesize}
\end{center}
\end{table*}

\begin{table*}[!h]
\begin{center}
\begin{footnotesize}

\begin{tabular}{|c|l|l|} \hline
{\sc Line} & $G^\cC$ & $G^C$ \cr \hline
23& $\PSL(3,5), \PSL(5,2)$ & $A_{16}, S_{16}$, affine \\ 
45& $\PSL(3,11)$ & $M_{11}, M_{12}, \PSL(2,11), \PGL(2,11), A_{12}, S_{12}$ \\ 
146& $A_{33}, S_{33}$ on pairs & $\PSL(2,31), \PGL(2,31), A_{32}, S_{32}$, affine \\ 
675& $\ngeq A_d$ & $HS, A_{176}, S_{176}$ \\ 
1206& $\ngeq A_d$ & no information obtained \\ 
1207& $\ngeq A_d$ & not 2-homogenenous \\ 
710,1184,1197,1205& $\ngeq A_d$ & $\PSL(2,c-1), \PGL(2,c-1), A_{c}, S_{c}$ \\ 
670,707,736,743,747, & $\ngeq A_d$ & $A_{c}, S_{c}$ \\
1179,1188,1194 &&\\ \hline
\end{tabular}

\caption{Candidate groups for surviving potentially quasiprimitive {\sc Lines} with $k>8$\label{tab:outputqpgrps}}
\end{footnotesize}
\end{center}
\end{table*}


\subsection{Surviving potentially $G$-normal {\sc Lines} with $k>8$}\label{subsec:outputgnormal}

Tables~\ref{tab:outputgnormal} and \ref{tab:outputgnormalgrps} display the output of
Algorithm~\ref{alg:gnormalbottomgroup} with the restriction that $k>8$. In particular, Table~\ref{tab:outputgnormal}
displays the parameters of those {\sc Lines} in {\sc ParameterListA}$(8)$ with $k>8$ for which {\sc GNormalBottomGroups}
was not an empty list in the output of Algorithm~\ref{alg:gnormalbottomgroup} (including those for which $c\geq 2500$),
with Table~\ref{tab:outputgnormalgrps} displaying information about the candidate top and bottom groups for each
surviving {\sc Line}.

\begin{table*}[!h]
\begin{center}
\begin{footnotesize}

\begin{tabular}{|cccccccc|}
\hline
{\sc Line} & $d\cdot c$ & $(x,y)$ & $(\gamma,\delta)$ & $k^{(v)}\cdot k^{(r)}$ 
& $b^{(v)}\cdot b^{(r)}$ & int type & $t_{\max}$ \cr
\hline
22& $16\cdot31$ & $(4,2)$ & $(2,1)$ & $4\cdot3$ & $124\cdot15$ & $(1^4,2^4)$ &2\\
46& $21\cdot81$ & $(28,7)$ & $(8,2)$ & $7\cdot5$ & $243\cdot10$ & $(2^7,3^7)$ &2\\
56& $18\cdot137$ & $(24,3)$ & $(8,1)$ & $6\cdot5$ & $411\cdot17$ & $(2^6,3^6)$ &2\\
60& $24\cdot139$ & $(18,3)$ & $(6,1)$ & $6\cdot5$ & $556\cdot23$ & $(1^{12},3^6)$ &2\\
64& $85\cdot43$ & $(5,10)$ & $(2,4)$ & $5\cdot6$ & $731\cdot21$ & $(1^{20},2^5)$ &2\\
91& $40\cdot157$ & $(60,15)$ & $(12,3)$ & $10\cdot7$ & $628\cdot13$ & $(1^{30},4^{10})$ &2\\
127& $32\cdot373$ & $(48,4)$ & $(12,1)$ & $8\cdot7$ & $1492\cdot31$ & $(1^8,3^{16})$ &2\\
128& $32\cdot373$ & $(48,4)$ & $(12,1)$ & $8\cdot7$ & $1492\cdot31$ & $(1^{24},4^8)$ &2\\
156& $64\cdot379$ & $(24,4)$ & $(6,1)$ & $8\cdot7$ & $3032\cdot63$ & $(1^8,2^{24})$ &2\\
157& $64\cdot379$ & $(24,4)$ & $(6,1)$ & $8\cdot7$ & $3032\cdot63$ & $(1^{32},3^8)$ &2\\
185& $192\cdot383$ & $(8,4)$ & $(2,1)$ & $8\cdot7$ & $9192\cdot191$ & $(1^{40},2^8)$ &2\\
673& $176\cdot3501$ & $(720,36)$ & $(20,1)$ & $72\cdot7$ & $8558\cdot175$ & $(3^{120},6^{24})$ &2\\
674& $176\cdot3501$ & $(720,36)$ & $(20,1)$ & $72\cdot7$ & $8558\cdot175$ & $(3^{144},9^8)$ &2\\
709& $240\cdot7410$ & $(2356,76)$ & $(31,1)$ & $152\cdot7$ & $11700\cdot239$ & $(4^{152},5^{76},19^4)$ &2\\
1198& $464\cdot17595$ & $(6840,180)$ & $(38,1)$ & $360\cdot7$ & $22678\cdot463$ & $(6^{360},9^{40})$ &2\\
1199& $464\cdot17595$ & $(6840,180)$ & $(38,1)$ & $360\cdot7$ & $22678\cdot463$ & 
$(3^{120},5^{144},6^{120},10^{72})$ &2\\
1200& $464\cdot17595$ & $(6840,180)$ & $(38,1)$ & $360\cdot7$ & $22678\cdot463$ & $(3^{120},6^{240},9^{80})$ 
&2\\
1201& $464\cdot17595$ & $(6840,180)$ & $(38,1)$ & $360\cdot7$ & $22678\cdot463$ & $(5^{360},10^{72})$ &2\\
1202& $464\cdot17595$ & $(6840,180)$ & $(38,1)$ & $360\cdot7$ & $22678\cdot463$ & $(5^{216},6^{120},9^{80})$ 
&2\\
1203& $464\cdot17595$ & $(6840,180)$ & $(38,1)$ & $360\cdot7$ & $22678\cdot463$ & $(3^{120},5^{216},9^{120})$ 
&2\\
1204& $464\cdot17595$ & $(6840,180)$ & $(38,1)$ & $360\cdot7$ & $22678\cdot463$ & $(5^{432},15^{24})$ &2\\
\hline
\end{tabular}

\caption{Surviving potentially $G$-normal {\sc Lines} with $k>8$\label{tab:outputgnormal}}
\end{footnotesize}
\end{center}
\end{table*}

\begin{table*}[!h]
\begin{center}
\begin{footnotesize}

\begin{tabular}{|c|l|l|}
\hline
{\sc Line} & $G^\cC$ & $G^C$ \cr\hline
22& affine $\leq \AGL(4,2)$ & $31:[15a]\leq \AGL(1,31),\,\, a\mid 2$ \\ 
46&  $A_7, S_7, \PSL(3,4), \PSSL(3,4),$ & affine $\leq \AGL(4,3)$ \\
 & $\PGL(3,4), \PGGL(3,4)$ & \\ 
56& $\PSL(2,17)$ & $137:[17a]\leq \AGL(1,137),\,\, a\mid 4$ \\ 
60& $\PSL(2,23)$ & $139:[23a]\leq \AGL(1,139),\,\, a\mid 6$ \\ 
64& $\PSL(4,4), \PSSL(4,4)$ & $43:[21a]\leq \AGL(1,43),\,\, a\mid 2$ \\ 
91& $\PSL(4,3), \PGL(4,3)$ & $157:[13a]\leq \AGL(1,157),\,\, a\mid 12$ \\ 
127,128& $\PSL(2,31), \AGL(1,32), \AGGL(1,32)$ & $373:[31a]\leq \AGL(1,373),\,\, a\mid 12$ \\ 
156,157& affine $\leq \AGL(6,2)$ & $379:[63a]\leq \AGL(1,379),\,\, a\mid 6$ \\ 
185& $\PSL(2,191)$ & $383:[191a]\leq \AGL(1,383),\,\, a\mid 2$ \\ 
673,674& $HS$ & not 3-homogeneous, $G^C_\alpha\,\,{\rm is}\,\,$\{2,3\}${\rm -group}$ \\ 
709& $\PSL(2,239)$ & not 2-homogeneous, $G^C_\alpha\,\,{\rm is}\,\,$\{2,3\}${\rm -group}$ \\ 
1202& $\PSL(2,463)$ & not 2-homogeneous \\ 
1198,1201,1204& $\PSL(2,463)$ & not 3-homogeneous \\ 
1199,1200,1203& $\PSL(2,463)$ & not 2-homogeneous, $G^C_\alpha\,\,{\rm is}\,\,$\{2,3\}${\rm -group}$ \\ 
\hline
\end{tabular}

\caption{Candidate groups for surviving potentially $G$-normal {\sc Lines} with $k>8$\label{tab:outputgnormalgrps}}
\end{footnotesize}
\end{center}
\end{table*}

\clearpage

\section{Analysing remaining {\sc Lines}}\label{sec:proofs}

In this section we consider each of the surviving cases $(\cS,G)$ appearing in the tables of
Section~\ref{sect:results} and complete the proof of Theorem~\ref{main}. Before treating the three
distinguished situations, we give two useful lemmas.

\begin{lemma}[Davies~\cite{Davies}]\label{lem:Davies} Let $g$ be a non-trivial automorphism of 
a linear space $\cS$ with constant line size $k$ and $r$ lines through a point.  Let $g$ have 
prime order $p$. Then $g$ has at most $\max(r+k-p-1,r)$ fixed points. Moreover, 
$|\Fix(h)|\leq k+r-3$ for any non-identity automorphism $h$ of such a linear space. \end{lemma}

\begin{lemma} \label{lem:KSXY} Assume that the {\sc Hypothesis} holds and that $\cC$ is $G$-normal and minimal. Let $K:= 
G_{(\cC)}$, $S:=\Soc(K)$, $X:=C_G(K)$, and $Y:=C_G(S)$. Then the following hold.
\begin{enumerate}
\item[(a)] Either (i) $Y\cap K=1$ and $S$ is nonabelian, or (ii) $Y\cap K=S$ and $S$ is elementary abelian.
\item[(b)] Either (i) $X\cap K=1$, or (ii) $X\cap K=K=S$ and $S$ is elementary abelian.
\item[(c)] Suppose in addition that $\cC$ is maximal. If there is a non-trivial intransitive subgroup 
$N\vartriangleleft G$ such that $N\cap S=1$, then there is a second $G$-normal partition $\cC'$ with 
$c$ classes of size $d$ such that $|C\cap C'|=1$ for each $C\in\cC,C'\in\cC'$.
\end{enumerate}
\end{lemma}

\begin{proof} (a) If $Y\cap K=1$ then in particular $S$ is nonabelian. So suppose $Y\cap K\ne 1$ and
note that $Y\cap K=C_K(S)$. Since $\cC$ is minimal, $G^C$ is primitive. By Theorem~\ref{thm:CP93},
$K\cong K^C$ where $C\in\cC$. Then $Y\cap K\cong(Y\cap K)^C$, a non-trivial normal subgroup that
centralises $S\cong S^C\subseteq\Soc(G^C)$. Since $G^C$ is primitive, the normal subgroups $(Y\cap
K)^C$ and $S^C$ are both transitive, and since the centraliser in $\Sym(C)$ of a transitive group is
semiregular, it follows that both $(Y\cap K)^C$ and $S^C$ are regular and are isomorphic to each
other. Moreover in this case $Y\cap K\subseteq\Soc(K)=S$, and hence $Y\cap K=S$ is elementary
abelian.

(b) Since $S\leq K$, we have $X\leq Y$. Thus if $Y\cap K=1$ then $X\cap K=1$. So we may assume by part (a) that $Y\cap
K=S$ and $S$ is elementary abelian. As in the previous paragraph, $S^C\cong S$ is self-centralising in $\Sym(C)$ and in
fact $S^C$ is a minimal normal subgroup of $G^C$, for $C\in\cC$. Since $X\cap K=C_K(K)=Z(K)$, we have $X\cap K\cong(X\cap
K)^C$ centralising $S\cong S^C$, and so $X\cap K\leq S$. By the minimality of $S^C$, either $X\cap K=1$ or $X\cap K=S$.
Since also $X\cap K=Z(K)$, in the case $X\cap K=S$ we have $S=K$.

(c) Suppose such a subgroup $N$ exists. We claim that $N\cap K=1$. If not, then $N\cap K$ is a
non-trivial normal subgroup of $K$ and so contains some minimal normal subgroup $U$ of $K$. Then
$U\leq N\cap S$, which is a contradiction. Hence $N\cap K=1$ and so $N^\cC\cong N$. Also $N^\cC$ is
normal in the primitive group $G^{\cC}$ (since $N$ is normal, and $\cC$ is maximal), and thus
$N^{\cC}$ is transitive. Let $C'$ be an $N$-orbit in $\cal P$. Then $C'$ contains an equal number of
points from each class of $\cC$ and for $C\in\cC$, $C\cap C'$ is an orbit for $N_C$ in $C$. Also,
since $N^C\vartriangleleft G^C$, $C\cap C'$ is a block of imprimitivity for the primitive group
$G^C$. Hence either $C\cap C'=C$ or $|C\cap C'|=1$. Since $N$ is intransitive it follows that
$|C\cap C'|=1$ and hence $N$ has $c$ orbits of size $d$, and so the set of $N$-orbits in $\cP$ forms
a $G$-normal partition with $c$ classes of size $d$.  \end{proof}

\subsection{{\sc Lines} with $k\leq 8$}\label{subsec:proofskleq8}

As mentioned above, line-transitive, point-imprimitive linear spaces with $k\leq 8$ have been classified by Camina and
Mischke~\cite{CaminaMischke96}. We apply this result in our proof of the following Proposition. 

\begin{proposition}\label{prop:k<9}
Suppose $k\leq 8$ so that one of {\sc Lines} $1-10,\, 12,\, 14-16,\, 35$ hold, as in 
Table~\ref{tab:outputkleq8}. Then one 
of the following holds.
\begin{enumerate}
\item[(a)] $\cS=\PG(2,q)$, $G$ lies in $Z_{q^2+q+1}.Z_3$, the normaliser of a Singer cycle, and either $q=4$ and {\sc 
Line} $1$ or $2$ holds, or $q=7$ and {\sc Line} $4$ or $5$ holds;
\item[(b)] $\cS$ is the Colbourn-McCalla design, $Z_{91}:Z_3\leq G\leq Z_{91}:Z_{12}$, and {\sc Line} $9$ or $10$ holds;
\item[(c)] $\cS$ is the Mills design, $G=Z_{91}:Z_3$, and {\sc Line} $9$ or $10$ holds;
\item[(d)] $\cS$ is one of the $467$ linear spaces constructed in 
\cite{NNOPP}, $G=\Aut(\cS)$, and {\sc Line} $35$ holds.
\end{enumerate}
In part (d) (see \cite{NNOPP},\cite{KPP93}) all the automorphism groups $\Aut(\cS)$ are of the form 
$N.Z$ with $|N|=3^6$ 
and $|Z|=13$. There are three possibilities for $N$, namely $Z_3^6,Z_9^3$, and a special $3$-group of exponent $3$; these 
correspond to $27,13$ and $427$ linear spaces respectively.
\end{proposition}

\begin{proof} First we prove that $\cC$ is $G$-normal. Suppose this is not so. Then we are in {\sc Case} 1 and $G$ is 
quasiprimitive. From Subsection~\ref{subsec:outputkleq8}, {\sc Line} 2 of Table~\ref{tab:outputkleq8} holds and the only 
possibility for $G$ is $G\cong G^\cC\cong\PSL(3,2)$. By \cite{CaminaMischke96}, $\cS={\rm 
PG}(2,4)$. However $\PSL(3,2)$ 
is not transitive on the points of ${\rm PG}(2,4)$ (see \cite{ATLAS}). Thus $\cC$ is $G$-normal.

Next we apply the classification in \cite{CaminaMischke96} to deduce that {\sc Lines} 
$3,6-8,12,14-16$ of
Table~\ref{tab:outputkleq8} do not occur. We consider the remaining {\sc Lines}. Suppose that {\sc Line} 1 or 2 of
Table~\ref{tab:outputkleq8} holds. By \cite{CaminaMischke96}, $\cS$ is a projective plane, and by 
\cite[p232]{KPP93},
$\cS={\rm PG}(2,4)$ and $G=Z_{21}$ or $Z_{21}.Z_3$ as in (a) and there are two $G$-normal partitions satisfying {\sc
Lines} 1 and 2. Next suppose that {\sc Line} 4 or 5 of Table~\ref{tab:outputkleq8} holds. By 
\cite{CaminaMischke96},
$\cS$ is a projective plane, and as there is a unique projective plane of order 7 (see
\cite{Hall},\cite{Hall-correction}), $\cS={\rm PG}(2,7)$. If {\sc Line} 4 holds then $N=\Soc(K)\cong Z_{19}$ is normal in
$G$. On the other hand if {\sc Line} 5 holds then $S=\Soc(K)\cong Z_3$ by Theorem~\ref{thm:CP93} (and see
Table~\ref{tab:outputgnormalgrpskleq8}). Then by Lemma~\ref{lem:KSXY}(a), $Y=C_G(S)$ satisfies $Y\cap K=S$. Since
$G/Y\leq\Aut(S)=Z_2$, $1\ne Y^\cC\leq G^\cC\leq\AGL(1,19)$ and so $G$ has a normal subgroup $M$ such 
that $S<M\leq Y$ and
$M/S\cong Z_{19}$. In this case the unique Sylow 19-subgroup $N$ of $M$ is normal in $G$. Thus in {\sc Line} 4 or {\sc
Line} 5 we have a normal subgroup $N\cong Z_{19}$ of $G$. Hence $G\leq N_{\PGL(2,7)}(N)=Z_{57}.Z_3$ as in (a) and there
are two $G$-normal partitions satisfying {\sc Lines} 4 and 5. Now suppose that one of {\sc Lines} 9 or 10 of
Table~\ref{tab:outputkleq8} holds.  By \cite{CaminaMischke96},\cite{ST} $\cS$ is either the 
Colbourn-McCalla design which
has automorphism group $Z_{91}:Z_{12}$ or the Mills design which has automorphism group $Z_{91}:Z_3$. Both designs have
two invariant partitions, one satisfying {\sc Line} 9 and the other satisfying {\sc Line} 10. For the Colbourn-McCalla
design, the line-transitive subgroups $G$ are all subgroups containing $Z_{91}:Z_3$. Thus (b) or (c) holds. Finally
suppose that {\sc Line} 35 holds. Then part (d) follows from \cite{NNOPP},\cite{KPP93}.

\end{proof}

\subsection{Potential quasiprimitive cases with $k>8$}\label{subsec:proofsqp}

In this section we discuss the quasiprimitive cases. In each case we have $K=G_ {(\cC)}= 1$ and $G\cong G^\cC$ is almost 
simple.

\begin{proposition}\label{prop:qp} There are no line-transitive point-imprimitive linear spaces with a quasi\-primitive
action corresponding to any of the {\sc Lines} of Table~\ref{tab:outputqp}, i.e., {\sc Lines} $23$, $45$, $146$, $670$,
$675$, $707$, $710$, $736$, $743$, $747$, $1179$, $1184$, $1188$, $1194$, $1197$, $1205$, $1206$, $1207$ cannot occur.
\end{proposition}

\begin{proof} Suppose one of the {\sc Lines} of Table~\ref{tab:outputqp} holds. We first deal with some special cases,
before giving a general argument for the remaining cases.

{\sc Line} 23: Here $G^\cC$ is either $\PSL(3,5)$ or $\PSL(5,2)$, and $G^C$ contains either $Z_2^4$ or $A_{16}$. Suppose
that $G^\cC=\PSL(3,5)$.  Then $G_C\cong G^{\cC}_C=5^2:\GL(2,5)$. However, no quotient of this group has a normal subgroup
$Z_2^4$ or $A_{16}$. Hence $G^\cC=\PSL(5,2)$, and this possibility will be considered below in the general argument.

{\sc Line} 675: Here $G^\cC$ does not contain $A_{3501}$, and $G^C$ is either $HS$, $A_{176}$ or $S_{176}$. Suppose that
$G^C\cong HS$. Now $G^{\cC}$ is a primitive group of degree $9 \times 389$, so let $p=389$ and $g$ be an element in
$G^{\cC}$ of order $p$. Suppose $g$ has some fixed points in $\cC$. Then it has at least 389 fixed points and $q$ cycles
of length $p$, where $q\leq 8$. By~\cite[Theorem 13.10]{Wielandt64}, since $G^{\cC}$ is not alternating or symmetric, the
number of fixed points should be at most $4q-4$. Thus $g$ fixes no element of $\cC$, and the almost simple group 
$G^{\cC}$ satisfies the
hypothesis (*) in ~\cite{LiebeckSaxl}. However, by checking the list in~\cite[Theorem 1.1(iii), Table 3]{LiebeckSaxl}, we find no
primitive group of degree $9\times 389$. Hence $G^C$ may only be $A_{176}$ or $S_{176}$, and these possibilities will be
considered below in the general argument. 

{\sc Lines} 1206, 1207: Here $\gamma=\delta=1$, and by Lemma~\ref{lem:subdegree} it follows that $G^\cC$ and $G^C$ are
both 2-transitive of degree 3128. Since 3128 is not a prime power, these groups are almost simple. Now 3127 is not a
prime power, and 3128 is not of the form $2^{n-1}(2^{n}\pm1)$ or $(q^n-1)/(q-1)$ for a prime power $q$, and hence the
only 2-transitive groups of degree 3128 are $A_{3128}$ and $S_{3128}$ (see \cite{Cam}). However this contradicts the fact
that $t_{\max}=2$. Thus {\sc Lines} 1206,1207 cannot occur.

We will now treat all remaining quasiprimitive cases in a somewhat uniform manner. The approach is to construct an
induced linear space on the fixed points of a certain subgroup, and then determine contradictions in the resulting
parameters, hence ruling out the remaining possibilities in Table~\ref{tab:outputqp}.

With reference to Table~\ref{tab:someqp}, we have for each {\sc Line} a prime $p$ which divides $|G|$, but which does not
divide the number of lines $b$. Thus, in each case, a Sylow $p$-subgroup $P$ of $G$ will fix some line, $\lambda$ say.
Let $F:=\Fix_\cP(P)$. For the line size $k$, define $k'$ to be the integer such that $0\leq k'<k$ and $k\equiv
k'\pmod{p}$, and similarly define $d_i'$ for each $i\in\spec\,\cS$.

Since $P$ fixes $\lambda$ setwise, $|F\cap\lambda|\geq k'$, and these values are displayed in Table~\ref{tab:someqp}. If
$k'=0$ or $1$, then we consider the sizes of the intersections of the classes with $\lambda$. Since $P$ fixes $\lambda$,
$P$ must preserve the intersection type, and so fixes the set $\cC_i$ of $d_i$ classes $C$ such that $|\lambda\cap C|=i$,
for each $i\in\spec\,\cS$. Within each $\cC_i$, $P$ must fix setwise at least $d_i'$ classes, and within each of these
fixed classes, $P$ must fix $i$ points, since $i<p$ for each {\sc Line}. So
$|F\cap\lambda|\geq\sum_{i\in\spec\,\cS}id_i'$, and these values are displayed in Table~\ref{tab:someqp} where
necessary.  For each {\sc Line} we have that $|F\cap\lambda|\geq 2$. Since for each {\sc Line} $v-k\not\equiv 0\pmod{p}$,
there exist points not on $\lambda$ which are fixed by $P$, and so $F\nsubseteq\lambda$.

Thus we may apply Corollary~\ref{cor:CaminaSiemons} to obtain, for each case, an induced linear space $\cS|_F=(F,\cL|_F)$
upon which $N_G(P)$ acts line-transitively. Then $\cS|_F$ has $v_0=|F|$ points, and the lines have size
$k_0=|F\cap\lambda|$. The lower bound for $k_0$ determined previously can be used to determine a lower bound for $v_0$,
since by considering the lines through a point $\alpha_0$ not on some line $\lambda_0$, we see that $v_0-1\geq
k_0(k_0-1)$. However by Lemma~\ref{lem:Davies}, we have that $v_0=|F|\leq k+r-3$. These inequalities, whose values
are displayed in Table~\ref{tab:someqp} for each {\sc Line}, lead to a contradiction in all but the following cases: {\sc
Line} 23 with $G^\cC=\PSL(5,2)$, {\sc Line} 710 with $G^C\geq \PSL(2,239)$, {\sc Line} 1184 with $G^C\geq \PSL(2,1303)$,
{\sc Line} 1205 with $G^C\geq \PSL(2,463)$.

\begin{table*}[!h]
\begin{center}
\begin{footnotesize}

\begin{tabular}{|c|l|c|c|c|c|c|c|}
\hline
 & & & \multicolumn{2}{c|}{$k_0\geq$} & $v-k$ && \\
{\sc Line} & Group & $p$& $k'$ & $\sum id_i'$ & mod $p$ & $k_0(k_0-1)+1$ & $k+r-3$ \\
\hline
23 &\mbox{$G^{\cC}=\PSL(5,2)$}&7&5&&1&21&54\\
45 &\mbox{$G^{\cC}=\PSL(3,11)$}&5&0&10&1&91&82\\
146 &\mbox{$G^{\cC}=A_{33},S_{33}$ on pairs}&29&23&&24&507&262\\
670 &\mbox{$G^C\geq A_{184}$}&181&86&&19&7311&1726\\
675 &\mbox{$G^C\geq A_{176}$}&173&158&&138&24807&1726\\
707 &\mbox{$G^C\geq A_{320}$}&317&187&&199&34783&2734\\
710 &\mbox{$G^C\geq A_{240}$} &233&132&&12&17293&2734\\
 &\mbox{$G^C\geq\PSL(2,239)$}&17&10&&3&91&2734\\
736 &\mbox{$G^C\geq A_{783}$}&773&235&&180&54991&3742\\
743 &\mbox{$G^C\geq A_{603}$}&601&436&&186&189661&3742\\
747 &\mbox{$G^C\geq A_{848}$}&839&728&&552&529257&6654\\
1179 &\mbox{$G^C\geq A_{696}$}&691&205&&671&41821&5758\\
1184& \mbox{$G^C\geq A_{1304}$}&1301&560&&762&313041&9678\\
 &\mbox{$G^C\geq \PSL(2,1303)$}&31&2&&6&6&9678\\
1188 &\mbox{$G^C\geq A_{544}$}&541&318&&334&100807&5198\\
1194 &\mbox{$G^C\geq A_{536}$}&523&410&&475&167691&5198\\
1197 &\mbox{$G^C\geq A_{1760}$}&1753&504&&1340&253513&12814\\
 &\mbox{$G^C\geq \PSL(2,1759)$}&293&211&&88&44311&12814\\
1205 &\mbox{$G^C\geq A_{464}$}&461&215&&16&46011&5758\\
 &\mbox{$G^C\geq \PSL(2,463)$}&7&0&14&1&183&5758\\
\hline
\end{tabular}

\caption{Relevant information for some potentially quasiprimitive cases\label{tab:someqp}}
\end{footnotesize}
\end{center}
\end{table*}

For {\sc Line} 23 ($G^\cC=\PSL(5,2)$), we have $k_0=5$ and $21\leq v_0\leq 54$ (else
$k_0(k_0-1)+1>v_0$ or $v_0>k+r-3$).  Now $v_0\equiv v\pmod{p}$ and also $(k_0-1)\mid (v_0-1)$. So
$v_0\equiv 6\pmod{7}$ and $v_0-1\equiv 0\pmod{4}$, and thus $v_0=41$. Now $N_G(P)$ is transitive on
$F$, and so $v_0\mid |N_G(P)|$. However $41\nmid |G|=|G^\cC|$, and we have a contradiction.

For {\sc Line} 710 ($G^C\geq\PSL(2,239)$), {\sc Line} 1184 ($G^C\geq\PSL(2,1303)$), and {\sc Line}
1205 ($G^C\geq\PSL(2,463)$) we have $d\not\equiv 0\pmod{p}$, and so in each case $P$ fixes some
class setwise, $C$ say. So $P\leq G_C$, and by Corollary~\ref{cor:CaminaSiemons},
$|F|=|\Fix_\cC(P)|\cdot |\Fix_C(P)|$ with $|\Fix_C(P)|\geq 3$. Now $G^C\geq\PSL(2,c-1)$, and in each 
{\sc
Line}, $c-1$ is prime, so $P^C<G^C\leq\PGL(2,c-1)$. Since $\PGL(2,c-1)$ is sharply 3-transitive, the
nontrivial subgroup $P^C$ can fix at most 2 points of $C$, contradicting $|\Fix_C(P)|\geq 3$. 

\end{proof}

\subsection{Potential $G$-normal cases with $k>8$}\label{subsec:proofsgnormal}

In this section we discuss the $G$-normal {\sc Lines}. Set $K:=G_{(\cC)}$ so that $K\ne 1$.

\begin{proposition}\label{prop:Gnormal} There are no line-transitive point-imprimitive linear spaces with an action
preserving a $G$-normal partition corresponding to any of the {\sc Lines} of Table~\ref{tab:outputgnormal}, i.e., {\sc 
Lines}
$22$, $46$, $56$, $60$, $64$, $91$, $127$, $128$, $156$, $157$, $185$, $673$, $674$, $709$, $1198-1204$ cannot occur.
\end{proposition}

\begin{proof} We give a general group theoretic argument for the {\sc Lines} of Table~\ref{tab:outputgnormal} where $c$
is a prime (with a slight variation for {\sc Line} 46, where $c$ is a prime power), to construct a non-trivial
intransitive normal subgroup which gives rise to a second $G$-normal partition with $c$ classes of size $d$. In each case
there is no corresponding {\sc Line} in Table~\ref{tab:outputgnormal}, and hence these {\sc Lines} are ruled out by
Lemma~\ref{lem:KSXY}(c). The remaining {\sc Lines} of Table~\ref{tab:outputgnormal} are then also dealt with in a
somewhat uniform manner.

Let $K,S,X,Y$ be defined as in Lemma~\ref{lem:KSXY}. Suppose one of {\sc Lines} $22$, $56$, $60$,
$64$, $91$, $127$, $128$, $156$, $157$, $185$ holds. Then the class size $c$ is a prime, $1\ne
K^C\trianglelefteq G^C$ and \mbox{$Z_c:Z_u$}$\cong G^C\leq \AGL(1,c)$ with $u\mid (c-1)$ as in
Table~\ref{tab:someGnorm}. Also, by Theorem~\ref{thm:CP93}, $K=K^C$, so $S=\Soc(K)\cong Z_c$, and
$G/Y\leq \Aut(S)\cong Z_{c-1}$. For all these {\sc Lines}, $G/K\cong G^\cC$ is non-cyclic and hence
$Y\nleq K$. Thus $Y^\cC\ne 1$. Since $Y\vartriangleleft G$, the induced group $Y^\cC$ is normal in
the primitive group $G^\cC$ and hence $Y^\cC\geq \Soc(G^\cC)$. By Lemma~\ref{lem:KSXY}(a), $Y\cap
K=S$. Thus there exists a normal subgroup $M$ of $G$ such that $S<M\leq Y$ and $M/S\cong
M^\cC=\Soc(G^\cC)$. Suppose that $G^\cC$ is affine (that is, in {\sc Lines} 22, 156, 157, and
sometimes in {\sc Lines} 127,128). Then $M/S=Z_2^{c'}$ where $d=2^{c'}$. Since $M$ centralises $S$,
it follows in these cases that $M=N\times S$ with $N$ the unique Sylow 2-subgroup of $M$. Thus
$N\vartriangleleft G$, $N$ is intransitive on $\cP$ with $c$ orbits of size $d$, and $N\cap S=1$.
For the remaining cases, $M/S\cong \Soc(G^\cC)=\PSL(n,q)$ where $d=\frac{q^n-1}{q-1}$. Since the
prime $c$ does not divide the order of the Schur multiplier of $\PSL(n,q)$ in any of these cases
(see \cite[page 302]{Gorenstein}), it follows that $M=N\times S$ with $N\cong \Soc(G^\cC)$,
$N\vartriangleleft G$, and $N\cap S=1$. Moreover $c\nmid |N|$, so $N$ is intransitive with $c$
orbits of size $d$. Thus for each {\sc Line}, by Lemma~\ref{lem:KSXY}(c), there is a $G$-normal
partition with $c$ classes of size $d$. This is a contradiction since there are no such {\sc Lines}
in Table~\ref{tab:outputgnormal}.

\begin{table*}[!h]
\begin{center}
\begin{footnotesize}

\begin{tabular}{|c|l|l|l|l|}
\hline
{\sc Line} & $G^{\cC}$ & $d$ & $G^C=c:u$ &int type \\
\hline
22 & \mbox{$\leq \AGL(4,2)$} & 16 & \mbox{31:[15a],\, $a\mid 2$} & $(1^4,2^4)$ \\
56 & \mbox{$\PSL(2,17)$} & 18 & \mbox{137:[17a],\, $a\mid 4$} & $(2^6,3^6)$ \\
60 & \mbox{$\PSL(2,23)$} & 24 & \mbox{139:[23a],\, $a\mid 6$} & $(1^{12},3^6)$ \\
64 & \mbox{$\geq\PSL(4,4)$} & 85 & \mbox{43:[21a],\, $a\mid 2$} & $(1^{20},2^5)$ \\
91 & \mbox{$\geq\PSL(4,3)$} & 40 & \mbox{157:[13a],\, $a\mid 12$} & $(1^{30},4^{10})$ \\
127,128 & \mbox{$2^5:[31a]$,\, $a\mid 5$} & 32 & \mbox{373:[31a],\, $a\mid 12$} & 
\mbox{$(1^8,3^{16})$,$(1^{24},4^8)$} \\
 & \mbox{or \PSL(2,31)} & & & \\
156,157 & \mbox{$\leq \AGL(6,2)$} & 64 & \mbox{379:[63a],\, $a\mid 6$} & \mbox{$(1^{8},2^{24})$, 
$(1^{32},3^{8})$} \\
185 & \mbox{$\PSL(2,191)$} & 192 & \mbox{383:[191a],\, $a\mid 2$} & $(1^{40},2^8)$ \\
\hline
\end{tabular}

\caption{Relevant information for some potentially $G$-normal cases\label{tab:someGnorm}}
\end{footnotesize}
\end{center}
\end{table*}

Now consider {\sc Line} 46. We follow the proof above (with $S,Y$ as there). Here $K\cong K^C\leq
\AGL(4,3)$ so $S=\Soc(K)\cong Z_3^4$ and $G/Y\leq GL(4,3)$. This means that $7\nmid |G/Y|$ and hence
$Y\nleq K$. Thus $Y^\cC\ne 1$. Arguing as in the proof above, $Y\cap K=S$ and we have a normal
subgroup $M$ of $G$ such that $S<M\leq Y$, and $M/S \cong \Soc(G^{\cC})\cong A_7$ or $\PSL(3,4)$. If
$Z(M')\ne 1$ then $Z(M')$ lies in the 3-part of the Schur multiplier of $\Soc(G^{\cC})$, and hence
$Z(M')\cong Z_3$. This means that $Z(M')<S$, and $Z(M')\vartriangleleft G$, which is a contradiction
since $G$ (even $G_C$) acts irreducibly on $S$. Thus $Z(M')=1$ and so $M=M'\times S$ with $M'$ 
intransitive, $M'\cap S=1$ and $M'\vartriangleleft G$. By Lemma~\ref{lem:KSXY}(c), there is a
$G$-normal partition with 81 classes of size 21, and this is a contradiction.

In each of the remaining {\sc Lines} $673,674,709,1198-1204$, we have a prime $p$ (displayed in
Table~\ref{tab:specialGnorm}) such that $p\mid |G^\cC|$ and $p\nmid b$. Thus a Sylow $p$-subgroup
$P$ of $G$ fixes some line, $\lambda$ say. Consider {\sc Lines} $673,674$. Here $p=3$. Let $Q:=P\cap
K$. Then $Q<G_\lambda$ and $Q$ is a Sylow 3-subgroup of $K$. Since $9\mid c$, all $P$-orbits in $C$
(for any $C\in \cC$) have length divisible by 9. However in each {\sc Line} there is a class $C$
such that $|\lambda\cap C|=3$, and this gives a contradiction. Thus one of the other {\sc Lines}
holds, and in particular $p$ divides $d-2$. Since $\Soc(G^\cC)=\PSL(2,d-1)$, it follows that $P$
fixes exactly two classes of $\cC$ setwise, say $C_1$ and $C_2$. Consider {\sc Line} 709. Here
$p=17$ and as $P\leq G_\lambda$, $P$ fixes setwise each of the four classes $C$ such that
$|\lambda\cap C|=19$. This is a contradiction. Hence $p=11$, and if the number $d_i$ of classes $C$
such that $|\lambda\cap C|=i$ is nonzero then $d_i$ must be congruent to 0,1 or 2 modulo 11 (as
otherwise $P$ would fix more than 2 classes setwise). In each of {\sc Lines} $1198-1204$ there is an
$i$ for which the condition fails.

\begin{table*}[!h]
\begin{center}
\begin{footnotesize}

\begin{tabular}{|l|l|l|c|l|}
\hline
{\sc Line} & $G^{\cC}$ & $d$ & $p$ & int type \\
\hline
673, 674  & $HS$ & 176 & 3 & \mbox{$(3^{120},6^{24})$,$(3^{144},9^8)$} \\
709 & $\PSL(2,239)$ & 240 & 17 & $(4^{152},5^{76},19^{4})$ \\
1198-1204 & $\PSL(2,463)$ & 464 & 11 & \mbox{see Table~\ref{tab:outputgnormal}} \\
\hline
\end{tabular}

\caption{Relevant information for some potentially $G$-normal cases\label{tab:specialGnorm}}
\end{footnotesize}
\end{center}
\end{table*}

\end{proof}

\subsection{Proof of Theorem~\ref{main} and concluding remarks}

Theorem~\ref{main} follows from Propositions~\ref{prop:k<9},~\ref{prop:qp} and ~\ref{prop:Gnormal}.

Our choice of $k^{(r)}_{\max}=8$ was an ambitious one, but we chose it in order to obtain at least a characterisation of 
line-transitive, point-imprimitive linear spaces that included all the known examples apart from Desarguesian projective 
planes. It turned out that a significant range of new theory was needed to complete this classification.

Disappointingly no new linear spaces were discovered in the project. However the outcomes include a
set of computational tools for searches for such linear spaces, underpinned by a broad range of
combinatorial and group theoretic results. The algorithms have been implemented as follows:
Algorithms 1-3 in C and Algorithms 4-8 in GAP 4~\cite{GAP4}. The algorithms can easily be modified
for searches over different ranges of parameters and additional group theoretic restrictions, when
available, may be added. For example, these tools are being used to extend the Camina-Mischke
classification to all line sizes up to and including line size 12.



\begin{thebibliography}{10}

\bibitem{BettenDelandtsheeretal02}
Anton Betten, Anne Delandtsheer, Alice~C. Niemeyer, and Cheryl~E. Praeger.
\newblock On a theorem of {W}ielandt for finite primitive permutation groups.
\newblock {\em J. Group Theory}, 6(4):415--420, 2003.

\bibitem{Block67}
Richard~E. Block.
\newblock On the orbits of collineation groups.
\newblock {\em Math. Z.}, 96:33--49, 1967.

\bibitem{Block68}
Richard~E. Block.
\newblock On automorphism groups of block designs.
\newblock {\em J. Combinatorial Theory}, 5:293--301, 1968.

\bibitem{Cam}
Peter~J. Cameron.
\newblock Finite permutation groups and finite simple groups.
\newblock {\em Bull. London Math. Soc.}, 13(1):1--22, 1981.

\bibitem{CaminaSiemons89a}
Alan Camina and Johannes Siemons.
\newblock Block transitive automorphism groups of {$2$}-{$(v,k,1)$} block
  designs.
\newblock {\em J. Combin. Theory Ser. A}, 51(2):268--276, 1989.

\bibitem{CaminaGagen84}
Alan~R. Camina and Terence~M. Gagen.
\newblock Block transitive automorphism groups of designs.
\newblock {\em J. Algebra}, 86(2):549--554, 1984.

\bibitem{CaminaMischke96}
Alan~R. Camina and Susanne Mischke.
\newblock Line-transitive automorphism groups of linear spaces.
\newblock {\em Electron. J. Combin.}, 3(1):Research Paper 3, approx.\ 16 pp.\
  (electronic), 1996.

\bibitem{CaminaPraeger93}
Alan~R. Camina and Cheryl~E. Praeger.
\newblock Line-transitive automorphism groups of linear spaces.
\newblock {\em Bull. London Math. Soc.}, 25(4):309--315, 1993.

\bibitem{CaminaPraeger01}
Alan~R. Camina and Cheryl~E. Praeger.
\newblock Line-transitive, point quasiprimitive automorphism groups of finite
  linear spaces are affine or almost simple.
\newblock {\em Aequationes Math.}, 61(3):221--232, 2001.

\bibitem{ATLAS}
J.~H. Conway, R.~T. Curtis, S.~P. Norton, R.~A. Parker, and R.~A. Wilson.
\newblock {\em Atlas of finite groups}.
\newblock Oxford University Press, Eynsham, 1985.
\newblock Maximal subgroups and ordinary characters for simple groups, With
  computational assistance from J. G. Thackray.

\bibitem{Davies}
D.~H. Davies.
\newblock {\em Automorphism of Designs}, PhD thesis, University of East Anglia,
  1987.

\bibitem{Delandtsheer89}
Anne Delandtsheer.
\newblock Line-primitive automorphism groups of finite linear spaces.
\newblock {\em European J. Combin.}, 10(2):161--169, 1989.

\bibitem{DelandtsheerDoyen89}
Anne Delandtsheer and Jean Doyen.
\newblock Most block-transitive {$t$}-designs are point-primitive.
\newblock {\em Geom. Dedicata}, 29(3):307--310, 1989.

\bibitem{DelandtsheerNiemeyeretal01}
Anne Delandtsheer, Alice~C. Niemeyer, and Cheryl~E. Praeger.
\newblock Finite line-transitive linear spaces: parameters and normal
  point-partitions.
\newblock {\em Adv. Geom.}, 3(4):469--485, 2003.

\bibitem{DixonMortimer}
John~D. Dixon and Brian Mortimer.
\newblock {\em Permutation groups}, volume 163 of {\em Graduate Texts in
  Mathematics}.
\newblock Springer-Verlag, New York, 1996.

\bibitem{FangLi93}
Wei~Dong Fang and Hui~Ling Li.
\newblock A generalization of the {C}amina-{G}agen theorem.
\newblock {\em J. Math. (Wuhan)}, 13(4):437--442, 1993.
\newblock (in Chinese).

\bibitem{FangLi91}
Weidong Fang.
\newblock {\em A generalization of the Camina-Gagen Theorem}, Masters thesis
  (in Chinese), Zhejiang University, China, 1991.

\bibitem{Gorenstein}
Daniel Gorenstein.
\newblock {\em Finite simple groups}.
\newblock University Series in Mathematics. Plenum Publishing Corp., New York,
  1982.
\newblock An introduction to their classification.

\bibitem{GAP4}
The~GAP Group.
\newblock {\em GAP -- Groups, Algorithms, and Programming, Version 4}, 2005,
  \verb+(http://www.gap-system.org)+.

\bibitem{Hall}
Marshall Hall, Jr.
\newblock Uniqueness of the projective plane with {$57$} points.
\newblock {\em Proc. Amer. Math. Soc.}, 4:912--916, 1953.

\bibitem{Hall-correction}
Marshall Hall, Jr.
\newblock Correction to ``{U}niqueness of the projective plane with {$57$}
  points.''.
\newblock {\em Proc. Amer. Math. Soc.}, 5:994--997, 1954.

\bibitem{HigmanMcLaughlin61}
D.~G. Higman and J.~E. McLaughlin.
\newblock Geometric {$ABA$}-groups.
\newblock {\em Illinois J. Math.}, 5:382--397, 1961.

\bibitem{LiLiu01}
Huiling Li and Weijun Liu.
\newblock Line-primitive {$2$}-{$(v,k,1)$} designs with {$k/(k,v)\le10$}.
\newblock {\em J. Combin. Theory Ser. A}, 93(1):153--167, 2001.

\bibitem{LiebeckSaxl}
Martin~W. Liebeck and Jan Saxl.
\newblock Primitive permutation groups containing an element of large prime
  order.
\newblock {\em J. London Math. Soc. (2)}, 31(2):237--249, 1985.

\bibitem{NNOPP}
Werner Nickel, Alice~C. Niemeyer, Christine~M. O'Keefe, Tim Penttila, and
  Cheryl~E. Praeger.
\newblock The block-transitive, point-imprimitive {$2$}-{$(729,8,1)$} designs.
\newblock {\em Appl. Algebra Engrg. Comm. Comput.}, 3(1):47--61, 1992.

\bibitem{KPP93}
Christine~M. O'Keefe, Tim Penttila, and Cheryl~E. Praeger.
\newblock Block-transitive, point-imprimitive designs with {$\lambda=1$}.
\newblock {\em Discrete Math.}, 115(1-3):231--244, 1993.

\bibitem{PraegerTuan02}
Cheryl~E. Praeger and Ngo~Dac Tuan.
\newblock Inequalities involving the {D}elandtsheer-{D}oyen parameters for
  finite line-transitive linear spaces.
\newblock {\em J. Combin. Theory Ser. A}, 102(1):38--62, 2003.

\bibitem{ST}
S.~D. Stoichev and V.~D. Tonchev.
\newblock The automorphism groups of the known {$2$}-{$(91,6,1)$} designs.
\newblock {\em C. R. Acad. Bulgare Sci.}, 41(4):15--16, 1988.

\bibitem{Wielandt64}
Helmut Wielandt.
\newblock {\em Finite permutation groups}.
\newblock Translated from the German by R. Bercov. Academic Press, New York,
  1964.

\end{thebibliography}
\def\cprime{$'$}


\clearpage

\appendix

\section{Appendix: {\sc ParameterList(8)}}\label{appx}

The following table summarises the output of Algorithm~\ref{alg1}.

\begin{footnotesize}

\tablefirsthead{\hline}
\tablehead{\hline\multicolumn{6}{|l|}
{\small\slshape continued from previous page}\\\hline
{\sc Line} & $v=d\cdot c$&$(x,y)$&$(\gamma,\delta)$&$k=k^{(v)}\cdot k^{(r)}$&$b^{(v)}\cdot 
b^{(r)}$\cr
\hline}
\tabletail{\multicolumn{6}{|r|}
{\small\slshape continued on next page}\\\hline}
\tablelasttail{\hline}
\begin{supertabular*}{\textwidth}{|cccccc|}
{\sc Line} & $v=d\cdot c$&$(x,y)$&$(\gamma,\delta)$&$k=k^{(v)}\cdot k^{(r)}$&$b^{(v)}\cdot b^{(r)}$ 
\cr
 \hline
1 & $21= 3 \cdot  7 $ &  $(3, 1)$ &  $(6, 2)$ & $5= 1 \cdot  5 $ &  $21\cdot  1$  
\\
2 & $21= 7 \cdot  3 $ &  $(1, 3)$ &  $(2, 6)$ & $5= 1 \cdot  5 $ &  $21\cdot  1$  
\\
3 & $25= 5 \cdot  5 $ &  $(1, 1)$ &  $(2, 2)$ & $4= 1 \cdot  4 $ &  $25\cdot  2$  
\\
4 & $57= 3 \cdot  19 $ &  $(9, 1)$ &  $(18, 2)$ & $8= 1 \cdot  8 $ &  $57\cdot  1$  
\\
5 & $57= 19 \cdot  3 $ &  $(1, 9)$ &  $(2, 18)$ & $8= 1 \cdot  8 $ &  $57\cdot  1$  
\\
6 & $81= 9 \cdot  9 $ &  $(1, 1)$ &  $(2, 2)$ & $5= 1 \cdot  5 $ &  $81\cdot  4$  
\\
7 & $85= 5 \cdot  17 $ &  $(4, 1)$ &  $(8, 2)$ & $7= 1 \cdot  7 $ &  $85\cdot  2$  
\\
8 & $85= 17 \cdot  5 $ &  $(1, 4)$ &  $(2, 8)$ & $7= 1 \cdot  7 $ &  $85\cdot  2$  
\\
9 & $91= 7 \cdot  13 $ &  $(2, 1)$ &  $(4, 2)$ & $6= 1 \cdot  6 $ &  $91\cdot  3$  
\\
10 & $91= 13 \cdot  7 $ &  $(1, 2)$ &  $(2, 4)$ & $6= 1 \cdot  6 $ &  $91\cdot  3$  
\\
11 & $136= 34 \cdot  4 $ &  $(1, 11)$ &  $(1, 11)$ & $10= 2 \cdot  5 $ &  $68\cdot  3$  
\\
12 & $169= 13 \cdot  13 $ &  $(2, 2)$ &  $(4, 4)$ & $8= 1 \cdot  8 $ &  $169\cdot  3$  
\\
13 & $196= 14 \cdot  14 $ &  $(1, 1)$ &  $(1, 1)$ & $6= 2 \cdot  3 $ &  $98\cdot  13$  
\\
14-15 & $225= 9 \cdot  25 $ &  $(3, 1)$ &  $(6, 2)$ & $8= 1 \cdot  8 $ &  $225\cdot  4$  
\\
16 & $225= 25 \cdot  9 $ &  $(1, 3)$ &  $(2, 6)$ & $8= 1 \cdot  8 $ &  $225\cdot  4$  
\\
17 & $232= 29 \cdot  8 $ &  $(2, 8)$ &  $(1, 4)$ & $12= 4 \cdot  3 $ &  $58\cdot  7$  
\\
18-19 & $441= 21 \cdot  21 $ &  $(3, 3)$ &  $(2, 2)$ & $12= 3 \cdot  4 $ &  $147\cdot  10$  
\\
20 & $456= 6 \cdot  76 $ &  $(15, 1)$ &  $(15, 1)$ & $14= 2 \cdot  7 $ &  $228\cdot  5$  
\\
21 & $456= 76 \cdot  6 $ &  $(1, 15)$ &  $(1, 15)$ & $14= 2 \cdot  7 $ &  $228\cdot  5$  
\\
22 & $496= 16 \cdot  31 $ &  $(4, 2)$ &  $(2, 1)$ & $12= 4 \cdot  3 $ &  $124\cdot  15$  
\\
23 & $496= 31 \cdot  16 $ &  $(2, 4)$ &  $(1, 2)$ & $12= 4 \cdot  3 $ &  $124\cdot  15$  
\\
24-26 & $561= 17 \cdot  33 $ &  $(6, 3)$ &  $(4, 2)$ & $15= 3 \cdot  5 $ &  $187\cdot  8$  
\\
27-30 & $561= 17 \cdot  33 $ &  $(12, 6)$ &  $(8, 4)$ & $21= 3 \cdot  7 $ &  $187\cdot  4$  
\\
31 & $561= 33 \cdot  17 $ &  $(3, 6)$ &  $(2, 4)$ & $15= 3 \cdot  5 $ &  $187\cdot  8$  
\\
32 & $561= 33 \cdot  17 $ &  $(6, 12)$ &  $(4, 8)$ & $21= 3 \cdot  7 $ &  $187\cdot  4$  
\\
33 & $638= 22 \cdot  29 $ &  $(4, 3)$ &  $(4, 3)$ & $14= 2 \cdot  7 $ &  $319\cdot  7$  
\\
34 & $638= 29 \cdot  22 $ &  $(3, 4)$ &  $(3, 4)$ & $14= 2 \cdot  7 $ &  $319\cdot  7$  
\\
35 & $729= 27 \cdot  27 $ &  $(1, 1)$ &  $(2, 2)$ & $8= 1 \cdot  8 $ &  $729\cdot  13$  
\\
36-37 & $820= 10 \cdot  82 $ &  $(9, 1)$ &  $(9, 1)$ & $14= 2 \cdot  7 $ &  $410\cdot  9$  
\\
38 & $820= 82 \cdot  10 $ &  $(1, 9)$ &  $(1, 9)$ & $14= 2 \cdot  7 $ &  $410\cdot  9$  
\\
39 & $946= 22 \cdot  43 $ &  $(2, 1)$ &  $(2, 1)$ & $10= 2 \cdot  5 $ &  $473\cdot  21$  
\\
40 & $946= 43 \cdot  22 $ &  $(1, 2)$ &  $(1, 2)$ & $10= 2 \cdot  5 $ &  $473\cdot  21$  
\\
41 & $1024= 32 \cdot  32 $ &  $(2, 2)$ &  $(1, 1)$ & $12= 4 \cdot  3 $ &  $256\cdot  31$  
\\
42 & $1548= 18 \cdot  86 $ &  $(5, 1)$ &  $(5, 1)$ & $14= 2 \cdot  7 $ &  $774\cdot  17$  
\\
43 & $1548= 86 \cdot  18 $ &  $(1, 5)$ &  $(1, 5)$ & $14= 2 \cdot  7 $ &  $774\cdot  17$  
\\
44 & $1596= 12 \cdot  133 $ &  $(36, 3)$ &  $(12, 1)$ & $30= 6 \cdot  5 $ &  $266\cdot  11$  
\\
45 & $1596= 133 \cdot  12 $ &  $(3, 36)$ &  $(1, 12)$ & $30= 6 \cdot  5 $ &  $266\cdot  11$  
\\
46 & $1701= 21 \cdot  81 $ &  $(28, 7)$ &  $(8, 2)$ & $35= 7 \cdot  5 $ &  $243\cdot  10$  
\\
47 & $1701= 81 \cdot  21 $ &  $(7, 28)$ &  $(2, 8)$ & $35= 7 \cdot  5 $ &  $243\cdot  10$  
\\
48 & $1936= 44 \cdot  44 $ &  $(1, 1)$ &  $(1, 1)$ & $10= 2 \cdot  5 $ &  $968\cdot  43$  
\\
49-50 & $2010= 15 \cdot  134 $ &  $(57, 6)$ &  $(19, 2)$ & $42= 6 \cdot  7 $ &  $335\cdot  7$  
\\
51-52 & $2010= 134 \cdot  15 $ &  $(6, 57)$ &  $(2, 19)$ & $42= 6 \cdot  7 $ &  $335\cdot  7$  
\\
53-55 & $2025= 45 \cdot  45 $ &  $(6, 6)$ &  $(4, 4)$ & $24= 3 \cdot  8 $ &  $675\cdot  11$  
\\
56 & $2466= 18 \cdot  137 $ &  $(24, 3)$ &  $(8, 1)$ & $30= 6 \cdot  5 $ &  $411\cdot  17$  
\\
57 & $2466= 137 \cdot  18 $ &  $(3, 24)$ &  $(1, 8)$ & $30= 6 \cdot  5 $ &  $411\cdot  17$  
\\
58 & $2640= 30 \cdot  88 $ &  $(3, 1)$ &  $(3, 1)$ & $14= 2 \cdot  7 $ &  $1320\cdot  29$  
\\
59 & $2640= 88 \cdot  30 $ &  $(1, 3)$ &  $(1, 3)$ & $14= 2 \cdot  7 $ &  $1320\cdot  29$  
\\
60 & $3336= 24 \cdot  139 $ &  $(18, 3)$ &  $(6, 1)$ & $30= 6 \cdot  5 $ &  $556\cdot  23$  
\\
61 & $3336= 139 \cdot  24 $ &  $(3, 18)$ &  $(1, 6)$ & $30= 6 \cdot  5 $ &  $556\cdot  23$  
\\
62-63 & $3655= 43 \cdot  85 $ &  $(10, 5)$ &  $(4, 2)$ & $30= 5 \cdot  6 $ &  $731\cdot  21$  
\\
64 & $3655= 85 \cdot  43 $ &  $(5, 10)$ &  $(2, 4)$ & $30= 5 \cdot  6 $ &  $731\cdot  21$  
\\
65-67 & $3808= 28 \cdot  136 $ &  $(40, 8)$ &  $(5, 1)$ & $48= 16 \cdot  3 $ &  $238\cdot  27$  
\\
68 & $3808= 136 \cdot  28 $ &  $(8, 40)$ &  $(1, 5)$ & $48= 16 \cdot  3 $ &  $238\cdot  27$  
\\
69-70 & $4096= 64 \cdot  64 $ &  $(12, 12)$ &  $(3, 3)$ & $40= 8 \cdot  5 $ &  $512\cdot  21$  
\\
71-72 & $4485= 39 \cdot  115 $ &  $(45, 15)$ &  $(6, 2)$ & $60= 15 \cdot  4 $ &  $299\cdot  19$  
\\
73-74 & $4485= 115 \cdot  39 $ &  $(15, 45)$ &  $(2, 6)$ & $60= 15 \cdot  4 $ &  $299\cdot  19$  
\\
75 & $4624= 68 \cdot  68 $ &  $(4, 4)$ &  $(1, 1)$ & $24= 8 \cdot  3 $ &  $578\cdot  67$  
\\
76-77 & $4761= 69 \cdot  69 $ &  $(3, 3)$ &  $(2, 2)$ & $21= 3 \cdot  7 $ &  $1587\cdot  34$  
\\
78-79 & $4761= 69 \cdot  69 $ &  $(9, 9)$ &  $(2, 2)$ & $36= 9 \cdot  4 $ &  $529\cdot  34$  
\\
80-82 & $5076= 36 \cdot  141 $ &  $(12, 3)$ &  $(4, 1)$ & $30= 6 \cdot  5 $ &  $846\cdot  35$  
\\
83 & $5076= 141 \cdot  36 $ &  $(3, 12)$ &  $(1, 4)$ & $30= 6 \cdot  5 $ &  $846\cdot  35$  
\\
84-85 & $5643= 209 \cdot  27 $ &  $(9, 72)$ &  $(2, 16)$ & $63= 9 \cdot  7 $ &  $627\cdot  13$  
\\
86-88 & $5776= 76 \cdot  76 $ &  $(20, 20)$ &  $(5, 5)$ & $56= 8 \cdot  7 $ &  $722\cdot  15$  
\\
89 & $5776= 361 \cdot  16 $ &  $(4, 96)$ &  $(1, 24)$ & $56= 8 \cdot  7 $ &  $722\cdot  15$  
\\
90-91 & $6280= 40 \cdot  157 $ &  $(60, 15)$ &  $(12, 3)$ & $70= 10 \cdot  7 $ &  $628\cdot  13$  
\\
92 & $6280= 157 \cdot  40 $ &  $(15, 60)$ &  $(3, 12)$ & $70= 10 \cdot  7 $ &  $628\cdot  13$  
\\
93 & $6816= 48 \cdot  142 $ &  $(9, 3)$ &  $(3, 1)$ & $30= 6 \cdot  5 $ &  $1136\cdot  47$  
\\
94 & $6816= 142 \cdot  48 $ &  $(3, 9)$ &  $(1, 3)$ & $30= 6 \cdot  5 $ &  $1136\cdot  47$  
\\
95-103 & $7176= 26 \cdot  276 $ &  $(33, 3)$ &  $(11, 1)$ & $42= 6 \cdot  7 $ &  $1196\cdot  25$  
\\
104 & $7176= 276 \cdot  26 $ &  $(3, 33)$ &  $(1, 11)$ & $42= 6 \cdot  7 $ &  $1196\cdot  25$  
\\
105 & $8100= 90 \cdot  90 $ &  $(1, 1)$ &  $(1, 1)$ & $14= 2 \cdot  7 $ &  $4050\cdot  89$  
\\
106 & $8836= 94 \cdot  94 $ &  $(2, 2)$ &  $(1, 1)$ & $20= 4 \cdot  5 $ &  $2209\cdot  93$  
\\
107 & $8856= 369 \cdot  24 $ &  $(4, 64)$ &  $(1, 16)$ & $56= 8 \cdot  7 $ &  $1107\cdot  23$  
\\
108 & $10296= 72 \cdot  143 $ &  $(6, 3)$ &  $(2, 1)$ & $30= 6 \cdot  5 $ &  $1716\cdot  71$  
\\
109 & $10296= 143 \cdot  72 $ &  $(3, 6)$ &  $(1, 2)$ & $30= 6 \cdot  5 $ &  $1716\cdot  71$  
\\
110-117 & $11040= 20 \cdot  552 $ &  $(174, 6)$ &  $(29, 1)$ & $84= 12 \cdot  7 $ &  $920\cdot  19$  
\\
118-119 & $11040= 96 \cdot  115 $ &  $(36, 30)$ &  $(6, 5)$ & $84= 12 \cdot  7 $ &  $920\cdot  19$  
\\
120-125 & $11040= 115 \cdot  96 $ &  $(30, 36)$ &  $(5, 6)$ & $84= 12 \cdot  7 $ &  $920\cdot  19$  
\\
126 & $11040= 552 \cdot  20 $ &  $(6, 174)$ &  $(1, 29)$ & $84= 12 \cdot  7 $ &  $920\cdot  19$  
\\
127-128 & $11936= 32 \cdot  373 $ &  $(48, 4)$ &  $(12, 1)$ & $56= 8 \cdot  7 $ &  $1492\cdot  31$  
\\
129 & $11936= 373 \cdot  32 $ &  $(4, 48)$ &  $(1, 12)$ & $56= 8 \cdot  7 $ &  $1492\cdot  31$  
\\
130-137 & $12096= 42 \cdot  288 $ &  $(42, 6)$ &  $(7, 1)$ & $60= 12 \cdot  5 $ &  $1008\cdot  41$  
\\
138-139 & $12096= 288 \cdot  42 $ &  $(6, 42)$ &  $(1, 7)$ & $60= 12 \cdot  5 $ &  $1008\cdot  41$  
\\
140 & $15400= 88 \cdot  175 $ &  $(20, 10)$ &  $(2, 1)$ & $60= 20 \cdot  3 $ &  $770\cdot  87$  
\\
141 & $15400= 175 \cdot  88 $ &  $(10, 20)$ &  $(1, 2)$ & $60= 20 \cdot  3 $ &  $770\cdot  87$  
\\
142-144 & $15841= 73 \cdot  217 $ &  $(21, 7)$ &  $(6, 2)$ & $56= 7 \cdot  8 $ &  $2263\cdot  36$  
\\
145 & $15841= 217 \cdot  73 $ &  $(7, 21)$ &  $(2, 6)$ & $56= 7 \cdot  8 $ &  $2263\cdot  36$  
\\
146 & $16896= 528 \cdot  32 $ &  $(11, 187)$ &  $(1, 17)$ & $110= 22 \cdot  5 $ &  $768\cdot  31$  
\\
147 & $18096= 48 \cdot  377 $ &  $(32, 4)$ &  $(8, 1)$ & $56= 8 \cdot  7 $ &  $2262\cdot  47$  
\\
148 & $18096= 377 \cdot  48 $ &  $(4, 32)$ &  $(1, 8)$ & $56= 8 \cdot  7 $ &  $2262\cdot  47$  
\\
149 & $19600= 140 \cdot  140 $ &  $(8, 8)$ &  $(1, 1)$ & $48= 16 \cdot  3 $ &  $1225\cdot  139$  
\\
150-153 & $19881= 141 \cdot  141 $ &  $(18, 18)$ &  $(4, 4)$ & $72= 9 \cdot  8 $ &  $2209\cdot  35$  
\\
154 & $20736= 144 \cdot  144 $ &  $(3, 3)$ &  $(1, 1)$ & $30= 6 \cdot  5 $ &  $3456\cdot  143$  
\\
155 & $24025= 155 \cdot  155 $ &  $(5, 5)$ &  $(2, 2)$ & $40= 5 \cdot  8 $ &  $4805\cdot  77$  
\\
156-157 & $24256= 64 \cdot  379 $ &  $(24, 4)$ &  $(6, 1)$ & $56= 8 \cdot  7 $ &  $3032\cdot  63$  
\\
158 & $24256= 379 \cdot  64 $ &  $(4, 24)$ &  $(1, 6)$ & $56= 8 \cdot  7 $ &  $3032\cdot  63$  
\\
159 & $25600= 160 \cdot  160 $ &  $(15, 15)$ &  $(3, 3)$ & $70= 10 \cdot  7 $ &  $2560\cdot  53$  
\\
160-161 & $31536= 54 \cdot  584 $ &  $(132, 12)$ &  $(11, 1)$ & $120= 24 \cdot  5 $ &  $1314\cdot  53$  
\\
162-163 & $31536= 584 \cdot  54 $ &  $(12, 132)$ &  $(1, 11)$ & $120= 24 \cdot  5 $ &  $1314\cdot  53$  
\\
164 & $35344= 188 \cdot  188 $ &  $(2, 2)$ &  $(1, 1)$ & $28= 4 \cdot  7 $ &  $8836\cdot  187$  
\\
165 & $36576= 96 \cdot  381 $ &  $(16, 4)$ &  $(4, 1)$ & $56= 8 \cdot  7 $ &  $4572\cdot  95$  
\\
166 & $36576= 381 \cdot  96 $ &  $(4, 16)$ &  $(1, 4)$ & $56= 8 \cdot  7 $ &  $4572\cdot  95$  
\\
167-171 & $44745= 95 \cdot  471 $ &  $(75, 15)$ &  $(10, 2)$ & $120= 15 \cdot  8 $ &  $2983\cdot  47$  
\\
172 & $44745= 471 \cdot  95 $ &  $(15, 75)$ &  $(2, 10)$ & $120= 15 \cdot  8 $ &  $2983\cdot  47$  
\\
173-174 & $48896= 128 \cdot  382 $ &  $(12, 4)$ &  $(3, 1)$ & $56= 8 \cdot  7 $ &  $6112\cdot  127$  
\\
175 & $48896= 382 \cdot  128 $ &  $(4, 12)$ &  $(1, 3)$ & $56= 8 \cdot  7 $ &  $6112\cdot  127$  
\\
176 & $58996= 172 \cdot  343 $ &  $(14, 7)$ &  $(2, 1)$ & $70= 14 \cdot  5 $ &  $4214\cdot  171$  
\\
177 & $58996= 343 \cdot  172 $ &  $(7, 14)$ &  $(1, 2)$ & $70= 14 \cdot  5 $ &  $4214\cdot  171$  
\\
178-183 & $66816= 116 \cdot  576 $ &  $(30, 6)$ &  $(5, 1)$ & $84= 12 \cdot  7 $ &  $5568\cdot  115$  
\\
184 & $66816= 576 \cdot  116 $ &  $(6, 30)$ &  $(1, 5)$ & $84= 12 \cdot  7 $ &  $5568\cdot  115$  
\\
185 & $73536= 192 \cdot  383 $ &  $(8, 4)$ &  $(2, 1)$ & $56= 8 \cdot  7 $ &  $9192\cdot  191$  
\\
186 & $73536= 383 \cdot  192 $ &  $(4, 8)$ &  $(1, 2)$ & $56= 8 \cdot  7 $ &  $9192\cdot  191$  
\\
187-188 & $86436= 294 \cdot  294 $ &  $(6, 6)$ &  $(1, 1)$ & $60= 12 \cdot  5 $ &  $7203\cdot  293$  
\\
189-190 & $117216= 198 \cdot  592 $ &  $(36, 12)$ &  $(3, 1)$ & $120= 24 \cdot  5 $ &  $4884\cdot  197$  
\\
191-192 & $117216= 592 \cdot  198 $ &  $(12, 36)$ &  $(1, 3)$ & $120= 24 \cdot  5 $ &  $4884\cdot  197$  
\\
193-196 & $117720= 135 \cdot  872 $ &  $(234, 36)$ &  $(13, 2)$ & $252= 36 \cdot  7 $ &  $3270\cdot  67$  
\\
197-198 & $117720= 872 \cdot  135 $ &  $(36, 234)$ &  $(2, 13)$ & $252= 36 \cdot  7 $ &  $3270\cdot  67$  
\\
199 & $147456= 384 \cdot  384 $ &  $(4, 4)$ &  $(1, 1)$ & $56= 8 \cdot  7 $ &  $18432\cdot  383$  
\\
200-222 & $185760= 120 \cdot  1548 $ &  $(208, 16)$ &  $(13, 1)$ & $224= 32 \cdot  7 $ &  $5805\cdot  119$  
\\
223 & $185760= 1548 \cdot  120 $ &  $(16, 208)$ &  $(1, 13)$ & $224= 32 \cdot  7 $ &  $5805\cdot  119$  
\\
224 & $197136= 444 \cdot  444 $ &  $(9, 9)$ &  $(1, 1)$ & $90= 18 \cdot  5 $ &  $10952\cdot  443$  
\\
225-226 & $229944= 429 \cdot  536 $ &  $(110, 88)$ &  $(5, 4)$ & $308= 44 \cdot  7 $ &  $5226\cdot  107$  
\\
227 & $229944= 536 \cdot  429 $ &  $(88, 110)$ &  $(4, 5)$ & $308= 44 \cdot  7 $ &  $5226\cdot  107$  
\\
228-241 & $307476= 351 \cdot  876 $ &  $(90, 36)$ &  $(5, 2)$ & $252= 36 \cdot  7 $ &  $8541\cdot  175$  
\\
242-243 & $307476= 876 \cdot  351 $ &  $(36, 90)$ &  $(2, 5)$ & $252= 36 \cdot  7 $ &  $8541\cdot  175$  
\\
244 & $309009= 249 \cdot  1241 $ &  $(255, 51)$ &  $(10, 2)$ & $357= 51 \cdot  7 $ &  $6059\cdot  124$  
\\
245-246 & $309009= 1241 \cdot  249 $ &  $(51, 255)$ &  $(2, 10)$ & $357= 51 \cdot  7 $ &  $6059\cdot  124$  
\\
247-249 & $314280= 324 \cdot  970 $ &  $(30, 10)$ &  $(3, 1)$ & $140= 20 \cdot  7 $ &  $15714\cdot  323$  
\\
250 & $314280= 970 \cdot  324 $ &  $(10, 30)$ &  $(1, 3)$ & $140= 20 \cdot  7 $ &  $15714\cdot  323$  
\\
251-360 & $362664= 207 \cdot  1752 $ &  $(612, 72)$ &  $(17, 2)$ & $504= 72 \cdot  7 $ &  $5037\cdot  103$  
\\
361-362 & $362664= 1752 \cdot  207 $ &  $(72, 612)$ &  $(2, 17)$ & $504= 72 \cdot  7 $ &  $5037\cdot  103$  
\\
363-366 & $370881= 609 \cdot  609 $ &  $(49, 49)$ &  $(2, 2)$ & $245= 49 \cdot  5 $ &  $7569\cdot  304$  
\\
367-370 & $431649= 657 \cdot  657 $ &  $(27, 27)$ &  $(2, 2)$ & $189= 27 \cdot  7 $ &  $15987\cdot  328$  
\\
371-380 & $450576= 336 \cdot  1341 $ &  $(108, 27)$ &  $(4, 1)$ & $270= 54 \cdot  5 $ &  $8344\cdot  335$  
\\
381 & $450576= 1341 \cdot  336 $ &  $(27, 108)$ &  $(1, 4)$ & $270= 54 \cdot  5 $ &  $8344\cdot  335$  
\\
382-388 & $485472= 312 \cdot  1556 $ &  $(80, 16)$ &  $(5, 1)$ & $224= 32 \cdot  7 $ &  $15171\cdot  311$  
\\
389 & $485472= 1556 \cdot  312 $ &  $(16, 80)$ &  $(1, 5)$ & $224= 32 \cdot  7 $ &  $15171\cdot  311$  
\\
390 & $544446= 522 \cdot  1043 $ &  $(42, 21)$ &  $(2, 1)$ & $210= 42 \cdot  5 $ &  $12963\cdot  521$  
\\
391 & $544446= 1043 \cdot  522 $ &  $(21, 42)$ &  $(1, 2)$ & $210= 42 \cdot  5 $ &  $12963\cdot  521$  
\\
392-669 & $572608= 184 \cdot  3112 $ &  $(544, 32)$ &  $(17, 1)$ & $448= 64 \cdot  7 $ &  $8947\cdot  183$  
\\
670-671 & $572608= 3112 \cdot  184 $ &  $(32, 544)$ &  $(1, 17)$ & $448= 64 \cdot  7 $ &  $8947\cdot  183$  
\\
672 & $602176= 776 \cdot  776 $ &  $(8, 8)$ &  $(1, 1)$ & $112= 16 \cdot  7 $ &  $37636\cdot  775$  
\\
673-674 & $616176= 176 \cdot  3501 $ &  $(720, 36)$ &  $(20, 1)$ & $504= 72 \cdot  7 $ &  $8558\cdot  175$  
\\
675 & $616176= 3501 \cdot  176 $ &  $(36, 720)$ &  $(1, 20)$ & $504= 72 \cdot  7 $ &  $8558\cdot  175$  
\\
676-677 & $632025= 795 \cdot  795 $ &  $(25, 25)$ &  $(2, 2)$ & $200= 25 \cdot  8 $ &  $25281\cdot  397$  
\\
678-679 & $681528= 584 \cdot  1167 $ &  $(24, 12)$ &  $(2, 1)$ & $168= 24 \cdot  7 $ &  $28397\cdot  583$  
\\
680 & $681528= 1167 \cdot  584 $ &  $(12, 24)$ &  $(1, 2)$ & $168= 24 \cdot  7 $ &  $28397\cdot  583$  
\\
681-682 & $799236= 894 \cdot  894 $ &  $(18, 18)$ &  $(1, 1)$ & $180= 36 \cdot  5 $ &  $22201\cdot  893$  
\\
683-706 & $1123200= 320 \cdot  3510 $ &  $(396, 36)$ &  $(11, 1)$ & $504= 72 \cdot  7 $ &  $15600\cdot  319$  
\\
707 & $1123200= 3510 \cdot  320 $ &  $(36, 396)$ &  $(1, 11)$ & $504= 72 \cdot  7 $ &  $15600\cdot  319$  
\\
708-709 & $1778400= 240 \cdot  7410 $ &  $(2356, 76)$ &  $(31, 1)$ & $1064= 152 \cdot  7 $ &  $11700\cdot  239$  
\\
710 & $1778400= 7410 \cdot  240 $ &  $(76, 2356)$ &  $(1, 31)$ & $1064= 152 \cdot  7 $ &  $11700\cdot  239$  
\\
711-712 & $2246140= 1060 \cdot  2119 $ &  $(390, 195)$ &  $(6, 3)$ & $910= 130 \cdot  7 $ &  $17278\cdot  353$  
\\
713 & $2246140= 2119 \cdot  1060 $ &  $(195, 390)$ &  $(3, 6)$ & $910= 130 \cdot  7 $ &  $17278\cdot  353$  
\\
714 & $2433600= 1560 \cdot  1560 $ &  $(16, 16)$ &  $(1, 1)$ & $224= 32 \cdot  7 $ &  $76050\cdot  1559$  
\\
715-729 & $2656900= 1630 \cdot  1630 $ &  $(150, 150)$ &  $(3, 3)$ & $700= 100 \cdot  7 $ &  $26569\cdot  543$  
\\
730-735 & $2756160= 783 \cdot  3520 $ &  $(648, 144)$ &  $(9, 2)$ & $1008= 144 \cdot  7 $ &  $19140\cdot  391$  
\\
736-737 & $2756160= 3520 \cdot  783 $ &  $(144, 648)$ &  $(2, 9)$ & $1008= 144 \cdot  7 $ &  $19140\cdot  391$  
\\
738 & $3227400= 1467 \cdot  2200 $ &  $(135, 90)$ &  $(3, 2)$ & $630= 90 \cdot  7 $ &  $35860\cdot  733$  
\\
739-740 & $3227400= 2200 \cdot  1467 $ &  $(90, 135)$ &  $(2, 3)$ & $630= 90 \cdot  7 $ &  $35860\cdot  733$  
\\
741-742 & $3449160= 603 \cdot  5720 $ &  $(2223, 234)$ &  $(19, 2)$ & $1638= 234 \cdot  7 $ &  $14740\cdot  301$  
\\
743-744 & $3449160= 5720 \cdot  603 $ &  $(234, 2223)$ &  $(2, 19)$ & $1638= 234 \cdot  7 $ &  $14740\cdot  301$  
\\
745-746 & $4310384= 848 \cdot  5083 $ &  $(312, 52)$ &  $(6, 1)$ & $728= 104 \cdot  7 $ &  $41446\cdot  847$  
\\
747 & $4310384= 5083 \cdot  848 $ &  $(52, 312)$ &  $(1, 6)$ & $728= 104 \cdot  7 $ &  $41446\cdot  847$  
\\
748-1178 & $4354176= 696 \cdot  6256 $ &  $(576, 64)$ &  $(9, 1)$ & $896= 128 \cdot  7 $ &  $34017\cdot  695$  
\\
1179-1180 & $4354176= 6256 \cdot  696 $ &  $(64, 576)$ &  $(1, 9)$ & $896= 128 \cdot  7 $ &  $34017\cdot  695$  
\\
1181-1183 & $5098640= 1304 \cdot  3910 $ &  $(120, 40)$ &  $(3, 1)$ & $560= 80 \cdot  7 $ &  $63733\cdot  1303$  
\\
1184 & $5098640= 3910 \cdot  1304 $ &  $(40, 120)$ &  $(1, 3)$ & $560= 80 \cdot  7 $ &  $63733\cdot  1303$  
\\
1185-1187 & $5317600= 544 \cdot  9775 $ &  $(1800, 100)$ &  $(18, 1)$ & $1400= 200 \cdot  7 $ &  $26588\cdot  543$  
\\
1188 & $5317600= 9775 \cdot  544 $ &  $(100, 1800)$ &  $(1, 18)$ & $1400= 200 \cdot  7 $ &  $26588\cdot  543$  
\\
1189-1193 & $5448976= 536 \cdot  10166 $ &  $(1976, 104)$ &  $(19, 1)$ & $1456= 208 \cdot  7 $ &  $26197\cdot  535$  
\\
1194 & $5448976= 10166 \cdot  536 $ &  $(104, 1976)$ &  $(1, 19)$ & $1456= 208 \cdot  7 $ &  $26197\cdot  535$  
\\
1195-1196 & $6193440= 1760 \cdot  3519 $ &  $(72, 36)$ &  $(2, 1)$ & $504= 72 \cdot  7 $ &  $86020\cdot  1759$  
\\
1197 & $6193440= 3519 \cdot  1760 $ &  $(36, 72)$ &  $(1, 2)$ & $504= 72 \cdot  7 $ &  $86020\cdot  1759$  
\\
1198-1204 & $8164080= 464 \cdot  17595 $ &  $(6840, 180)$ &  $(38, 1)$ & $2520= 360 \cdot  7 $ &  $22678\cdot  463$  
\\
1205 & $8164080= 17595 \cdot  464 $ &  $(180, 6840)$ &  $(1, 38)$ & $2520= 360 \cdot  7 $ &  $22678\cdot  463$  
\\
1206-1207 & $9784384= 3128 \cdot  3128 $ &  $(32, 32)$ &  $(1, 1)$ & $448= 64 \cdot  7 $ &  $152881\cdot  3127$  
\\
 \hline
\end{supertabular*}\\
\vfill\eject

\end{footnotesize}

{\bf Acknowledgements}

\noindent
The sixth author wishes to thank the School of Mathematics and 
Statistics at the University of Western Australia for their hospitality 
during his stay in Perth.

{\bf Authors' addresses}

\noindent
Anton Betten\,\, betten@math.colostate.edu\\
Colorado State University\\
Department of Mathematics\\
Fort Collins, CO 80523 \\
U.S.A.

\noindent
Anne Delandtsheer\,\, adelandt@ulb.ac.be\\
Universit\'e Libre de Bruxelles\\
Facult\'e des Sciences appliqu\'ees, CP165/11\\
av. F.D. Roosevelt, 50\\
B 1050 Brussels\\
Belgium

\noindent
Maska Law\,\, maska@maths.uwa.edu.au\\
Alice C. Niemeyer\,\, alice@maths.uwa.edu.au\\
Cheryl E. Praeger\,\, praeger@maths.uwa.edu.au\\
The University of Western Australia\\ 
School of Mathematics and Statistics\\ 
35 Stirling Highway, Crawley 6009\\
Australia

\noindent
Shenglin Zhou\,\, slzhou@stu.edu.cn\\
Shantou University\\
Department of Mathematics\\
Shantou, Guangdong 515063\\
People's Republic of China

\end{document}